\documentclass[a4paper,10pt]{article}
\usepackage[english]{babel}%
\usepackage[T1]{fontenc}
\usepackage[utf8]{inputenc}
\RequirePackage{fix-cm}
\usepackage{geometry}
\usepackage{titlesec}
\titleformat*{\section}{\sloppy\bfseries\Large\hyphenchar\font=-1}
\usepackage{amsmath, amsfonts, amssymb, amsthm, mathrsfs}
\usepackage{fontawesome}
\usepackage{marvosym}
\usepackage{multirow}
\usepackage[title]{appendix}
\usepackage{epsfig,graphicx,color}
\usepackage[rgb]{xcolor}%
\usepackage{framed}
\colorlet{shadecolor}{gray!2}
\usepackage{textcomp}
\usepackage{manyfoot}
\usepackage{booktabs}
\usepackage{dsfont}
\usepackage{natbib}
\usepackage{doi}
\usepackage{hyperref}%
\usepackage{varioref}
\usepackage{cleveref}
\hypersetup{
    colorlinks,
    linkcolor={brown!50!black},
    citecolor={blue!50!black},
    urlcolor={blue!80!black}
}
\usepackage{url}%
\usepackage{mathtools}
\usepackage{environ}%
\usepackage{bm}%
\usepackage{verbatim}
\usepackage{fancyhdr}
\usepackage{setspace}
\usepackage[all]{xy}
\UseCrayolaColors
\usepackage[position=top,font={sf,md,up}]{subfig}
\usepackage{orcidlink}
\newcommand{\set}[1]{\ensuremath{\left\{ #1 \right\}}}
\newcommand{\svigi}[1]{\ensuremath{\left( #1 \right)}}
\newcommand{\uN}{{\ensuremath{\zeta(N)}}}

\newcommand{\id}{\ensuremath{{\rm{d}}}}%

\newcommand{\IR}{\mathds R}
\newcommand{\IZ}{\mathbb Z}
\newcommand{\N}{\mathds N}
\newcommand{\e}{\ensuremath{\varepsilon}}%
\newcommand{\en}{\ensuremath{\varepsilon_N}}%
\newcommand{\II}{\ensuremath{\mathbb I}}%

\newcommand{\ch}{{\ensuremath{\mathds{C}}}}%
\newcommand{\dd}{{\ensuremath{\mathbb{D}}}}%
\newcommand{\IN}{\ensuremath{\mathds N}}%
\newcommand{\IB}{\ensuremath{\mathds B}}%

\newcommand{\h}{{\ensuremath{\mathds h} }}
\newcommand{\one}[1]{\ensuremath{\mathds{1}_{\left\{ #1 \right\}}}}%
\newcommand{\prob}[1]{\ensuremath{\mathds{P}\left( #1 \right)}}%
\newcommand{\EE}[1]{\ensuremath{\mathds{E}\left[ #1 \right]}}%
\newcommand{\OO}[1]{\ensuremath{{O}\left( #1 \right)}}  
\newcommand{\m}{\ensuremath{\text{\sf m}}}%
\newcommand{\ao}{\ensuremath{{\alpha_1}}}%
\newcommand{\at}{\ensuremath{{\kappa}}}%
\newcommand{\ak}{\ensuremath{\left(\frac{1}{k^\alpha} - \frac{1}{{(k+1)}^\alpha} \right)}}
\newcommand{\ako}{\ensuremath{\left(\frac{1}{k^\ao} - \frac{1}{{(k+1)}^\ao} \right)}}
\newcommand{\akt}{\ensuremath{\left(\frac{1}{k^\at} - \frac{1}{{(k+1)}^\at} \right)}}
\newcommand{\breyta}[1]{\textcolor{black}{#1}}%

\newcommand{\norm}[1]{\ensuremath{\left\lVert #1 \right\rVert } }
\theoremstyle{thmstyleone}
\newtheorem{thm}{Theorem}[section]
\newtheorem{propn}[thm]{Proposition}%
\newtheorem{lemma}[thm]{Lemma}%

\newtheorem{assumption}[thm]{Assumption}
\theoremstyle{thmstyletwo}
\newtheorem{remark}[thm]{Remark}

\theoremstyle{thmstylethree}
\newtheorem{defn}[thm]{Definition}
\makeatletter
\newcommand{\customlabel}[2]{%
   \protected@write \@auxout {}{\string \newlabel {#1}{{#2}{\theequation}{#2}{#1}{}} }%
   \hypertarget{#1}{#2}
}
\makeatother

\begin{document}

\title{\flushleft\textsf{\textbf{Gene genealogies  in haploid
populations evolving according to sweepstakes reproduction}}}
\author{{\sc Bjarki}
{\sc Eldon}\footnote{\href{mailto:beldon11@gmail.com}{beldon11@gmail.com}}\orcidlink{0000-0001-9354-2391}
}
\date{}%

\setcounter{page}{0}
\pagenumbering{arabic}
\setcounter{page}{1}
\maketitle

\renewcommand{\abstractname}{\vspace{-\baselineskip}}
\vspace{-\baselineskip}\vspace{-\baselineskip}

\abstract{ Recruitment dynamics, or the distribution of the number of
offspring among individuals, is fundamental to ecology and
evolution. We take sweepstakes reproduction to mean a skewed (heavy
right-tailed) offspring number distribution without natural selection
being involved.  Sweepstakes   may be generated by 
chance matching of  reproduction with favorable environmental conditions.
Gene genealogies generated by sweepstakes reproduction are in the
domain of attraction of  multiple-merger
coalescents where a random number of lineages merges at such times. 
We consider population genetic models of sweepstakes reproduction for
haploid panmictic populations of both  constant  ($N$),  
and  varying population size, and  evolving in a
random environment. We construct our models so that we can recover the
observed number of new mutations in a given sample without requiring
strong assumptions regarding the population size or the mutation rate.
Our main results are {\it (i)} continuous-time coalescents that are
either the Kingman coalescent or specific families of Beta- or
Poisson-Dirichlet coalescents; when combining the results the
parameter $\alpha$ of the Beta-coalescent ranges from 0 to 2, and the
Beta-coalescents may be incomplete due to an upper bound on the number
of potential offspring an arbitrary individual may produce; {\it (ii)}
in large populations we measure time in units proportional to either
$ N/\log N$ or $N$ generations; {\it (iii)} incorporating fluctuations
in population size leads to time-changed multiple-merger coalescents where the
time-change does not depend on $\alpha$; {\it (iv)} using simulations
we show that in some cases approximations of functionals 
of  a given coalescent do not match the ones of 
 the ancestral process in the domain of attraction of the given
coalescent; {\it (v)} approximations of functionals 
obtained by conditioning on the population ancestry (the ancestral
relations of all gene copies at all times)  are
broadly similar (for the models considered here) to the approximations
obtained without conditioning on the population ancestry. }

\smallskip 

{{\it Key words and phrases:}  Multiple-merger coalescent, High fecundity,
Site-frequency Spectrum,  Quenched trees, Annealed trees, Poisson-Dirichlet
distribution, Beta-coalescent}

\maketitle
\tableofcontents
\section{Introduction}\label{sec:intro}%

Inferring evolutionary histories of natural populations is one of the
main aims of population genetics. Inheritance, or the transfer of
genetic information from a parent to an offspring, is the
characteristic of organisms that makes inference possible.
Inheritance of ancestral gene copies leaves a (unknown) `trail' of
ancestral relations; in other words gene copies are connected by a
gene tree, or gene genealogy. The probabilistic modeling of the random
ancestral relations of sampled gene copies is one way to approximate
the trail going backwards in time to learn about the past events that
shaped observed genetic variation.  The mathematical study of random
gene genealogies generates a framework for constructing inference
methods \citep{Donnelly1995,W07,B09}.  Absent natural selection and complex
demography the shape or structure of gene genealogies, and so
predictions about genetic variation, will primarily be influenced by
the offspring number distribution. The importance of the offspring
number distribution for understanding evolution and ecology of natural
populations is clear.

Gene genealogies in populations where one can ignore the potential
occurrence of large (proportional to the population size) number of
offspring per individual in a sufficiently large population, are in
the domain of attraction of the Kingman coalescent \citep{K82,K82b}.
In population genetics a coalescent is a probabilistic description of
the random gene genealogies (random ancestral relations) of sampled
gene copies.  In the Kingman coalescent at most two gene copies
ancestral to the sampled ones reach a common ancestor at the same time
assuming the number of sampled gene copies is `small enough'
\citep{WT2003,BCS2014,Melfi2018b,Melfi2018}.  The classical
Wright-Fisher model (corresponding to offspring choosing a parent
independently and uniformly at random and with replacement;
~\cite{fisher22,wright31:_evolut_mendel}), the standard reference model
in population genetics, is one such population model.  However, it may
not be a good choice for highly fecund populations characterised by
sweepstakes reproduction \citep{Eldon2020,HP11}.

For understanding the ecology and evolution of highly fecund
populations the question of sweepstakes reproduction is central
\citep{HP11,Arnasonsweepstakes2022}. We say that a population evolves
according to sweepstakes reproduction \citep{H82,HP11,H94,B94} when
natural selection is not involved in mechanisms generating a skewed
offspring number distribution.  The chance matching of reproduction
with favorable environmental conditions in populations characterised
by Type III survivorship may be one such mechanism.  Sweepstakes
reproduction has also been referred to as random sweepstakes
\citep{Arnasonsweepstakes2022}.   In highly fecund
populations, in particular ones characterised by Type III
survivorship, individuals may have the capacity to produce numbers of
offspring proportional to the population size, even in a large
population.  Such events may affect the evolution of the population
should they occur often enough. Processes tracking type frequencies in
populations characterised by sweepstakes reproduction 
 are in the domain of attraction of Fleming-Viot measure-valued jump
 diffusions 
 \citep{fleming79:_some_markov,ethier93:_flemin_viot,BB09}; 
  and gene genealogies of gene copies in
samples from such populations are in the domain of attraction of    multiple-merger coalescents, where the 
 number of ancestral lineages  involved when mergers occur is random 
\citep{S03,S99,DK99,P99,S00}.  Multiple-merger coalescents also result
from recurrent strong bottlenecks \citep{BBMST09,EW06}, and strong
positive selection resulting in recurrent selective sweeps
\citep{DS05}.  Multiple-merger coalescents predict patterns of genetic
diversity different from the ones predicted by the Kingman coalescent
\citep{Blath2016,BBE2013a,SW08,EW06}, and may be essential for
explaining genetic diversity and evolutionary histories across a broad
spectrum of highly fecund natural populations
\citep{HP11,B94,A04,Eldon2020,AH2015,https://doi.org/10.1111/mec.16774,Vendrami2021,CHRISTIE2010,
Arnasonsweepstakes2022,hedgecock_etal07_marinebiology}.
This motivates us, and it is the aim here, to investigate gene
genealogies based on extensions of current models
\citep{schweinsberg03,Eldon2026biorxiv} of sweepstakes reproduction.

Population genetic models of sweepstakes reproduction for haploid and
diploid populations have been rigorously formulated and studied
\citep{S99,MS01,schweinsberg03,SW08,MS03,BLS15,HM12}.  In particular,
~\cite{EW06} and~\cite{BBE13} consider populations evolving according
to sweepstakes reproduction in a random environment. ~\cite{EW06}
consider a haploid Moran and Wright-Fisher-type models, in which, with
probability $p_N$ (where $p_{N}\to 0$ as $N\to \infty$ where $N$ is
the population size of a haploid panmictic population), a single
individual contributes a fixed fraction of the total population, and
the remaining fraction of the offspring choose a parent from the other
$N-1$ available parents independently and uniformly at random and with
replacement; with probability $1-p_N$ ordinary Wright-Fisher
reproduction occurs.  \cite{HM12} extend the Moran model to be in the
domain of attraction of multiple-merger coalescents. ~\cite{BBE13}
construct ancestral recombination graphs admitting simultaneous
mergers based on a diploid version of the model studied by
~\cite{EW06}.

Here we investigate coalescents derived from population genetic models
of sweepstakes reproduction in a haploid panmictic population evolving
in a random environment.  This article contributes {\it (i)}
continuous-time coalescents that are either the Kingman coalescent or
specific families of multiple-merger coalescents and with time in
arbitrarily large populations measured in units proportional to either
$N/\log N$ or $N$ generations; {\it (ii)} numerical comparisons of
functionals of coalescents and pre-limiting ancestral processes in the
domain of attraction of the given coalescent; {\it (iii)} comparisons
of functionals of ancestral processes when conditioning resp.\ not
conditioning on the population ancestry. Here, `population ancestry'
refers to a record of the ancestral relations of the entire population
from the present  and infinitely far into the past.

 \breyta{To briefly summarize the results, we
extend a  model of the evolution (absent selection)  of a haploid
panmictic population evolving in a constant environment   studied  by~\cite{Eldon2026biorxiv} 
to  a random environment.  In the random environment,   most of the
time (i.e.\ with high probability)   individuals
  produce `small' numbers of potential
offspring, but occasionally   the environment  turns 
favorable (there is an increased chance) for
producing a significant (at least on the order of the population size) numbers
of potential offspring.   The probability of producing  $k$ numbers of
potential offspring is assumed to  decay roughly like   $k^{-a} - {(k+1)}^{-a}$   for
$k = 1,2,\ldots, \zeta(N)$ for some given $a > 0$ and  $\zeta(N) > 0$.
The quantity $\zeta(N)$ is taken to be an  upper bound on the number of  potential
offspring  any given individual can produce. 
Depending on how $\zeta(N)/N$ behaves in a large population ($N\to
\infty$),     we obtain specific examples of continuous-time
multiple-merger coalescents (Beta$(\gamma;2-\alpha,\alpha)$- and
Poisson-Dirichlet$(\alpha,0)$-coalescents with an atom at zero).  In
contrast to   Beta- and Poisson-Dirichlet-coalescents without an
atom at zero \citep{schweinsberg03},  time in our models  is measured in
units proportional to (at least) $N/\log N$ generations, where $N$ is
the population size (in a constant-size population).  Combining the
results gives a Beta$(\gamma;2-\alpha,\alpha)$-coalescent where 
$0 < \alpha < 2$ and $0 < \gamma \le 1$ is a truncation parameter
(cf.\ also   \citep{Eldon2026biorxiv}).   We
have  extended  the range of the  $\alpha$ parameter
of the Beta$(2-\alpha,\alpha$)-coalescent  introduced by 
~\cite{schweinsberg03}.    Using simulations we show  that including an atom at zero   has
implications for the effect of population growth on gene genealogies,
in the sense that the effects do not vanish as the skewness of the
offspring distribution increases ($\alpha$ 
 decreases), in contrast to what has
been observed when applying population growth to gene genealogies
governed by Beta-coalescents without an atom at zero
\citep{freund2020cannings}.  We use simulations to compare coalescent
trees with the pre-limiting gene genealogies (using relative branch
lengths); the agreement,  especially when the probability
distribution for the  numbers of potential  offspring 
is highly skewed, 
 is poor.  We also use simulations to compare gene
genealogies obtained by conditioning resp.\ not conditioning on the
population ancestry, a record of the ancestral relations of the gene
copies in the population.  The idea is that the complete sample gene
genealogy of the gene copies in any given sample is,  for a given
population ancestry and  once the
identity of the sampled gene copies is known, fixed.     The simulation results
  indicate  that the two approaches lead to different   predictions
  about patterns of genetic variation in a sample. }

In \S~\ref{back}
we provide a brief background to coalescent processes and to models of
sweepstakes, in \S~\ref{mathresults}   we state our main results given in
Theorems~\ref{thm:haplrandomalpha} and~\ref{thm:hapl-alpha-random-one}
\breyta{for constant population size} (see \S~\ref{reshapl}), 
in Theorems~\ref{thm:convtimechangedLambdacoal}  and
~\ref{thm:time-changed-deltanullpoidr}  \breyta{(see \S~\ref{sec:lambda-coal-with})  concerning  varying  population
size}; in \S~\ref{sec:numerics} we
give numerical examples comparing functionals of stated processes,  
\S~\ref{sec:conclusion} contains  concluding remarks. Proofs of our mathematical
results may be found in
\S~\ref{sec:proofs-1}.~\breyta{  A brief discussion of the
Beta$(2-\beta,\beta)$-Poisson-Dirichlet$(\alpha,0)$-coalescent may be
found in Appendix~\ref{sec:betapoissondirichlet}.
Appendices~\ref{sec:relat-expect-sfs}--\ref{sec:further-graphs}
contain further numerical examples. }

\section{Background}
\label{back}%

For ease of reference we 
state   standard notation used throughout.
\begin{defn}[Standard notation]\label{def:notation}
Write  $\N \equiv  \set{1,2,\ldots}$, $\N_{0}
\equiv \N\cup\set{0}$. Let $N\in \N$ denote the population size (when constant).  
Asymptotics refer to the ones in an arbitrarily large population,
i.e.\ taking $N \to \infty$, unless otherwise noted.  Write
$[n] \equiv  \set{1,2,\ldots, n}$ for any $n \in \N$.  For $x\in \IR$ (the
set of reals) and $m \in \N_{0}$ recall the falling factorial
\begin{equation}
\label{eq:1}
 {(x)}_m \equiv   x(x-1)\cdots (x-m+1),\quad  {(x)}_0 \equiv  1
\end{equation}
For positive sequences ${(x_{n})}_{n}$ and ${(y_{n})}_{n}$ we write  $x_{n} \in o(y_{n})$ if
$\limsup_{n\to \infty} x_{n}/y_{n}= 0$, and $x_{n}\in O(y_{n})$ if
$\limsup_{n\to \infty} x_{n}/y_{n} < \infty $ (so that
$o(y_{n}) \subset O(y_{n})$), and $x_{n} \sim y_{n}$ if
$\lim_{n\to \infty}x_{n}/y_{n} = 1$ (assuming $y_{n} > 0$ for all
$n$).   We extend the notation  $x_{n}\sim y_{n}$ for asymptotic
equivalence  and  write
\begin{equation}
\label{simseq}
 x_{n} \overset{c}{\sim}  y_{n}
\end{equation}
if $\lim_{n\to\infty} x_{n}/y_{n} = c$ for some constant $c >0$ that
may change depending on the context (when we use this notation we
simply intend to emphasize the conditions under which~\eqref{simseq} holds
given $(x_{n})$ and $(y_{n})$); $K, C,C^{\prime}, c, c^{\prime}$ 
will denote unspecified positive constants.  Define
\begin{equation}
\label{eq:one}
\one{ E} \equiv  1
\end{equation}
if a given condition/event $E$ holds, and $\one{ E} \equiv  0$ otherwise.
The abbreviation i.i.d.\ will stand for independent and identically
distributed (random variables).
\end{defn}

The introduction of the coalescent \citep{K82,K82b,H83b,T83} advanced
the field of population genetics.  A coalescent
$\set{\xi} \equiv \{\xi(t); t\ge 0 \}$ is a Markov chain taking values
in the partitions of $\IN$, such that the restriction
$\set{\xi^{n}} \equiv \{\xi^{n}(t); t\ge 0 \} $ to $[n]$ for a fixed
$n\in \IN$ takes values in $\mathcal E_{n}$, the set of partitions of
$[n]$.  The only transitions are the merging of blocks of the current
partition,  and the time between transitions is a random
exponential with rate given by the generator of the associated
semigroup. Each block in a partition represents an ancestor to the
elements (the leaves) in each of the block, in the sense that distinct
leaves (corresponding to sampled gene copies) $i$ and $j$ are in the
same block at time $t \ge 0$ if and only if they share a common
ancestor at time $t$ in the past \citep{MS01}.  At time zero,
$\xi^{n}(0) \equiv \{ \{1\}, \ldots, \{n\}\}$, and the time
$\inf\{t \ge 0 : \xi^{n}(t) = \{ [n]\} \}$, where the partition
$\{[n]\}$ contains only the block $[n]$, is the time of the most
recent common ancestor of the $n$ sampled gene copies.  By
$\set{ \xi^{n,N}} \equiv \{ \xi^{n,N}(r); r \in \N_{0} \}$ we denote the
(pre-limiting) \emph{ancestral process}.  The ancestral process for a
sample of $n$ gene copies is a  Markov sequence (a
Markov process with a countable state space and evolving in
discrete-time) taking values in $\mathcal{E}_{n}$, where each block
represents the ancestor of the leaves (gene copies; arbitrarily
labelled) it contains, starting from $\{ \{1\}, \ldots, \{n\}\}$, and
the only transitions are the merging of blocks of the current
partition (we exclude further elements such as recombination and focus
on gene genealogies of a single  contiguous non-recombining segment of
a chromosome in a haploid  panmictic population).  We will also use
the term `coalescent' for the block-counting process associated with a
given coalescent.  
  A central quantity in proving convergence of
\set{\xi^{n,N}(\lfloor t/c_{N}\rfloor); t \ge 0} is the coalescence probability \citep{S99}.
\begin{defn}[The coalescence probability]\label{def:cNhapl}
Define  $c_{N}$ as the probability that two given   gene
copies of the same generation 
 derive  from the same parent  gene copy.
\end{defn}
One obtains a continuous-time coalescent whenever $c_{N}\to 0$ as
$N\to \infty$ with time measured in units of $\lfloor 1/c_{N} \rfloor$
generations \citep{schweinsberg03,MS01,S03}.  Let 
 $\nu_{i}$ be  the random number
of \emph{surviving} offspring of the $i$th (arbitrarily labelled)
individual in a haploid panmictic population.  When the population is
of   constant size $N$ it holds that  $\nu_{1} + \cdots + \nu_{N} =
N$.  Moreover,  since then 
$\EE{\nu_{1}} = 1$, taking $\mathds{V}(\nu_{1})$ as the variance of
$\nu_{1}$, then  by Definition~\ref{def:cNhapl}
\begin{equation}
\label{eq:6}
   c_N =  N\frac{\EE{\nu_1(\nu_1 - 1)} }{N(N-1) } = \frac{\mathds{V}(\nu_{1})}{N-1}
\end{equation}
\citep[Equation 4]{MS01}.

For a given population model, one aims to identify the limiting
process $\{\xi^{n}(t); t \ge 0\}$ to which
$\{ \xi^{n,N}( \lfloor t/c_N \rfloor ), t \ge 0 \}$ converges in
finite dimensional distributions as $N \to \infty$.
The Kingman coalescent \citep{K82,K82b}, in which each pair of blocks
in the current partition merges at rate 1, and transitions involving
more than two blocks are not possible, holds for a large class of
population models \citep{mohle1998robustness}.  A more general class
of coalescents, in which a random number of blocks merges each time,
are generally referred to as multiple-merger coalescents
\citep{P99,DK99,S99,S00}. They arise for example from population
models of sweepstakes reproduction 
\citep{HM12,schweinsberg03,SW08,EW06,BLS15,HM11b}.  Lambda-coalescents
($\Lambda$-coalescents) are multiple-merger coalescents where mergers
occur asynchronously \citep{Gnedin2014}.  They are characterised by a
finite measure $\Lambda_{+}$ on $(0,1]$ \citep{P99}.  
 In a
$\Lambda$-coalescent,  a given group of \breyta{$2  \le k \le m$} blocks merges
at the rate \citep[Equation~4]{S99}
\begin{equation}
\label{eq:19}
\lambda_{m,k} =   c\one{k=2} +  c^{\prime}\int_{0}^{1} x^{k-2}{(1-x)}^{m-k}\Lambda_{+}{\rm d}x  
\end{equation}
\citep{S99,DK99,P99}.  The Kingman coalescent is thus a special case
of a $\Lambda$-coalescent where $c =1$ and $\Lambda_{+}=0$ in
~\eqref{eq:19}. \breyta{~\cite{BB09}  discuss   the connection
between $\Lambda$-coalescents and forward-in-time Fleming-Viot
measure-valued diffusions  \citep{fleming79:_some_markov,ethier93:_flemin_viot}   tracing
the frequency of genetic types in populations evolving according to
sweepstakes reproduction.}

  We will verify Case~\ref{item:2} of  Theorem~\ref{thm:haplrandomalpha} and
Case~\ref{item:9} of  Theorem~\ref{thm:hapl-alpha-random-one} given in \S~\ref{mathresults}
by checking the conditions of \citep[Proposition~3]{schweinsberg03}
(see also \citep{S99} and \citep{MS01}).
\begin{propn}[\citep{schweinsberg03}, Proposition~3;  conditions for convergence to a $\Lambda$-coalescent]\label{prp:conditionsconvergenceLambda}
Suppose, with $c_{N}$ defined in Definition~\ref{def:cNhapl} (see 
~\eqref{eq:6})
\begin{subequations}
\begin{align}
\lim_{N\to \infty}c_{N}  & = 0, \label{eq:cntozero} \\
\lim_{N\to \infty}\frac{\EE{ {(\nu_{1})}_{2} {(\nu_{2})}_{2}} }{ N^{2}c_{N} } & = 0, \label{eq:simtozero} \\
\label{eq:existenceLambda}
\lim_{N\to \infty}\frac{N}{c_{N}} \prob{\nu_{1} > Nx} &  =  \int_{x}^{1} y^{-2}\Lambda_{+}(dy),
\end{align}
\end{subequations}
all hold where in~\eqref{eq:existenceLambda}   $\Lambda_{+}$ 
 is a finite measure  on (the Borel subsets of)  $(0,1]$ and  $0 < x < 1$ is fixed.
Then $\left\{ \xi^{n,N}( \lfloor t/c_{N} \rfloor); t \ge 0 \right\}$ converges
(in the sense of convergence of finite-dimensional distributions) to
$\left\{ \xi^{n}( t); t \ge 0 \right\}$ with transition rates as in~\eqref{eq:19}
with $\Lambda_{+} > 0$.
\end{propn}

On the other hand, suppose \breyta{ \citep[Proposition~2]{schweinsberg03},
\citep[Section~4, Equation~14]{mohle2000total}}
\begin{equation}
\label{condi}
 \begin{split}
\lim_{N\to \infty} \frac{\EE{ {(\nu_1)}_{3} }
}{N^2c_N} & = 0 \\
\end{split}
\end{equation}
Then $c_N \to 0$, and
$\{\xi^{n,N}(\lfloor t/c_{N} \rfloor); t \ge 0\}$ converges in the
sense of convergence of finite-dimensional distributions (in
\breyta{our case } 
 equivalent to weak convergence in the $J_{1}$-Skorokhod
topology) to the Kingman coalescent with transition rates
~\eqref{eq:19} where $c = 1$ and $\Lambda_{+} = 0$
\citep{mohle2000total,mohle1998robustness}.

One  particular family   of
$\Lambda$-coalescents has been much  investigated
(e.g.\ \citep{BBC05,DKW2014,BBS08,BBS07,Gnedin2014asymptotics,Birkner2024}).
\begin{defn}[The Beta$(\gamma, 2-\alpha,\alpha)$-coalescent]\label{def:betacoal}
The Beta$(\gamma,  2-\alpha,\alpha)$-coalescent is a $\Lambda$-coalescent with
transition rates  as in~\eqref{eq:19}  with $c = 0$,   and   $\Lambda_+$, 
 for  $0 < \alpha  < 2$   and $0 < \gamma  \le 1$,  is 
\begin{equation}
\label{eq:3}
\text{d} \Lambda_+(x) =  \frac{1}{B(\gamma, 2-\alpha,\alpha)} \one{0 < x \le \gamma }  x^{1-\alpha} {(1-x)}^{\alpha - 1}  \text{d}x
\end{equation}
In~\eqref{eq:3}  $B(p,   a,b) \equiv   \int_{0}^{1}\one{0 < t \le p }
t^{a-1}{(1-t)}^{b-1}dt$ for   \breyta{some given } 
$a,b > 0$ and $0 < p \le 1$ all fixed   is the (lower  incomplete when $p < 1$)   beta function.
\end{defn}
In Definition~\ref{def:betacoal} we take $0 < \alpha < 2$. The family
of Beta$(\gamma,2-\alpha,\alpha)$-coalescents with $\gamma = 1$ and
$1 \le \alpha < 2$ can be shown to follow from a particular model of
sweepstakes reproduction  \citep[Equation~11]{schweinsberg03}.
We will study extensions  \breyta{of 
 \citep[Equation~3.2]{Eldon2026biorxiv}, in turn an extension of
 \citep[Equation~11]{schweinsberg03},  such that by combining the
 results we can take $0 < \alpha < 2$.  } 
Moreover, we will consider a \emph{incomplete}
Beta$(\gamma,2-\alpha,\alpha)$-coalescent (recall
Definition~\ref{def:betacoal}) where $0 < \gamma < 1$.

Xi-coalescents ($\Xi$-coalescents) \citep{S00,S03,MS01} extend
$\Lambda$-coalescents to simultaneous mergers where ancestral lineages
may merge in two or more groups simultaneously.  Xi-coalescents arise,
for example, 
from models of  diploid   populations  evolving according to
sweepstakes  \citep{BLS15,BBE13,MS03} (see also \citep{SW08}),    recurrent strong bottlenecks
\citep{BBMST09}, and strong positive selection resulting in  recurrent selective sweeps 
 \citep{DS05,SD2005}.  Write
\begin{equation}
\label{eq:simplex}
\Delta_{+} := \set{ (x_1, x_2, \ldots ) :   x_1 \ge x_2 \ge \ldots \ge 0,
 \sum\nolimits_{j}x_{j} \le 1}  \setminus \set{(0,\ldots ) }
\end{equation}
 Let  $\Xi_+$ denote a finite measure on $\Delta_{+}$.  
Then, with $n\ge 2$ blocks in a given partition,
$k_1, \ldots, k_r \ge 2$, $r \in \N $, and
$s = n - \sum_{j=1}^r k_r \ge 0$, there exists a coalescent
$\set{\xi^{n}(t); t \ge 0}$ restricted to the partitions of $[n]$
where the rate at which $k_{1} + \cdots + k_{r} \in \{2,\ldots, n\}$
blocks merge in $r$ groups with group $j$ of size $k_{j}$ is
\breyta{($s = n - k_{1} - \cdots - k_{r}$)}
\begin{equation}
\label{xirate}
\begin{split}
\lambda_{n; k_1, \ldots, k_r;s} & =  c\one{r=1, k_1=2} \\
& \quad   +  c^{\prime} \int_{\Delta_{+}} \sum_{\ell = 0}^s \sum_{i_1 \neq \ldots \neq i_{r+\ell} } \binom{s}{\ell} x_{i_1}^{k_1} \cdots  x_{i_r}^{k_r}x_{i_{r+1}}\cdots x_{i_{r+\ell}}{\left( 1 -    \sum_{j=1}^\infty x_j  \right)}^{s - \ell} \\
& \quad  \times \frac{1}{\sum_{j=1}^\infty x_j^2 } \Xi_+ (dx) 
\end{split}
\end{equation}
\citep{S00,MS01}. Birkner et al  \citep[Theorem~A.5]{BLS15}  summarise equivalent conditions for
convergence, in the sense of finite-dimensional distributions, to a
$\Xi$-coalescent.  One condition is the existence of the limits
($k_{1},\ldots, k_{r} \in \N$ and $2 \le \min\set{k_{1}, \ldots,
k_{r}}$   for all  $r\in \N$)
\begin{equation}
\label{xmoments}
\phi_{n; k_1, \ldots, k_r} :=  \lim_{N\to \infty}  \frac{  \EE{ {(\nu_1)}_{k_1} \cdots {(\nu_r)}_{k_r}  }} {c_{N} N^{k_1 + \cdots + k_r - r } }  
\end{equation}
\citep{MS01}.
\begin{propn}[\citep{schweinsberg03} Proposition~1; convergence to a $\Xi$-coalescent]\label{prp:convergenceXicoal}
Let $r \in \N$ and $k_{1},\ldots, k_{r} \ge 2$. Suppose that the
limits in~\eqref{xmoments} exist.  Then
$\{\xi^{n,N}( \lfloor t/c_{N}\rfloor)\} $ converge (in the sense of
convergence of finite-dimensional distributions) to $\{\xi^{n}\} $
with transition rates as in~\eqref{xirate}.
\end{propn}

A discrete-time ($c_{N} \overset c \sim 1$ as $N\to \infty$)
$\Xi$-coalescent associated with the Poisson-Dirichlet distribution
with parameter $(\alpha,0)$ for $0 < \alpha < 1$ is obtained from the
same model \citep[Equation~11]{schweinsberg03} of sweepstakes
reproduction as gives rise to the Beta$(2-\alpha,\alpha)$-coalescent
\citep[Theorem~4d]{schweinsberg03}.  The Poisson-Dirichlet distribution
\citep{Kingman1975} has found wide applicability, including in
population genetics (cf.\ e.g.\
\citep{feng2010poisson,bertoin2006random}).  We will focus on the
two-parameter Poisson-Dirichlet$(\alpha, \theta)$ distribution,
 for $0 < \alpha < 1$ and
$\theta > -\alpha$ restricted to $\theta = 0$
\citep{schweinsberg03}. One can sample from the Poisson-Dirichlet$(\alpha, \theta)$
distribution in the following way. Let $ {(W_{k})}_{k \in \N}$ be a
sequence of independent random variables where $W_{k}$ has the beta
distribution with parameters $1-\alpha$ and $\theta + k\alpha$.
Define
 \begin{displaymath}
V_{1} := W_{1}, \quad V_{m} :=  (1-W_{1})\cdots (1 - W_{m-1})W_{m}, \quad   m \ge 2.
\end{displaymath}
Then $\sum_{m=1}^{\infty} V_{m} = 1 $
almost surely by construction, i.e.\ the length $V_{m}$ of every
random fragment of a stick of length 1 is added back to the sum.
The distribution of $(V_{1}, V_{2}, \ldots)$ goes by the name of the
two-paremeter GEM distribution, and will be denoted  GEM$(\alpha,\theta)$ (cf.\
\citep{feng2010poisson,ewens79,engen78:_abund,mccloskey65}).
One refers to the law of $(V_{(1)}, V_{(2)}, \ldots )$, where
$V_{(1)} \ge V_{(2)} \ge \ldots$ is ${(V_{j})}_{j\in \IN}$ ordered in
descending order, as the two-parameter Poisson-Dirichlet distribution
\citep{bertoin2006random}.  
An alternative construction of PD$(\alpha,\theta)$ for the case
$\theta = 0$ and $0 < \alpha < 1$ uses  the ranked points of a Poisson
point process on $(0,\infty)$ with characteristic measure of the form
$\nu_{\alpha}((x,\infty)) = \one{x > 0}Cx^{-\alpha}$~\cite[Section~1.3]{schweinsberg03}.

\begin{defn}[\citep{schweinsberg03}; The Poisson-Dirichlet$(\alpha,0)$-coalescent]\label{def:poisson-dirichlet}
Let $0 < \alpha < 1$ be fixed,  write
$x := (x_1, x_2, \ldots)$ for $x \in \Delta_{+}$ (recall
~\eqref{eq:simplex}), and $(x,x) := \sum_{j=1}^\infty x_j^2$. Let
$F_\alpha$ be a probability measure on $\Delta_{+}$ associated with
the Poisson-Dirichlet$(\alpha,0)$-distribution, and $\Xi_\alpha$ a
measure on $\Delta_{+}$ given by
\begin{displaymath}
\Xi_\alpha({\rm d}x) =   (x,x)F_\alpha( {\rm d}x)
\end{displaymath}
A Poisson-Dirichlet$(\alpha,0)$-coalescent is a discrete-time
$\Xi$-coalescent with $\Xi$-measure $\Xi_\alpha$ and no atom at zero.
The transition probability of merging blocks in $r$ groups of size
$k_{1}, \ldots, k_{r} \ge 2$ with current number of blocks
$b \ge k_{1} + \cdots + k_{r}$ and $s = b - k_{1} - \cdots - k_{r}$ is
(recall~\eqref{eq:1} in Definition~\ref{def:notation}),
\begin{equation}
\label{eq:34}
p_{b; k_{1}, \ldots, k_{r}; s} =  \frac{\alpha^{r+s-1}(r+s-1)! }{(b-1)!}\prod_{i=1}^{r} {(k_{i}-1 - \alpha)}_{k_{i}-1}
\end{equation}
so that $p_{2;2;0} = 1-\alpha$.  
\end{defn}
\begin{remark}[The transition probability of the Poisson-Dirichlet
coalescent]
According to \citep{schweinsberg03} the transition
probability of the Poisson-Dirichlet$(\alpha,0)$-coalescent for an $r$-merger of
$k_{i}\ge 2$ blocks in merger $i$ with $ k_{1} + \cdots + k_{r} \le b$
and $s = b - k_{1} - \cdots - k_{r}$, where $b$ is the current number
of blocks, is 
\begin{equation}
\label{eq:25}
p_{b;k_{1},\ldots, k_{r};s} =  \frac{\alpha^{r+s-1}(r+s-1)!}{(b-1)!}\prod_{i=1}^{r}{(k_{i}-\alpha)}_{k_{i}} 
\end{equation}
\citep[Equation~13]{schweinsberg03}.   
For $r = 1$ and $k_{1} = b$~\eqref{eq:25} becomes 
\begin{displaymath}
p_{b; b; 0} =  (b-\alpha)\prod_{i=1}^{b-1}\frac{b-i-\alpha}{b-i}
\end{displaymath}
In \S~\ref{sec:verifying-PDtransition-probability} we verify~\eqref{eq:34}.
\end{remark}

The Poisson-Dirichlet$(\alpha,0)$-coalescent, and the
Beta$(1,2-\alpha,\alpha)$-coalescent with $1 \le \alpha < 2$ (also
referred to as the Beta$(2-\alpha,\alpha)$-coalescent) result from the
following model of the evolution of a haploid population when the 
distribution of the number of potential offspring is as in
~\eqref{eq:4} \citep{schweinsberg03}.

\begin{defn}[Evolution of the population]\label{hschwpop}
Consider a haploid, panmictic population evolving in discrete
(non-overlapping) generations.  In any given generation each of the
current individuals independently produces a random number of
potential offspring according to some given law.  If the total number
of potential offspring produced in this way is at least some given
number  \breyta{$M$},  then $M$ of them sampled uniformly at random without
replacement survive to maturity and replace the current individuals;
the remaining potential offspring perish.  Otherwise we will assume an
unchanged population over the generation (all the potential offspring
perish before reaching maturity).
\end{defn}
The details of what happens when the total number of potential
offspring is less than the desired number can be shown to be
irrelevant in the limit by tuning the model such that the mean number
of potential offspring produced by any given individual is greater
than 1 \citep[Lemma~5]{schweinsberg03}.

Suppose a  population evolves  according to
Definition~\ref{hschwpop}. Let $X$ denote the random number of potential
offspring 
produced by an arbitrary  individual,  $\alpha, C > 0$ fixed,   and  
\begin{equation}
\label{eq:4}
\lim_{x\to \infty} x^{\alpha}\prob{X \ge x}  = C
\end{equation}
\citep[Equation~11]{schweinsberg03}.  Then
$\{\xi^{n,N}( \lfloor t/c_{N}\rfloor ); t \ge 0 \}$ converges (in the
sense of convergence of finite-dimensional distributions) to the
Kingman coalescent when $\alpha \ge 2$, to the
Beta$(\gamma, 2-\alpha,\alpha)$-coalescent as in
Definition~\ref{def:betacoal} with $\gamma = 1$ when
$1 \le \alpha < 2$, and to the
Poisson-Dirichlet$(\alpha,0)$-coalescent as in
Definition~\ref{def:poisson-dirichlet} when $0 < \alpha < 1$ but with
time measured in generations ($c_{N} \overset c\sim 1$)
\citep{schweinsberg03}.

\begin{thm}[\citep{schweinsberg03}; Theorem~4, Lemmas~6, 13, and 16]\label{thm:Schwthm4}
Recall Definition~\ref{def:notation}, and $c_{N}$ from
Definition~\ref{def:cNhapl}.  Suppose a haploid population of size $N$
evolves according to Definition~\ref{hschwpop} with law on the random
number of potential offspring as given in~\eqref{eq:4}.  Then
$\{\xi^{n,N}( \lfloor t/c_{N}\rfloor ); t \ge 0 \}$ convergences in
the sense of convergence of finite-dimensional distributions to
$\{\xi^{n}\} \equiv \set{\xi^{n}(t); t \ge 0}$ as given in each case:
\begin{enumerate} 
\item if  $ \alpha \ge 2$ then
 $\set{\xi^{n}}$ is the Kingman coalescent; 
\item\label{item:3}  if  $1 \le \alpha < 2$ then  $\set{\xi^{n}}$ is the  Beta$(\gamma,2-\alpha,\alpha)$-coalescent as in Definition~\ref{def:betacoal} with $\gamma=1$; 
\item if  $0 < \alpha < 1$ then 
$\set{\xi^{n}}$  is the Poisson-Dirichlet$(\alpha,0)$-coalescent as in
Definition~\ref{def:poisson-dirichlet}.
\end{enumerate}
By Theorem~4d, and Lemmas~13 and 16 in \citep{schweinsberg03},
 as $N\to \infty$,
\begin{equation}
\label{eq:14}
\frac{1}{c_N}  \overset{c}{\sim}  \begin{cases} 1 &  \text{when $0 < \alpha < 1$,} \\
 \log N &  \text{when $ \alpha = 1 $,} \\
 N^{\alpha - 1} &  \text{when $1 < \alpha < 2$,} \\
 N / \log N &  \text{when $\alpha = 2$,} \\
 N  &  \text{when $ \alpha > 2$} \\
\end{cases}
\end{equation}
Let $E = \{\text{given groups of $k_{1}, \ldots, k_{r} \ge 2$ blocks
merge}\}$.  By Lemma~6 in \citep{schweinsberg03}, where
$X_1, \ldots, X_N$ denotes the i.i.d.\ random number of potential
offspring produced by the current $N$ individuals according to the law
in~\eqref{eq:4}, with $S_N = X_1 + \cdots + X_N$ and
$k_{1},\ldots, k_{r} \ge 2$,
\begin{equation}
\label{eq:SchwLm6}
 \frac{1}{ c_{N}} \prob{E} = \frac{\EE{ {(\nu_{1})}_{k_{1}}\cdots (\nu_{r})_{k_{r}} }}{c_{N} N^{k_{1} + \cdots + k_{r} - r} }  \sim     \frac{N^{r}}{c_{N}} \EE{ \frac{ {(X_1)}_{k_{1}}\cdots {(X_{r})}_{k_{r}} }{S_N^{k_{1} + \cdots + k_{r}}}\one{S_N \ge N}}
\end{equation}
Moreover,   \citep[Lemma~6; Equation 17]{schweinsberg03}
\begin{displaymath}
c_N \sim  N\EE{ \frac{ {(X_1)}_2}{S_N^2}\one{S_N \ge N}}
\end{displaymath}
\end{thm}

 The law 
~\eqref{eq:4} requires that $X$ can be arbitrarily large. We will
consider extensions of~\eqref{eq:4} that relax this requirement. To
state our results 
we  require the following  example of a $\Lambda$-coalescent with rate as in
~\eqref{eq:19} with $c,c^{\prime} > 0$.    
\begin{defn}[The $\delta_{0}$-Beta$(\gamma,2-\alpha,\alpha)$-coalescent]\label{def:KingmanBeta}
The $\delta_{0}$-Beta$(\gamma,2-\alpha,\alpha)$-coalescent is a
$\Lambda$-coalescent with $\Lambda$-measure
$\Lambda = c\delta_{0} + c^{\prime}\Lambda_{+}$ taking values in the
partitions of $[n]$ with rate as in~\eqref{eq:19} with
$c,c^{\prime} > 0$ and measure $\Lambda_{+}$ as in~\eqref{eq:3} in
Definition~\ref{def:betacoal}.
\end{defn}
We also require the following extension of the
Poisson-Dirichlet$(\alpha,0)$-coalescent (recall
Definition~\ref{def:poisson-dirichlet}).
\begin{defn}[The $\delta_{0}$-Poisson-Dirichlet$(\alpha,0)$-coalescent]\label{def:kingman-poisson-dirichlet}
The $\delta_{0}$-Poisson-Dirichlet$(\alpha,0)$-coalescent is a
continuous-time  $\Xi$-coalescent with $\Xi$-measure $\Xi = \delta_{0} + \Xi_{+}$  
taking values in the partitions of $[n]$, 
where  $\Xi_{+} =  \Xi_{\alpha}$ is as in
Definition~\ref{def:poisson-dirichlet}.  Let $n \ge 2$ denote the current
number of blocks in a partition,  $k_{1},\ldots, k_{r} \ge 2$ with
$2\le k_{1} + \cdots + k_{r} \le n$ for some $r \in \N$ denoting the
merger sizes of $r$ (simultaneous) mergers, and $s = n -
k_{1} - \cdots - k_{r}$. Then the rate at which such mergers occur is
\begin{equation}
\label{eq:ratekpd}
\lambda_{n;k_{1}, \ldots, k_{r}; s} =   c\one{r=1,k_{1} = 2} +  c^{\prime} p_{n;k_{1},\ldots, k_{r};s}
\end{equation}
where  $p_{n;k_{1},\ldots, k_{r};s}$ is as in~\eqref{eq:34}. 
\end{defn}


\section{Results}\label{mathresults}

In this section we collect the results. In \S~\ref{reshapl} we state
the main mathematical results in Theorems~\ref{thm:haplrandomalpha}, 
~\ref{thm:hapl-alpha-random-one},
~\ref{thm:convtimechangedLambdacoal}, and
~\ref{thm:time-changed-deltanullpoidr}.     In \S~\ref{sec:numerics} we give
 numerical examples comparing functionals (mean relative branch
lengths) of the $\set{\xi^{n}}$ and $\set{\xi^{n,N}}$, 
 and  comparing functionals of $\set{\xi^{n,N}}$ to the
ones obtained by conditioning on the population ancestry (see
\S~\ref{sec:numerics}).


First we state the population model underlying the mathematical
results. We will adapt the formulation of  \citep{Eldon2026biorxiv}.   Suppose a population is and  evolves as in Definition~\ref{hschwpop}.
Let $X$ be the random number of potential offspring produced by an
arbitrary individual.  Let $\zeta: \N\to \N$ be a deterministic function
 and $a > 0$ be fixed.  For all  $k \in [\zeta(N)] =  \{ 1,2, \ldots, \uN \}$
let the probability mass function of the law of $X$ be bounded by
\begin{equation}
\label{eq:PXiJ}
g_{a}(k) \left( \frac{1}{k^a} - \frac{1}{ {(1+k)}^a} \right)  \leq
\prob{X = k} \leq f_{a}(k) \left( \frac{1}{k^a} -
\frac{1}{{(1+k)}^a} \right), 
\end{equation}
where   $f_{a}$ and $g_{a}$ are positive functions on
$\IN$. The quantity $\uN$ is
an upper bound on the random number of potential offspring produced by
any one  individual such that $\prob{X \le \uN} = 1$. We assign
any mass outside $[\uN]$ to $\{ X= 0 \}$.  The model in
\eqref{eq:PXiJ} is an extension of the one in \eqref{eq:4} since, if
$X$ has  law \eqref{eq:4} then
$\prob{X = x} = \prob{X \ge x} - \prob{X \ge x+1} \overset{c}{\sim}
x^{-\alpha} - (x+1)^{-\alpha}$ for $x$ large.  We will identify
conditions on $g_{a}$ and $f_{a}$ required for convergence of the
(time-rescaled) ancestral process  to a non-trivial
limit. Write 
  \begin{equation}
\label{eq:28}
X \vartriangleright  \mathds{L}(a, \uN)
\end{equation} 
when the law of $X $ is according to \eqref{eq:PXiJ} with $a$ and
$\uN$ as given each time.

We define
\begin{subequations}
\begin{align}
\notag
\underline{g_{a}}(k) := \inf_{i\ge k } g_{a}(i), \quad  \overline{g_{a}}(k) := \sup_{i\ge k } g_{a}(i), \\ \notag
\underline{g_{a}} :=\sup_{k} \inf_{i\ge k } g_{a}(i), \quad  \overline{g_{a}} := \inf_{k}\sup_{i\ge k } g_{a}(i), \\   \label{eq:isfg}
\underline{f_{a}}(k) := \inf_{i\ge k } f_{a}(i), \quad \overline{f_{a}}(k) := \sup_{i\ge k } f_{a}(i), \\ \notag
\underline{f_{a}} :=\sup_{k} \inf_{i\ge k } f_{a}(i), \quad 
\overline{f_{a}} := \inf_{k}\sup_{i\ge k } f_{a}(i), \\ \notag
g_{a}^{(\infty)}  := \lim_{k\to \infty} g_{a}(k), \quad 
f_{a}^{(\infty)}  := \lim_{k\to \infty} f_{a}(k) 
\end{align}
\end{subequations}
We assume $g_{a}(k) \le f_{a}(k)$, $\underline{g_{a}} > 0$,
$\overline{f_{a}} < \infty$, and that $f_{a}^{(\infty)},
g_{a}^{(\infty)} > 0$   \breyta{ exist}.
 Then, for all  $k\in [\uN]$,
\begin{equation}
\label{eq:12}
\left( \frac{1}{k^a} - \frac{1}{(1+\uN)^a} \right)\underline{g_{a}}(k) \le  \prob{X  \ge  k} \le \left( \frac{1}{k^a} - \frac{1}{ (\uN+1)^a} \right)\overline{f_{a}}(k)
\end{equation}
\begin{assumption}[$\EE{X} > 1$]
\label{a:mN1plus}
Suppose $X\vartriangleright \mathds L(a,\zeta(N))$ (recall
\eqref{eq:28}).  The functions $f_{a}$ and $g_{a}$ in \eqref{eq:PXiJ}
are such that $\EE{X} > 1$ for all $ N \in \mathds N$.
\end{assumption}

Assumption~\ref{a:mN1plus} results in
$\prob{X_{1} + \cdots + X_{N} < N}$ decreasing exponentially in $N$,
where $X_1, \ldots, X_N$ are independent copies of $X$
\citep[Lemma~5]{schweinsberg03}.

For $X \vartriangleright \mathds L(a,\uN)$,
using the lower bound in \eqref{eq:20} and \eqref{eq:21} in
Lemma~\ref{lm:ineqs} with
$a^{\prime} \equiv  a\one{0 < a < 1} + \one{a \ge 1}$
\begin{equation}
\label{EX}
\begin{split}
\EE{X} & \ge   \underline{g_{a}} \sum_{k=1}^{ \uN }k\left( \frac{1}{k^{a}} - \frac{1}{(k+1)^{a}} \right)  \ge  \underline{g_{a}} a^{\prime}  \sum_{k=1}^{\uN}\frac{k }{(1+k)^{1+ a}}  \\
& =   \underline{g_{a}}a^{\prime}  \sum_{k=1}^{\uN}(1+k)^{-a} -    \underline{g_{a}}a^{\prime} \sum_{k=1}^{\uN}(1+k)^{-1 - a}   \\
 & \ge  \underline{g_{a}}a^{\prime}  \int_{1}^{\uN+1 }(1+x)^{-a} dx -    \underline{g_{a}}a^{\prime} \int_{0}^{\uN }(1+x)^{-1-a} dx  \\
 & = \begin{cases}  \frac{\underline{g_{a}}a^{\prime}}{1-a} \left( (2 + \uN)^{1-a} - 2^{1-a}\right)  +  \frac{\underline{g_{a}}a^{\prime} }{a} \left( (1 + \uN)^{-a} - 1 \right) & \text{when $a \neq 1$} \\
  \underline{g_{a}}a^{\prime}( \log(\uN + 2) - \log 2 ) +  \underline{g_{a}}a^{\prime}\frac{ \uN  }{1 + \uN}  & \text{when $a = 1$}\\
\end{cases} 
\end{split}
\end{equation}

For ease of reference   we state  notation used from here on.
\begin{defn}[Notation]
\label{def:not}
Let $X_1, \ldots, X_N$ denote the independent  random number of potential offspring
(recall Definition~\ref{hschwpop}) produced in an arbitrary generation
by the   current  $N$ individuals.     Write 
\begin{subequations}
\begin{align}
S_N &  := X_1 + \cdots + X_N  \label{SN} \\
\widetilde{S}_N &  :=    X_2 + \cdots + X_N \label{SN2} \\
m_{N} &  := \EE{X_{1} } \to m_{\infty} \label{eq:5} \\
\label{eq:26}
\breyta{ \frac{ x_{n} }{y_{n} }} &  \gneqq 0
\end{align}
\end{subequations}
where $m_{\infty}$ in \eqref{eq:5}  is the limit as $N\to \infty$, and \eqref{eq:26}
 \breyta{ for positive sequences $(x_{n})$ and $(y_{n})$ will
 mean that 
either  $x_{n}/y_{n}$ converges to some  $\ell > 0$,  or 
$x_{n}/y_{n}$ diverges as $n\to \infty$}.
\end{defn}

\subsection{Beta- and Poisson-Dirichlet-coalescents }
\label{reshapl}

Before we state the results we define random environments
(Definitions~\ref{def:haplrandomalpha} and
\ref{def:alpha-random-one}); the population then evolves in the given
environment as in Definition~\ref{hschwpop}. The environments are
modelled as simple mixture distributions on the law of the number of
potential offspring (recall \eqref{eq:PXiJ}).  Depending on the
specific scenario each time, we obtain a continuous-time $\Lambda$- or
$\Xi$-coalescent with an atom at zero, and with time measured in units
proportional to either  $N/\log N$ or $N$  generations.

\begin{defn}[Type $A$ random environment]
\label{def:haplrandomalpha}
Suppose a population evolves according to Definition~\ref{hschwpop}.
Fix $0 < \alpha < 2$ and $2 \le \at $. Recall \eqref{eq:28} and the
$X_{1}, \ldots, X_{N}$ from Definition~\ref{def:not}.  Write $E $
for the event when 
$ X_i \vartriangleright \mathds L(\alpha, \uN)$ for all 
 $i\in [N]$,   and $E^{\sf c}$ for the event when $\at$ replaces $\alpha$
in $ E$ ($X_{i} \vartriangleright \mathds L(\kappa,\zeta(N))$ for all
$i\in [N]$);   $\uN$ is fixed for each $N\in \N$.   Let
$(\varepsilon_{N})_{N\in \N}$ be a positive sequence with
$0 < \varepsilon_{N} < 1$ for all $N$.  It may hold that
$\varepsilon_{N}\to 0$ as $N\to \infty$.  Suppose 
\begin{displaymath}
\prob{ E}  = \varepsilon_N, \quad \prob{ E^{\sf c}}  = 1 - \varepsilon_N
\end{displaymath}
\end{defn}%
Definition~\ref{def:haplrandomalpha} says that the population evolves
in an environment where changes to the environment affect all
individuals equally so that at any given time the
$X_{1}, \ldots, X_{N}$ are i.i.d.  The bound $\uN$ stays the same
between $ E $ and $ E^{\sf c} $.  We will identify conditions on
$\varepsilon_{N}$ (see \eqref{eq:33} in
Lemma~\ref{lm:cNhaplrandomall}) so that the ancestral process
$\set{\xi^{n,N}\svigi{\lfloor t/c_{N} \rfloor }; t \ge 0}$ converges
(in the sense of convergence of finite-dimensional distributions) 
to a non-trivial limit as $N\to \infty$. Moreover, the
$X_{1}, \ldots, X_{N}$ in Definition~\ref{def:haplrandomalpha} are
i.i.d.\ copies of $X$ where (recall \eqref{eq:one} in
Definition~\ref{def:notation})
\begin{displaymath}
X \vartriangleright \mathds L\left( \one{E} \alpha + \one{E^{\sf c} }\kappa  ,\zeta(N) \right)
\end{displaymath}

We will also consider a scenario where the $X_{1}, \ldots, X_{N}$ stay
independent but may not always be identically distributed.
\begin{defn}[Type $B$ random environment]
\label{def:alpha-random-one}
Suppose a population evolves according to Definition~\ref{hschwpop}.
Fix $0 < \alpha \le 1$ and $\kappa \ge 2$.  Recall \eqref{eq:28} and
write $ E_{1}$ for the event ($[N] = \set{1,2,\ldots, N}$ by
Definition~\ref{def:notation}) when 
there exists exactly one  $i \in [N]$ such that $X_i \vartriangleright
\mathds L(\alpha, \uN)$,  and  $X_{j} \vartriangleright
\mathds L(\at, \uN)$  for all $j \in [N]\setminus \{i\}$. When
$E_{1}$ occurs, the `lucky' individual (the index $i$)  is  picked uniformly at random.
When event $ E _{1}^{\sf c}$ is in force
 $\at$ replaces $\alpha$ in $ E_{1}$ ($X_{i} \vartriangleright
\mathds L(\kappa, \zeta(N))$ for all $i\in [N]$)  from
Definition~\ref{def:haplrandomalpha}.  Let
$(\overline \varepsilon_{N})_{N\in \IN}$ be a sequence with
$0 < \overline \varepsilon_{N} < 1$ for all $N\in \N$.   It may hold that
$\overline \varepsilon_{N}\to 0$ as $N\to \infty$.  Suppose 
$\prob{ E_{1} } = \overline \varepsilon_N$, and  $\prob{ E_{1}^{\sf c}}  = 1 - \overline \varepsilon_N$.
\end{defn}
Definition~\ref{def:alpha-random-one} says that when the environment
turns favorable for producing potential offspring through $\alpha$
(event $E_{1}$ occurs), exactly one individual sampled uniformly at
random will do so, i.e.\ with probability $1/N$ it holds that 
\begin{displaymath}
X_{i} \vartriangleright \mathds L\left( \one{E_{1}} \alpha +
\one{E_{1}^{\sf c} }\kappa  ,\zeta(N) \right),\quad X_{j}
\vartriangleright \mathds L(\kappa, \zeta(N)) \text{ for all  $j\in [N]\setminus \set{i}$}
\end{displaymath}
for all $i\in [N]$.  
Since the `lucky' individual is sampled at random when event $E_{1}$
occurs the $X_{1},\ldots, X_{N}$ are exchangeable and we can use the
results of \citep{MS01} for proving convergence of $\set{\xi^{n,N}}$.
As in Definition~\ref{def:haplrandomalpha} the bound $\uN$ stays the
same between $ E_{1}$ and $ E_{1}^{\sf c}$.  We will identify
conditions on $\overline\varepsilon_{N}$ (see \eqref{eq:51} in
Lemma~\ref{lm:cNrandomalphaone}) such that 
$\set{\xi^{n,N}\svigi{\lfloor t/c_{N} \rfloor }; t \ge 0}$ converges
to a nontrivial limit.  In Theorems~\ref{thm:haplrandomalpha} and
\ref{thm:hapl-alpha-random-one} the law for the number of potential
offspring is as in \eqref{eq:PXiJ}. Write (for $\kappa \ge 2$)
\begin{equation}
\label{eq:CNmap}
C_{\kappa}^{N} :=    \one{\kappa > 2}N + \one{\kappa = 2} {N}{/ \log N} 
\end{equation}
 Section~\ref{sec:proof-thm-randomallhaploidnewmodel}
contains a proof of Theorem~\ref{thm:haplrandomalpha}.
\begin{shaded*}
\begin{thm}[Evolution \breyta{in a type $A$ random environment (Definition~\ref{def:haplrandomalpha})}]
\label{thm:haplrandomalpha}%
Suppose a population evolves as defined in Definitions~\ref{hschwpop}
and \ref{def:haplrandomalpha} \breyta{(type $A$ random
environment)}  with the law for the
$X_{1},\ldots, X_{N}$ given by \eqref{eq:PXiJ}, that
Assumption~\ref{a:mN1plus} holds, $g_{\alpha}^{(\infty)} = f_{\alpha}^{(\infty)}$, and
$\underline{g_{\kappa}} (2) = \overline{f_{\kappa}}(2)$ (recall
\eqref{eq:isfg}).   Then
$\{ \xi^{n,N}(\lfloor t/c_{N} \rfloor ); t \ge 0 \}$ converges in the
sense of convergence of finite-dimensional distributions to
$\{\xi^{n}\} \equiv \{ \xi^{n}(t); t \ge 0 \}$ where $\{ \xi^{n} \}$
is as specified in each case.
\begin{enumerate}
\item\label{item:1}  Suppose    $\uN / N
\to 0$ and  $1 \le \alpha < 2$.       Then 
$\{ \xi^{n}\}$  is  the Kingman coalescent.
\item\label{item:2} Suppose $1 \le \alpha < 2$ and $\uN/N \gneqq 0$
(recall \eqref{eq:26} in Definition~\ref{def:not}).  Take
$(\varepsilon_{N})_{N}$ from Definition~\ref{def:haplrandomalpha} 
    \breyta{with}    \breyta{  $\varepsilon_{N} = cN^{\alpha -
  2}\svigi{\one{\kappa > 2} + \one{\kappa = 2}\log N  }$ for some
  fixed $c>0$ .}  
     Then
$\{\xi^{n}\}$ is the
$\delta_{0}$-Beta$(\gamma, 2-\alpha,\alpha)$-coalescent defined in
Definition~\ref{def:KingmanBeta} with $1 < m_{\infty} < \infty$
(recall \eqref{eq:5}). The transition rate for a $k$-merger of a  given group of
$k \in \{2,3,\ldots, n\}$ blocks is
\begin{align}
\label{eq:ratesarandall}
\lambda_{n,k} = \left( \one{k=2} C_{\kappa}  +  { m_{\infty}^{-\alpha} \alpha c f_{\alpha }^{(\infty)} }  B(\gamma, k-\alpha,n - k+ \alpha)\right)/  {C_{\kappa, \alpha, \gamma}}
\end{align}
with   $B(p,a,b) \equiv \int_{0}^{1} \one{0 < t \le p}
t^{a-1}(1-t)^{b-1}dt$ denoting   the beta function ($a,b > 0$, $0< p
\le 1$)  and  (recalling  $\overline{f_{\kappa}}(2),
f_{\alpha}^{(\infty)}$ from   \eqref{eq:isfg})
\begin{subequations}
\begin{align}
\label{eq:27}
\gamma & =  \one{\frac{\uN}{N} \to K } \frac{K}{m_{\infty} + K} + \one{ \frac{\uN}{N} \to \infty }  \\
\label{eq:23}
C_{\kappa, \alpha,\gamma} &  =   C_{\kappa}   +  \frac{ \alpha c f_{\alpha}^{(\infty)} }{m_{\infty}^{\alpha}}B(\gamma,2-\alpha,\alpha)  \\
C_{\kappa} & =    \frac{2\overline{ f_{\kappa}}(2)}{ m_{\infty}^{2}} \left(  \one{\kappa = 2} + \one{\kappa > 2} \frac{c_{\kappa}}{2^{\kappa}(\kappa - 2)(\kappa - 1)} \right) \label{eq:37}   \\
&      \underline{ g_{\kappa}} +  \underline{ g_{\kappa}}  2^{1-\kappa}/(\kappa - 1)  < m_{\infty} <     \overline{ f_{\kappa}} +   \overline{ f_{\kappa}} /(\kappa - 1)  \label{eq:49}
\end{align}
\end{subequations}
where in \eqref{eq:37}   we have $\kappa + 2 < c_{\kappa} <
\kappa^{2}$ for  $\kappa > 2$.  
\item \label{item:4} Suppose $0 < \alpha < 1$ and
$\uN/N^{1/\alpha} \to \infty$ and $\varepsilon_{N} = c/C_{\kappa}^{N}$
with $c > 0$ fixed.  Then $\{ \xi^{n} \}$ is the
$\delta_{0}$-Poisson-Dirichlet$(\alpha,0)$ coalescent defined in
Definition~\ref{def:kingman-poisson-dirichlet} with transition rates (recall
\eqref{eq:ratekpd})
\begin{equation}
\label{eq:39}
\lambda_{n; k_{1}, \ldots, k_r; s} = \frac{C_{\kappa} }{C_{\kappa} + c(1-\alpha)  }\one{r=1, k_{1}=2} +  \frac{c }{C_{\kappa} + c(1-\alpha) }p_{n;k_{1}, \ldots, k_{r};s}
\end{equation}
with  $C_{\kappa}$  as in \eqref{eq:37} and  $p_{2;2;0} = 1-\alpha$ so that $\lambda_{2;2;0} = 1$.   
\end{enumerate}
In all cases (recall \eqref{simseq} in Definition~\ref{def:notation}
and $c_{N}$ from Definition~\ref{def:cNhapl} and $C_{\kappa}^{N}$ from
\eqref{eq:CNmap}) as $N\to \infty$
\begin{equation}
\label{eq:cNthmrandall}
C_{\kappa}^{N}  c_{N}  \overset{c}{\sim}      1  
\end{equation}
\end{thm}
\end{shaded*}
In Theorem~\ref{thm:haplrandomalpha} the population evolves
according to Definition~\ref{def:haplrandomalpha}, where everyone produces
potential offspring through $\alpha$ (recall $\alpha < \at$) when the
environment turns favorable. In Theorem~\ref{thm:hapl-alpha-random-one}
the population evolves according to Definition~\ref{def:alpha-random-one}
where exactly one randomly picked individual produces potential
offspring through $\alpha$ when the conditions are favorable; thus the
$X_{1},\ldots, X_{N}$  may not always be identically
distributed.   \breyta{ Suppose  $\svigi{\overline \varepsilon _{N}}_{N}$
from
Definition~\ref{def:alpha-random-one}   takes the form as in  \eqref{eq:51} in
Lemma~\ref{lm:cNrandomalphaone},
\begin{displaymath}
\overline \varepsilon _{N} = \begin{cases}
cN^{\alpha - 1} & \text{ when $0 < \alpha < 1$ and  $\kappa > 2$ } \\  
cN^{\alpha - 1}\log N & \text{ when  $0 < \alpha < 1$ and  $\kappa = 2$ } \\
\varepsilon & \text{when  $\alpha = 1$ and  $\kappa > 2$}
\end{cases}
\end{displaymath}
with $c > 0$ and $0 < \varepsilon < 1$ both fixed.
}

Section~\ref{sec:proof-thm-random-alpha-one} contains a proof
of Theorem~\ref{thm:hapl-alpha-random-one}.
\begin{shaded*}
\begin{thm}[Evolution \breyta{in a type $B$ random environment} (Definition~\ref{def:alpha-random-one})]
\label{thm:hapl-alpha-random-one}
Suppose a population evolves as defined in Definitions~\ref{hschwpop}
 and \ref{def:alpha-random-one} \breyta{(type $B$ random
 environment)},  Assumption~\ref{a:mN1plus} holds, and $g_{\alpha}^{(\infty)} = f_{\alpha}^{(\infty)}$ (recall~\eqref{eq:isfg}).   
Take $(\overline \varepsilon_{N})_{N}$ from
Definition~\ref{def:alpha-random-one} as in \eqref{eq:51} in
Lemma~\ref{lm:cNrandomalphaone}. Then
$\{\xi^{n,N}(\lfloor t/c_{N} \rfloor); t \ge 0\}$ converges in the
sense of convergence of finite-dimensional distributions  to
$\{ \xi^{n}\} \equiv \{\xi^{n}(t); t \ge 0 \}$ where $\{ \xi^{n}\}$ is
as specified in each case.
\begin{enumerate}
\item \label{item:7} Suppose $\uN/ N\to 0$.  Then
$\{\xi^{n}\}$ is the Kingman coalescent.
 \item \label{item:9} Suppose $0 < \alpha \le 1$ and $\uN/N \gneqq 0$
 (recall \eqref{eq:26} in Definition~\ref{def:not}).  Then $\{\xi^{n}\}$ is
 the $\delta_{0}$-Beta$(\gamma, 2-\alpha,\alpha)$-coalescent defined in
 Definition~\ref{def:KingmanBeta} with $1 < m_{\infty} < \infty$ (recall
 \eqref{eq:5}).  The transition rates are as in
 \eqref{eq:ratesarandall} with $c_{\alpha}$ replacing $c$ where
 $c_{\alpha} = \one{\alpha = 1}\varepsilon + \one{0 < \alpha < 1}c$
 where $0 < \varepsilon < 1$ and $c > 0 $ as in \eqref{eq:51}.
\end{enumerate}
Furthermore, \eqref{eq:cNthmrandall} holds in both cases.   
\end{thm}
\end{shaded*}

\begin{remark}[The behaviour of  $\uN$]
\label{rm:uNincreasingwithN}
When $\uN/N \to \infty$ in Theorems~\ref{thm:haplrandomalpha} and
\ref{thm:hapl-alpha-random-one} (or $\uN/N^{1/\alpha} \to \infty$
as in Case~\ref{item:4} of Theorem~\ref{thm:haplrandomalpha}) it is with
the understanding that $\log(\uN) \overset{c}{\sim} \log N $ (as
$N\to \infty$), so that $\uN$ is at most of the form $N^{1 + \eta}$
for $0 < \eta \ll 1$ (see
Lemma~\ref{lm:law-large-numbers-alpha-random-one}).  One could also
think of $\uN \overset{c}{\sim} N\log N$ (or $N^{1/\alpha}\log N$ in
Case~\ref{item:4} of Theorem~\ref{thm:haplrandomalpha}).  We do not see
$\uN$  increasing exponentially (or faster) with $N$.
 \end{remark}

\subsubsection{\breyta{A general random environment}}\label{sec:gener-rand-envir}
\breyta{
This section is a response to the question \newline
\textit{`If I take 2 (or m, integer) different reproduction models
(Cannings models) that have known coalescent limits and mix them from
generation to generation, do they have a coalescent limit (with which
time scaling), what is the limit and is it e.g. a Lambda coalescent
with Lambda a sum of the Lambdas of the underlying different
reproduction models (assuming their limits are Lambda-coalescents)?'}
\newline
We have considered random environments each split into two
possibilities.   When warranted by the biology of the study system at
hand,  one may want to consider  a random
environment split into  several (say $m\ge 2$)  events, where  each
event  entails  reproduction according to a given Cannings model with
known  coalescent limit.
For example,
one may extend  the   type $A$ random environment
(Definition~\ref{def:haplrandomalpha}) to one where, with $0<\alpha<1$
and $1 \le \beta < 2$ both fixed, 
\begin{displaymath}
X_{1},\ldots,X_{N} \vartriangleright \begin{cases}
\mathds L (\alpha,  \zeta(N)) & \text{with probability $\varepsilon_{N}$} \\
\mathds L (\beta,  \zeta(N)) & \text{with probability $\varepsilon_{N}^{\prime}$} \\
\end{cases}
\end{displaymath}
such that    $0 < \varepsilon_{N} + \varepsilon_{N}^{\prime} < 1$ for
all $N\in \N$,  and
with $\kappa \ge 2$ fixed, 
\begin{displaymath}
X_{1},\ldots, X_{N} \vartriangleright \mathds L(\kappa,\zeta(N))
\end{displaymath}
with probability $1 - \varepsilon_{N} - \varepsilon_{N}^{\prime}$.
Choosing  $C_{\kappa}^{N} \varepsilon_{N} \overset{c}{\sim}  1$, and
$ N^{\beta - 2} \varepsilon_{N}^{\prime}  \overset{c}\sim 1  $ (recall~\eqref{eq:33}
in Lemma~\ref{lm:cNhaplrandomall}), our calculations  show  we would
obtain a continuous-time  $\Xi$-coalescent with an atom at zero and  with the  transition rate
being a linear combination of the rates given in
~\eqref{eq:ratesarandall} and~\eqref{eq:39} in
Theorem~\ref{thm:haplrandomalpha} with timescale  proportional to
$C_{\kappa}^{N}$~\eqref{eq:CNmap}.   In contrast,  the continuous-time
coalescent
considered in Appendix~\ref{sec:betapoissondirichlet} runs on a
timescale proportional to  $N^{\beta-1}$ generations ($1 < \beta < 2$)
and so runs into the problem  of  recovering  mutations observed in a
given sample of DNA gene copies, especially for estimates of $\beta$
close to 1. Barring any additional  complications such as  separation of
timescales due to e.g.\   selfing \citep{M98} or diploidy
\citep{MS03,BBE13,BLS15}, a central condition required to hold  to prove convergence  to
a continuous-time $\Xi$-coalescent is the  existence of the limits in
\eqref{xmoments} \citep[Equation~16]{MS01} (see
\citep[Equation~6]{BLS15} for the diploid version).   Consider  a
random environment consisting of $m$ events $E^{(i)}$,  where event $E^{(i)}$  occurs with probability
$\varepsilon_{N}^{(i)}$ such that  when  $E^{(i)}$ occurs the
(haploid panmictic)  population  reproduces  according  to some given  Cannings
reproduction scheme  with timescale $c_{N}^{(i)}$   such that
$\sum_{i}\varepsilon_{N}^{(i)} = 1$ and  for at least one index $i$ 
 it holds $c_{N}^{(i)} \to 0$.  Moreover,  suppose each Cannings model
 has a known  coalescent limit  with transition rates
 $\lambda_{n;k_{1},\ldots,k_{r}}^{(i)}$  \eqref{xirate}  driven by a measure
 $\Xi^{(i)}$  determined by the Cannings model specified in event
 $E^{(i)}$.           Then 
\begin{displaymath}
 \frac{\EE{\svigi{\nu_{1}}_{k_{1}}\cdots \svigi{\nu_{r}}_{k_{r}}}}{ N^{k_{1}+\cdots + k_{r} - r} c_{N} } =   \sum_{i=1}^{m} \frac{\EE{\svigi{\nu_{1}}_{k_{1}}\cdots \svigi{\nu_{r}}_{k_{r}} | E^{(i)}}}{ N^{k_{1}+\cdots + k_{r} - r} c_{N} }\varepsilon_{N}^{(i)}
\end{displaymath}
and Condition~\eqref{xmoments} becomes equivalent to checking
that 
\begin{displaymath}
\lim_{N\to \infty} \frac{\EE{\svigi{\nu_{1}}_{k_{1}}\cdots \svigi{\nu_{r}}_{k_{r}} | E^{(i)} }}{ N^{k_{1}+\cdots + k_{r} - r} c_{N} }\varepsilon_{N}^{(i)}
\end{displaymath}
exists for all $i\in [m]$ and $r\in \N$ and
$k_{1}, \ldots, k_{r} \ge 2$ (and with the  $\varepsilon_{N}^{(i)}$
suitably tuned; recall that $\nu_{1},\ldots,\nu_{N}$ are the random
numbers of offspring of the current $N$ individuals in a haploid
panmictic population of constant size $N$).    Convergence to a continuous-time
$\Xi$-coalescent with transition rates as a linear combination of the
individual rates $\lambda_{n;k_{1},\ldots,k_{r}}^{(i)}$ then follows
by \citep[Theorem~2.1]{MS01}.  However, even though extending  the
simple random environments  behind  Theorems~\ref{thm:haplrandomalpha}
and \ref{thm:hapl-alpha-random-one} to more general settings is
straightforward,    the biology of the study
system at hand (and a general principle of parsimony to avoid
overfitting) should dictate the modeling.}

\subsection{Comparing  processes}
\label{sec:numerics}

In this section we use simulations to compare functionals of
coalescents $\set{\xi^{n}}$ and ancestral processes $\set{\xi^{n,N}}$
in the domain of attraction of $\set{\xi^{n}}$ as given in
Theorems~\ref{thm:haplrandomalpha} and
\ref{thm:hapl-alpha-random-one}. We are interested in \breyta{
comparing  $\set{\xi^{n}}$ and $\set{\xi^{n,N}}$.}

Let $|A|$ be the number of elements in  a given finite  set $A$, write
\\ 
$\tau^{N}(n) := \inf \set{j \in \N_{0} : |\xi^{n,N}(j)| = 1}$, and 
$\tau(n) := \inf \set{t \ge 0 : | \xi^{n}(t)| = 1
}$. \breyta{For any given ancestral process
$\set{\xi^{n,N}}$ and coalescent $\set{\xi^{n}}$}  consider the
functionals
\begin{equation}
\label{eq:36}
\begin{split}
L_{i}^{N}(n) := \sum_{j=0}^{\tau^{N}(n)}  | \set{\xi \in \xi^{n,N}(j) : |\xi|  = i }|, &  \quad   L^{N}(n) :=  \sum_{j=0}^{\tau^{N}(n)} | \xi^{n,N}(j)| \\
L_{i}(n) := \int_{0}^{\tau(n)} |\set{\xi \in \xi^{n}(t) : |\xi| = i }| dt, &  \quad L(n) :=  \int_{0}^{\tau(n)} | \xi^{n}(t)| dt
\end{split}
\end{equation}
\citep{BLS15}.   \breyta{Viewing
$\set{\xi^{n,N}}$ and $\set{\xi^{n}}$ as gene genealogies (trees) the
random variables  
$L_{i}^{N}(n)$ resp.\  $L_{i}(n)$ can be interpreted as the random  length of branches  supporting  $i \in [n-1]$
leaves when the  gene genealogy with $n$ leaves  is governed by  a
given ancestral process  $\set{\xi^{n,N}}$ resp.\ coalescent
$\set{\xi^{n}}$. }

\begin{figure}[htp]
\centering
\includegraphics[scale=1]{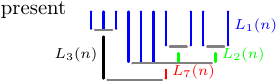}
\caption{\breyta{ The  functionals $L_{i}(n)$  ($L_{i}^{N}(n)$)   for an example
gene genealogy.} }
\label{fig:illL}
\end{figure}
\breyta{Given the example gene genealogy for $n=10$  in  Figure~\ref{fig:illL},
 $L_{1}(n)$ (resp.\ $L_{1}^{N}(n)$)  would be the total  length of the
 blue branches,  $L_{2}(n)$ the green branches,  $L_{3}(n)$ the black 
 branch, and  $L_{7}(n)$ the red branch.      }

It holds that 
$L^{N}(n) = L_{1}^{N}(n) + \cdots + L_{n-1}^{N}(n)$, and
$L(n) = L_{1}(n) + \cdots + L_{n-1}(n)$.
Define, for all  $i \in [n-1]$,
\begin{equation}
\label{eq:54}
\begin{split}
R_{i}(n) &  :=  \frac{L_{i}(n)}{L(n) }, \quad  R_{i}^{N}(n)   :=  \frac{L_{i}^{N}(n)}{ L^{N}(n)   } \\
\end{split}
\end{equation}
Both $R_{i}^{N}(n)$ and $R_{i}(n)$ are well defined since $L(n) > 0$
and $L^{N}(n) \ge n$ both almost surely.  We will compare
approximations of $\EE{R_{i}(n)}$ to approximations of
$\EE{R_{i}^{N}(n)}$ for corresponding processes.  Interpreting
$\set{\xi^{n}}$ and $\set{\xi^{n,N}}$ as `trees' one may view
$R_{i}(n)$ and $R_{i}^{N}(n)$ as random `relative branch lengths'.
Under the infinitely-many sites mutation model \citep{K1969,W1975} the
site-frequency spectrum, or the count of mutations of given
frequencies in a sample, corresponds directly to branch lengths
\citep{EBBF2015}.

We will use the notation (recall \eqref{eq:54})
\begin{equation}
\label{eq:functionals}
\begin{split}
 & \varrho_{i}^{N}(n) \equiv  \EE{R_{i}^{N}(n)}, \quad  \varrho_{i}(n) \equiv \EE{R_{i}(n)}, \quad \rho_{i}^{N}(n)  \equiv \EE{ \widetilde R_{i}^{N}(n) } \\
& \overline \varrho_{i}^{N}(n) \approx  \EE{R_{i}^{N}(n)}, \quad   \overline \varrho_{i}(n) \approx \EE{R_{i}(n)}, \quad    \overline \rho_{i}^{N}(n) \approx  \EE{ \widetilde R_{i}^{N}(n) } 
\end{split}
\end{equation}
where    $\varrho_{i}^{N}(n)$, $\varrho_{i}(n)$, and $\rho_{i}^{N}(n)$
are   functionals we are interested in, and $\overline
\varrho_{i}^{N}(n)$,  $\overline \varrho_{i}(n)$, and  $\overline
\rho_{i}^{N}(n)$ are the  corresponding approximations.  The
quantities  $\rho_{i}^{N}(n)$ are relative branch lengths read off fixed
complete sample trees  
(conditioning on the population ancestry and averaging over the
ancestries; \breyta{see \S~\ref{sec:comparingB}}).

  Various evolutionary histories may yield similar
site-frequency spectra, thus giving rise to non-identifiability issues
\citep{Myers2008,Bhaskar2014,Rosen2018,freund2023interpreting}.  The
site-frequency spectrum 
may, nevertheless, be used in 
conjunction with other statistics in inference
\citep{freund2021impact}.

Whenever $\set{\xi^{n,N}}$ converges to $\set{\xi^{n}}$ it should
follow that $\EE{R_{i}^{N}(n)}$ converges to $\EE{R_{i}(n)}$ as
$N\to \infty$ (and with $n$ fixed).  Moreover, one would hope that the
functionals would  be in good agreement for  all `not too
small' $N$.

In Appendix~\ref{sec:relat-expect-sfs} we
consider the  quantity (recall $L(n) >0$ almost surely) 
\begin{equation}
\label{eq:varphi}
\varphi_{i}(n) :=   \frac{\EE{L_{i}(n)} }{\EE{L(n)} }, \quad i\in [n-1], 
 \end{equation}
 for the $\delta_{0}$-Beta$(\gamma,2-\alpha,\alpha)$-coalescent (see
 Figure~\ref{fig:varphiA}).  Values of $\varphi_{i}(n)$ can be
 computed exactly   \breyta{using} a recursion  derived for
 $\Lambda$-coalescents \citep{BBE2013a}.

For sampling $\set{\xi^{n,N}}$ we will use the following example of
\eqref{eq:PXiJ}.  
Let $(p_{k}(a))_{k\in \N}$   be probability weights with
\begin{align}
\label{eq:40}
p_{k}(a) & =    \one{k \in [\uN]} h_{a}(\uN)  \left(  \frac{1}{k^{a}} -  \frac{1}{(1+k)^{a}}   \right), 
\end{align}
where $h_{a}(\uN)$ is such that $\sum_{k=1}^{\uN}p_{k}(a) = 1$ (i.e.\
$h_{a}(\uN) = (1 - (1+\uN)^{-a})^{-1}$  and $h_{a}(\uN) \to 1$ as
$\uN\to \infty$). If  $\uN < \infty$ for $N$ fixed  \eqref{eq:40}
is a special case of \eqref{eq:PXiJ} with $f_{a}(k) = g_{a}(k) = h_{a}(\uN)$
for all $k$  (and of \eqref{eq:4} when
$\uN = \infty$ for all $N$).   For the numerical experiments    we  will assume
that the sample comes from a population evolving according to
Definition~\ref{hschwpop}, and to either Definition~\ref{def:haplrandomalpha} or
\ref{def:alpha-random-one} as specified each time.   
\breyta{Briefly,  Definition~\ref{hschwpop}  describes the
evolution of a haploid panmictic population  evolving in
non-overlapping generations through the generation of potential
offspring,   and   Definitions~\ref{def:haplrandomalpha} (type $A$
random environment)  and
\ref{def:alpha-random-one} (type $B$ random environment)   describe the random environments used, the
difference being that  when Definition~\ref{def:alpha-random-one} is
in force the  numbers of potential offspring (always independent)  may not always be
identically distributed. 
  }
The probability weights of the random number of potential offspring are as
in \eqref{eq:40} with $a$ and $\uN$  as given in each case.
Case~\ref{item:2} of Theorem~\ref{thm:haplrandomalpha}
and Case~\ref{item:9} of Theorem~\ref{thm:hapl-alpha-random-one} then
show that choosing $(\varepsilon_{N})_{N}$, 
$(\overline\varepsilon_{N})_{N}$,  and $\uN$ suitably the
(time-rescaled)   ancestral
process converges to a
$\delta_{0}$-Beta$(\gamma,2-\alpha,\alpha)$-coalescent (recall
Definition~\ref{def:KingmanBeta}) with transition rates as in
\eqref{eq:ratesarandall} with $f_{\alpha}^{(\infty)} = 1$. With
$\kappa = 2$ the transition rate for a $k$-merger of 
 $k \in \set{2,3,\ldots, n}$ blocks is 
\begin{equation}
\label{eq:lambdankb}
\lambda_{n,k} =  \frac{1}{C_{\kappa, \alpha,\gamma}}  \binom{n}{k} \left( \frac 2 {m_{\infty}^{2}}  \one {k=2}  +  \frac{\alpha c_{\alpha}^{\prime} }{m_{\infty}^{\alpha} }  B(\gamma, k-\alpha, n-k+\alpha)\right)
\end{equation}
with $0 < \alpha < 2$,     $m_{\infty}$ from \eqref{eq:5} and   $\gamma$ from 
\eqref{eq:27} and  $B(\gamma,k-\alpha,n-k+\alpha)$ as in Definition~\ref{def:betacoal} and 
\begin{displaymath}
\begin{split}
 C_{\kappa, \alpha,\gamma}  &  =  2m_{\infty}^{-2}   +   \alpha c_{\alpha}^{\prime}  {m_{\infty}^{-\alpha}} B(\gamma,2-\alpha,\alpha)  
\end{split}
\end{displaymath}
so that $\lambda_{2,2} = 1$.  When $\alpha = 1$ we need to decide,
which of the \breyta{random environments} (Definitions~\ref{def:haplrandomalpha} or
\ref{def:alpha-random-one})  holds; let $c_{\alpha}^{\prime} := c$ if
Definition~\ref{def:haplrandomalpha} holds, and
$c_{\alpha}^{\prime} = c_{\alpha}$ if
Definition~\ref{def:alpha-random-one} holds where
$c_{\alpha} = \one{0 < \alpha < 1}c + \one{\alpha = 1}\varepsilon$ and
$0 < \varepsilon < 1$ is fixed.

We see, with
$\at \ge 2$,  
\begin{displaymath}
\begin{split}
  1  +  \frac{2^{1- \at }}{ \at - 1} + \OO{\frac{1}{\uN^{ \at - 1 }} } \le    \sum_{k=1}^{\uN} k\left( \frac{1}{k^{\at }} -    \frac{1}{(1+k)^{ \at }} \right)  \le 1 +  \frac{1}{\at - 1}   +         \OO{\frac{1}{\uN^{ \at - 1 }} }  
\end{split}
\end{displaymath}
Then, when a population evolves according to Definition~\ref{hschwpop}
and Definition~\ref{def:haplrandomalpha} with $\varepsilon_{N}$ as in
\eqref{eq:33} when Definition~\ref{def:haplrandomalpha} holds, and
$\overline \varepsilon_{N}$ as in \eqref{eq:51} when
Definition~\ref{def:alpha-random-one} holds we see that we can
approximate \breyta{$m_{\infty} = \lim_{N\to
\infty}\EE{X_{1}}$}  (recall \eqref{eq:5}) with {\sf m}
where $\kappa \ge 2$ and 
\begin{equation}
\label{eq:57}
\text{\sf m} =     1 +   \frac{1 + 2^{1- \at }}{ 2(\at - 1)}
\end{equation}
Thus, the $\delta_{0}$-Beta$(\gamma,2-\alpha,\alpha)$-coalescent
when $\kappa = 2$ involves the parameters,
$c_{\alpha}^{\prime}, \alpha,  \gamma$, in addition to $m_{\infty}$
(here approximated by  {\sf m} as in \eqref{eq:57}).

We also compare $\overline \varrho_{i}^{N}(n)$ to   $\overline \varrho_{i}(n)$  as predicted by the
$\delta_{0}$-Poisson-Dirichlet$(\alpha,0)$-coalescent (recall
Case~\ref{item:4} of Theorem~\ref{thm:haplrandomalpha}), where the
distribution of the number of potential offspring is as in
\eqref{eq:57} with $\uN = N^{1/\alpha}\log N$ (so that
$\uN/N^{1/\alpha} \to \infty$ as $N\to \infty$ as required).  Recall
from Case~\ref{item:4} of Theorem~\ref{thm:haplrandomalpha} the
transition rates of the
$\delta_{0}$-Poisson-Dirichlet$(\alpha,0)$-coalescent,
\begin{displaymath}
\lambda_{n;k_{1}, \ldots, k_{r}; s} =   \frac{C_{\kappa}}{C_{\kappa} + c(1-\alpha) }\one{r=1,k_{1}=2} +  \frac{c }{C_{\kappa} + c(1-\alpha) }p_{n;k_{1}, \ldots, k_{r};s}
\end{displaymath}
where $C_{\kappa}$ is as in \eqref{eq:37} with $m_{\infty}$
approximated as in \eqref{eq:57}.


\breyta{ In \S~\ref{sec:appr-eriNn}  we describe  an    algorithm  for
approximating
$\EE{R_{i}^{N}(n)}$ ($\varrho_{i}^{N}(n)$,  annealed mean relative branch lengths  predicted
by a given ancestral process, 
recall~\eqref{eq:36}, ~\eqref{eq:54}, 
~\eqref{eq:functionals}),    in \S~\ref{sec:appr-qeewid-r_inn}  for approximating   $\EE{\widetilde
R_{i}^{N}(n)}$ ($\rho_{i}^{N}(n)$,  quenched mean relative branch
lengths),  and in  \S~\ref{sec:sampl-from-delta-pd}  
 for  sampling from the
$\delta_{0}$-Poisson-Dirichlet$(\alpha,0)$-coalescent.}  
We used {\tt R} \citep{Rsystem} and the shell tool {\tt parallel}
\citep{tange11:_gnu_paral} for generating the graphs.

\subsubsection{\breyta{Approximating
$\EE{R_{i}^{N}(n)}$}}\label{sec:appr-eriNn}

 Recall \S~\ref{sec:numerics}, in particular   \eqref{eq:54}.   
In this section we describe the \breyta{(non-quenched,
annealed)}   algorithm for  sampling branch lengths
of a   gene
genealogy of a sample of size $n$  from a finite haploid panmictic   population of constant
size $N$.   Suppose the
population evolves according to Definition~\ref{hschwpop}.  Our algorithm
returns a realisation  $\left(\ell_{1}^{N}(n), \ldots,
\ell_{n-1}^{N}(n)\right)$  of  $\left( L_{i}^{N}(n), \ldots,
L_{n-1}^{N}(n) \right)$.  Since we are working with a  Markov sequence (a
Markov process   evolving in discrete time), and we are only interested in
the site-frequency spectrum,  we only  keep track of
the current block sizes  $\{ b_{1}(g), \ldots,
b_{m}(g) \}$, where $b_{i}(g)$ is the  size
(number of leaves the block is ancestral to)  of
block $i$ at time $g$, so that  $b_{i}(g) \in [n]$, and
$b_{1}(g) + \cdots + b_{m}(g) = n$ when there are
$m$ blocks.      Let $\overline \varrho_{i}^{N}(n)$  (initialised  to zero) 
denote an  estimate of
$\EE{R_{i}^{N}(n)}$ obtained in the following way:
\begin{enumerate}
\item   $\svigi{\ell_{1}^{N}(n), \ldots, \ell_{n-1}^{N}(n) } \leftarrow  (0, \ldots, 0) $  
\item   set the block sizes $b_{1}(0),\ldots, b_{n}(0)$  to 1 
\item while there are at least two blocks  repeat the following steps
in order
\begin{enumerate}
\item update the current branch lengths 
$\ell _{b}^{N}(n)
\leftarrow \ell_{b}^{N}(n) + 1$ for $b \in \{ b_{1}(g), \ldots,
b_{m}(g) \} $ given that there are $m$ blocks at time $g$
\item sample a realisation $x_{1}, \ldots, x_{N}$ of the juvenile
numbers  $X_{1}, \ldots, X_{N}$
\item sample number of blocks per family according to a multivariate
hypergeometric with parameters $m$ (the current number of blocks)  and  $x_{1}, \ldots, x_{N}$
\item merge blocks at  random  according to the numbers
sampled in (c) 
\end{enumerate}
\item update the estimate $\overline \varrho_{i}^{N}(n)$  of  $\EE{R_{i}^{N}(n)}$:
\begin{displaymath}
\overline \varrho_{i}^{N}(n) \leftarrow \overline \varrho_{i}^{N}(n) +    \frac{\ell_{i}^{N}(n) }{ \ell_{1}^{N}(n) + \cdots + \ell_{n-1}^{N}(n)  }
\end{displaymath}
\item after repeating steps (1) to (4) a given number ($M$)  of times
return  $ \overline \varrho_{i}^{N}(n) /M $ 
\end{enumerate}
\breyta{See   Figures~\ref{fig:compareERiNe2A},
~\ref{fig:increasingnkbalessone},  and
Figures~\ref{fig:deltanullBetaERiNaddA} 
and~\ref{fig:deltanullBetaERiNaddBB} in Appendix~\ref{sec:sampl-discr-trees}   for examples.  }

\subsubsection{\breyta{Approximating $\EE{\widetilde R_{i}^{N}(n)}$}}\label{sec:appr-qeewid-r_inn}

In this section we briefly describe an algorithm for approximating 
$\EE{ \widetilde R_{i}^{N}(n) }$ (quenched mean relative branch
lengths, $\rho_{i}^{N}(n)$, recall~\eqref{eq:functionals}).  
Figure~\ref{fig:quenchedannealedkpdbounded} \breyta{and
Figure~\ref{fig:QAaddDD} in Appendix~\ref{sec:estimatequenched} hold  examples}.  
Recall $[n]$ and $\IN$ from Definition~\ref{def:notation}.  Let
$(A_{i}(g))_{g\in \IN\cup \{0\}, i\in [N]}$ denote the ancestry of the
population where $A_{i}(g) \in [N]$ is the level of the immediate
ancestor of the individual occupying level $i$ at time $g$. We set
$A_{i}(0) = i$ for $i \in [N]$.  If $A_{i}(g) = A_{j}(g)$ for
$i \neq j$ the individuals occupying levels $i$ and $j$ at time $g$
derive from the same immediate ancestor. If the individual on level
$i$ at time $g$ produced $k$ surviving offspring then
$A_{j_{1}}(g+1) = \cdots = A_{j_{k}}(g+1) = i$.  Each individual 
 `points' to its' immediate ancestor.

\begin{equation}
\label{eq:38}
\begin{split}
& \qquad   \text{levels}  \\
\text{time} & \qquad  1  \qquad \ldots \qquad  \ell  \qquad \ldots \qquad N \\ 
g & \qquad   i  \qquad \ldots \qquad j \qquad \ldots \qquad k \\
\end{split}
\end{equation}
\breyta{In \eqref{eq:38},  the individual  occupying level   $1$ resp.\ $\ell$
resp.\ $N$ at time  $g$  derives  from the  immediate ancestor  occupying level
$i$ resp.\ $j$ resp.\ $k$ at time $g-1$.}

A `complete' sample tree is one where the leaves have a common
ancestor.   Let
\begin{displaymath}
 r_{1}^{N}(n,\mathds{A}^{(N,n)}), \ldots,
r_{n-1}^{N}(n,\mathds{A}^{(N,n)})
\end{displaymath}
denote the realised relative
branch lengths of a complete sample  tree whose ancestry is given by   $\mathds{A}^{(N,n)}$, then we estimate, for $i \in \{1,2,\ldots, n-1\}$,
\begin{equation}
\label{eq:estimate}
\EE{\widetilde  R_{i}^{N}(n) } \approx  \frac 1M \sum_{j=1}^{M}   r_{i}^{N}\left(n, \mathds{A}_{j}^{(N,n)}  \right)
\end{equation}
where $M$ is the number of experiments, the number of realised ancestries $ \mathds{A}^{(N,n)}$.

We summarize the algorithm to estimate  $\EE{ \widetilde R_{i}^{N}(n) }$. 
\begin{enumerate}
\item For each  experiment:
\begin{enumerate}
\item initialise the ancestry to $A_{i}(0) = i$ for $i\in [N]$
\item until a complete sample  tree is found:
\begin{enumerate}
\item draw a random sample, i.e.\ sample $n$ of $N$ levels at the most recent  time 
\item check if the tree of the given sample is complete,  if not discard the sample and  record the ancestry of  a new set of surviving offspring: 
\begin{enumerate}
\item   sample \breyta{numbers of } potential offspring  $X_{1},\ldots, X_{N}$  
\item   given $X_{1},\ldots, X_{N}$
potential offspring sample the surviving offspring uniformly at random
without replacement and update the ancestry; if the individual on level $i$ at time
$g$ produced $k$ surviving offspring then
$A_{j_{1}}(g+1) = \cdots = A_{j_{k}}(g+1) = i$.
\end{enumerate}
\end{enumerate}
\item given a complete tree read  the branch lengths off the tree  and merge  blocks according to the
ancestry; suppose two  blocks have  ancestors
on  levels $i$ and $j$ at time $g$,  if  $A_{i}(g) = A_{j}(g)$ the blocks 
are merged;
\item given the branch lengths  of a complete tree update the estimate
of $\EE{ \widetilde R_{i}^{N}(n) }$
\end{enumerate}
\item from $M$ realised  ancestries  $\mathds{A}_{1}^{(N,n)}, \ldots,
\mathds{A}_{M}^{(N,n)}$    return an estimate \eqref{eq:estimate}   of
$\EE{\widetilde R_{i}^{N}(n)  }$
\end{enumerate}

\subsubsection{\breyta{Sampling from the  $\delta_{0}$-Poisson-Dirichlet$(\alpha,
0)$-coalescent}}\label{sec:sampl-from-delta-pd}

In this section we briefly describe an algorithm for samping from the
$\delta_{0}$-Poisson-Dirichlet$(\alpha,0)$ coalescent (recall
Definition~\ref{def:kingman-poisson-dirichlet}). See
Figures~\ref{fig:kpdbounded},  and Figures~\ref{fig:dpdaddCC} and
~\ref{fig:dpdexamples} in Appendix~\ref{sec:samplingKPD}   for examples of approximations $\overline \varrho_{i}(n)$ of
$\EE{R_{i}(n)}$ obtained using the algorithm.

By Case~\ref{item:4} of  Theorem~\ref{thm:haplrandomalpha} the transition rate (where $2 \le
k_{1},\ldots, k_{r} \le n$, $\sum_{i}k_{i} \le n$, $s = n -
\sum_{j}k_{j}$) is 
\begin{equation}
\label{eq:2}
\begin{split}
\lambda_{n;k_{1},\ldots,k_{r};s} &  =  \one{r=1,k_{1}=2}  \binom{n}{2} \frac{C_{\kappa}}{C_{\kappa} + c(1-\alpha) } \\
& +  \binom{n}{k_{1} \ldots k_{r}\, s}  \frac{1}{\prod_{j=2}^{n}\svigi{ \sum_{i}\one{k_{i}=j}}! } \frac{c}{C_{\kappa} + c(1-\alpha) } p_{n;k_{1}, \ldots, k_{r};s}  \\
\end{split}
\end{equation}
where $p_{n;k_{1},\ldots, k_{r};s}$ is as in \eqref{eq:34}.  For each
$m \in \set{2,3,\ldots, n}$ and $r \in [\lfloor m/2 \rfloor]$ we list
all possible ordered merger sizes
$2 \le k_{1} \le \cdots \le k_{r} \le m$ with $\sum_{i}k_{i} \le m$;
the algorithms for listing mergers borrow from algorithms for listing
partitions of integers.  The time until a merger when $n$ blocks is
then an exponential with rate the sum of
$\lambda_{n;k_{1}, \ldots, k_{r};s}$ in~\eqref{eq:2} over the ordered
mergers, and the actual merger size(s) can be efficiently sampled
given a listing of all the mergers ordered in descending order
according to the rate.

\subsubsection{\breyta{Comparing $\overline \varrho_{i}^{N}(n) $ and $\overline
\varrho_{i}(n)$}}\label{sec:comparingA}


\breyta{Recall  $\overline \varrho_{i}^{N}(n)$ (approximations of annealed mean relative
branch lengths associated with an ancestral process)   and  $\overline
\varrho_{i}(n) $  (approximations of annealed mean relative branch lengths associated
with a coalescent)    from
~\eqref{eq:36},~\eqref{eq:54},~\eqref{eq:functionals}. }
In Figures~\ref{fig:compareERiNe2A}--\ref{fig:kpdbounded} we compare
$\overline \varrho_{i}^{N}(n)$ and $\overline \varrho_{i}(n)$ (black
lines)   when the population evolves as 
in Definition~\ref{hschwpop} and Definition~\ref{def:haplrandomalpha}
(\breyta{type $A$ random environment;} Figures~~\ref{fig:compareERiNe2A},~\ref{fig:kpdbounded}) and
Definition~\ref{def:alpha-random-one}
(\breyta{type $B$ random environment;}   Figure~\ref{fig:increasingnkbalessone}) when  $\set{\xi^{n}}$ is the 
$\delta_{0}$-Beta$(\gamma,2-\alpha,\alpha)$-coalescent 
(Figures~~\ref{fig:compareERiNe2A}--\ref{fig:increasingnkbalessone};
\breyta{see also Figures~\ref{fig:deltanullBetaERiNaddA}
and~\ref{fig:deltanullBetaERiNaddBB} in Appendix~\ref{sec:sampl-discr-trees}})
or the $\delta_{0}$-Poisson-Dirichlet$(\alpha,0)$-coalescent
(Figure~\ref{fig:kpdbounded}; \breyta{see also
Figure~\ref{fig:dpdaddCC} in Appendix~\ref{sec:samplingKPD}}).   When
$1 \le \alpha \le 3/2$ (Figure~\ref{fig:compareERiNe2A}) there is
discrepancy between  $\overline\varrho_{i}(n)$ (black lines) and
 $\overline \varrho_{i}^{N}(n) $   (coloured lines) regardless of the bound $\uN$.
In Figure~\ref{fig:compareERiNe2A} the case $\gamma = 1$
(Figures~\ref{fig:strjalNe3KBA},
\breyta{\ref{fig:strjalNe3KBAaddA},~\ref{fig:strjalNe3KBB},~\ref{fig:strjalNe3KBC}}) is compared
to the case $\uN = N\log N$ (so that one would have $\uN/N \to \infty$
but $\limsup_{N\to \infty} m_{N} < \infty$ by
Lemma~\ref{lm:randomalphamNbounded} taking
$\varepsilon_{N} = N^{\alpha - 2}\log N$ since $\kappa = 2$, recall
~\eqref{eq:33}), and the case $\gamma = 1/(1 + \breyta\m)$
(Figures~\ref{fig:strjalNe3KBD}--\ref{fig:strjalNe3KBFaddB}) to the case
$\uN = N$.  The estimates $\overline \varrho_{i}^{N}(n)$ do not agree with  $\overline \varrho_{i}(n)$  
 for the values of $\alpha$ considered when the population
evolves according to Definition~\ref{def:alpha-random-one}
(Figure~\ref{fig:increasingnkbalessone}); the agreement is somewhat
better when the coalescent is the
$\delta_{0}$-Poisson-Dirichlet$(\alpha,0)$-coalescent
(Figure~\ref{fig:kpdbounded}).  In Figure~\ref{fig:kpdbounded}
$\uN = N^{1/\alpha}\log N$ so that $\uN / N^{1/\alpha} \to \infty$ as
required for convergence to the
$\delta_{0}$-Poisson-Dirichlet$(\alpha,0)$-coalescent. However, this
means that as $\alpha$ approaches 0 the realised number of potential
offspring any single individual may produce rapidly increases (see
also Figure~\ref{fig:dpdexamples} in Appendix~\ref{sec:samplingKPD} for
examples comparing estimates of $\EE{R_{i}(n)}$ as predicted by the
$\delta_{0}$-Poisson-Dirichlet$(\alpha,0)$-coalescent for various
parameter values).  We have that $\varepsilon_{N}$ is proportional to
$N^{\alpha - 2}$ when Definition~\ref{def:haplrandomalpha} holds as
required by Lemma~\ref{lm:cNhaplrandomall} (and
$\overline \varepsilon_{N} \propto N^{\alpha - 1} $ when
Definition~\ref{def:alpha-random-one} holds as required by
Lemma~\ref{lm:cNrandomalphaone}); thus in both cases as $\alpha$
decreases the probability of seeing large families in the ancestral
process    decreases; however, one expects the large families to be (on
average) larger as $\alpha$ decreases  (in
Figure~\ref{fig:comparexactsim} in Appendix~\ref{sec:relat-expect-sfs}
we check that the algorithm for estimating $\EE{R_{i}(n)}$ correctly
approximates  $\varphi_{i}(n)$~\eqref{eq:varphi}).




 \breyta{A $U$-shaped
site-frequency spectra (excess of low-frequency and high-frequency
derived alleles relative to predictions of the Kingman coalescent) is
sometimes seen as a characteristic of multiple-merger coalescents, in
particular $\Lambda$-coalescents \citep{freund2023interpreting}.
Indeed, $\Lambda$-coalescents can predict $U$-shaped spectra
\citep{BBE2013a}.  Our results (see Figure~\ref{fig:varphiA} in Appendix~\ref{sec:relat-expect-sfs}) do,
however, reveal a broader range of shapes, including $L$-shaped
spectra (displaying an excess of low-frequency alleles) similar to
those predicted by a time-changed Kingman-coalescent as derived from a
panmictic population evolving according to the Wright-Fisher model and
continually increasing in size \citep{Donnelly1995}.  Thus, the
statistical power of identifying multiple-merger coalescents from
time-changed Kingman-coalescents using the site-frequency spectrum may
have been somewhat diminished \citep{EBBF2015,Koskela2018}. In
particular, Figure~\ref{fig:varphiA} illustrates the effect of the
cutoff parameter $\gamma$ on the site-frequency spectrum, even for
small $\alpha$ the site-frequency spectrum as predicted by the
$\delta_{0}$-Beta$(\gamma,2-\alpha,\alpha)$-coalescent stays
$L$-shaped if $\gamma$ is small enough.    In cases involving the
Poisson-Dirichlet coalescents and reflecting high skewness of the
offspring distribution (Figures~\ref{fig:betapoissondiriA},
\ref{fig:dpdexamples}, and \ref{fig:adddpdFF}) the site-frequency
spectrum exhibits a ``plateau'' (in particular
Figures~\ref{fig:betapoissonB}, \ref{fig:betapoissonC},
\ref{fig:betapoissonE}, \ref{fig:betapoissonF}, ~\ref{fig:adddpdFFg},
~\ref{fig:dpdexamplesc}).   Recall that  a 
Poisson-Dirichlet$(\alpha,0)$-coalescent is a $\Xi$-coalescent
admitting  simultaneous multiple mergers;  when $n$ blocks  there can be up to
$\lfloor n/2 \rfloor$ simultaneous mergers, thus as the skewness
increases we can expect to see more  simultaneous mergers.  In
contrast  the Beta-coalescents admit only asynchronous mergers (at
most one set   of blocks merges at any given time).  Moreover,     in
 the  cases where the ``plateau'' appears the 
relative length of  external branches
$(\overline \varrho_{1}(n)$) is increased.     }


\subsubsection{\breyta{Comparing  $\overline\varrho_{i}^{N}(n)$ and $\overline
\rho_{i}^{N}(n)$}}\label{sec:comparingB}

\begin{figure}[htp]
\centering
\captionsetup[subfloat]{labelfont={scriptsize,sf,md,up},textfont={scriptsize,sf}}
\subfloat[$\alpha = \gamma = 1$]{\label{fig:strjalNe3KBA}\includegraphics[angle=0,scale=.5]{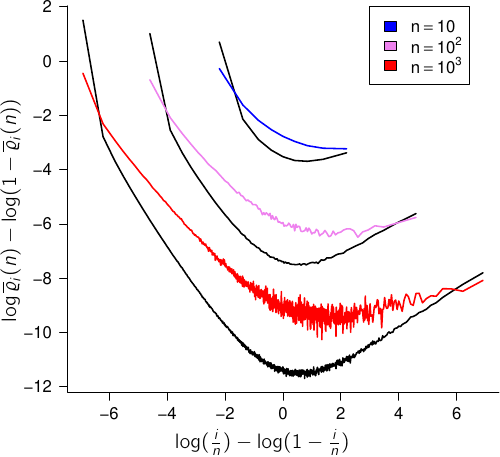}}
\subfloat[$\alpha = 1.5$, $\gamma =  \tfrac{1}{1 + \m} $]{\label{fig:strjalNe3KBF}\includegraphics[angle=0,scale=.5]{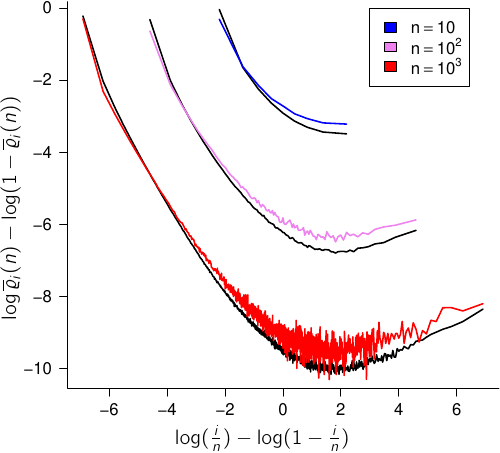}}
\caption{Definition~\ref{def:haplrandomalpha} \breyta{(type
$A$ random environment)}  and the 
$\delta_{0}$-Beta$(\gamma,2-\alpha,\alpha)$-coalescent.  Comparing
$\overline\varrho_{i}(n)$  (\breyta{estimates of mean relative branch lengths
predicted by the  $\delta_{0}$-Beta$(\gamma,2-\alpha,\alpha)$-coalescent},  black lines)   and
$\overline\varrho_{i}^{N}(n) $  (\breyta{estimates of  mean
relative branch lengths predicted by $\set{\xi^{n,N}}$ in attraction
of the  $\delta_{0}$-Beta$(\gamma,2-\alpha,\alpha)$-coalescent,
recall Figure~\ref{fig:illL} and  \eqref{eq:functionals}})     when   $N = 10^{3}$,  
$\kappa = 2$, $c = 1$, and $\gamma$ as shown;     black lines are 
$\overline \varrho_{i}(n)$   for sample size $n$ as
shown with rates as in \eqref{eq:ratesarandall} in Theorem~\ref{thm:haplrandomalpha},  coloured lines  are
 $\overline \varrho_{i}^{N}(n)$ for a sample from a population of finite size $N$
evolving according to Definition~\ref{hschwpop} and Definition~\ref{def:haplrandomalpha}
with the potential offspring distributed as in \eqref{eq:40} and with
$\varepsilon_{N} = cN^{\alpha - 2}\log N$ as in \eqref{eq:33} in
Lemma~\ref{lm:cNhaplrandomall}; the case $\gamma = 1$ is compared to
$\uN = N\log N$, and $\gamma = 1/(1 + \m) $ to $\uN = N$ with
$\m$ as in \eqref{eq:57} approximating $m_{\infty}$;  $\overline \varrho_{i}^{N}(n)$
 from $10^{4}$ experiments. \breyta{Further graphs are in
 Figure~\ref{fig:deltanullBetaERiNaddA} in
 Appendix~\ref{sec:sampl-discr-trees}.   The scale of the ordinate
 (vertical axis) may vary between graphs.  } }
\label{fig:compareERiNe2A}
\end{figure}

\begin{figure}[htp]
\centering
\captionsetup[subfloat]{labelfont={scriptsize,sf,md,up},textfont={scriptsize,sf}}
\subfloat[$\alpha = 0.01$, $\gamma = 1$]{\includegraphics[angle=0,scale=.5]{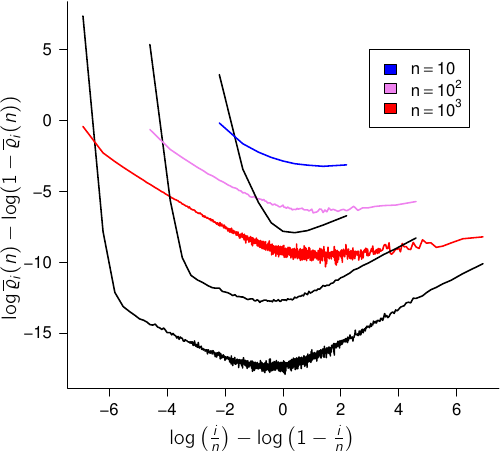}}
\subfloat[$\alpha = 0.5$, $\gamma = \tfrac{1}{1 + \m}$]{\includegraphics[angle=0,scale=.5]{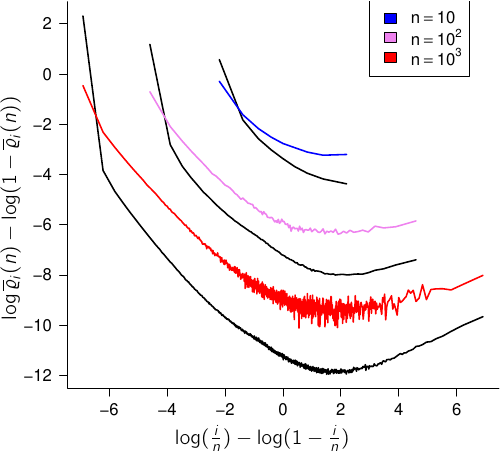}}
\caption{Definition~\ref{def:alpha-random-one} \breyta{(type
$B$ random environment)}  and the
$\delta_{0}$-Beta$(\gamma,2-\alpha,\alpha)$-coalescent. Comparing
$\overline\varrho_{i}(n)$  (\breyta{estimates of mean
relative branch lengths predicted by the
$\delta_{0}$-Beta$(\gamma,2-\alpha,\alpha)$-coalescent,}  black lines)
and $\overline \varrho_{i}^{N}(n)$ (\breyta{estimates of mean relative branch
lengths predicted by $\set{\xi^{n,N}}$ in attraction of the
$\delta_{0}$-Beta$(\gamma,2-\alpha,\alpha)$-coalescent,}  recall~\eqref{eq:functionals})  when $N=10^{3}$, $\kappa = 2$, $c = 1$, and
with $\alpha$, $\gamma$, and sampled size $n$ as shown; black lines
are  $\overline \varrho_{i}(n)$ with rates as in
~\eqref{eq:ratesarandall}, coloured lines are estimates of
$\EE{R_{i}^{N}(n)}$ for a sample from a population evolving according
to Definition~\ref{hschwpop} and Definition~\ref{def:alpha-random-one}
with the potential offspring distributed as in~\eqref{eq:40} and with
$\overline\varepsilon_{N} = cN^{\alpha - 1}\log N$ as in~\eqref{eq:51}
in Lemma~\ref{lm:cNrandomalphaone}; the case $\gamma = 1$ is compared
to $\uN = N\log N$, and the case $\gamma = 1/(1 + \m)$
compared to $\uN = N$ with $\m$  as in
~\eqref{eq:57} approximating $m_{\infty}$;  $\overline \varrho_{i}^{N}(n)$ from $10^{4}$
experiments. \breyta{Further graphs are in
Figure~\ref{fig:deltanullBetaERiNaddBB} in
Appendix~\ref{sec:sampl-discr-trees}.   The scale of the ordinate
 (vertical axis) may vary between graphs. } }
~\label{fig:increasingnkbalessone}
\end{figure}

\begin{figure}[htp]
\centering
\captionsetup[subfloat]{labelfont={scriptsize,sf,md,up},textfont={scriptsize,sf}}
\subfloat[$c = 1$]{\includegraphics[angle=0,scale=.6]{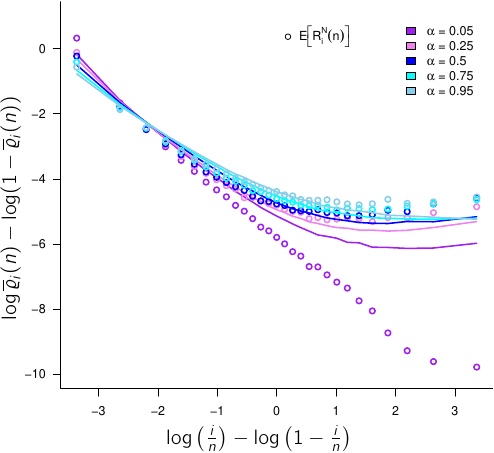}}
\subfloat[$c = 100$]{\includegraphics[angle=0,scale=.6]{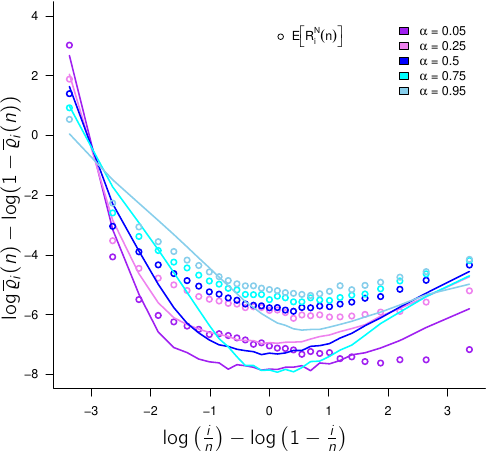}}
\caption{The $\delta_{0}$-Poisson-Dirichlet$(\alpha,0)$-coalescent.
Approximations  $\overline \varrho_{i}(n)$ \breyta{\eqref{eq:functionals}}   of $\EE{R_{i}(n)}$ 
(\breyta{mean relative branch lengths predicted by the    $\delta_{0}$-Poisson-Dirichlet$(\alpha,0)$-coalescent,} lines) predicted by the
$\delta_{0}$-Poisson-Dirichlet$(\alpha,0)$-coalescent compared to
$\overline\varrho_{i}^{N}(n)$ ( \breyta{mean relative branch
lengths predicted by $\set{\xi^{n,N}}$ in the attraction of the
$\delta_{0}$-Poisson-Dirichlet$(\alpha,0)$-coalescent,}  circles)  when the population evolves according
to Definition~\ref{def:haplrandomalpha} \breyta{(type $A$
random environment)}  with
$\uN = \infty$, $\varepsilon = c(\log N)/N$ for $c$ as shown,
$N=3000$, $\alpha$ as shown and $\kappa = 2$. The scale of the
ordinate (y-axis) may vary between the graphs.   
\breyta{See  \S~\ref{sec:sampl-from-delta-pd}} for an algorithm
for sampling from the
$\delta_{0}$-Poisson-Dirichlet$(\alpha,0)$-coalescent
\breyta{Further graphs are in 
Figure~\ref{fig:dpdaddCC} in Appendix~\ref{sec:samplingKPD}.   The scale of the ordinate
 (vertical axis) may vary between graphs.  }  }
\label{fig:kpdbounded}
\end{figure}

\breyta{Recall  $\overline \varrho_{i}^{N}(n)$ (approximations of
annealed mean relative branch lengths)  and  $\overline
\rho_{i}^{N}(n) $  (approximations of quenched  mean relative branch
lengths)     from
~\eqref{eq:36}, \eqref{eq:54}, \eqref{eq:functionals}. }
The approximations $\overline\varrho_{i}^{N}(n)$ in
Figures~\ref{fig:compareERiNe2A}--\ref{fig:kpdbounded} result from
ignoring the ancestral relations of the individuals in the population.
An alternative way of approximating `branch lengths' is to condition
on the ancestral relations of the individuals in the population
\citep{wakeley2012gene,Diamantidis2024}.  This idea follows from the
fact that as the population evolves ancestral relations are generated
through inheritance.  This means that the gene copies of a sample are
related through a single complete gene genealogy. Given the
\breyta{\it population ancestry}, the ancestral relations of
the entire population \breyta{of gene copies}  at the present and all previous times, and the
identity of the sampled gene copies, the complete sample tree is fixed
and known.  Denote by $ \EE{ \widetilde R_{i}^{N}(n) }$ the
\emph{quenched} analogue of $\EE{R_{i}^{N}(n)}$ obtained by
conditioning on the population ancestry, \breyta{reading one
sample tree from each ancestry,  }    and averaging over
independent ancestries   \breyta{to get an estimate of  $ \EE{
\widetilde R_{i}^{N}(n) }$}.  \breyta{Our approach is
conceptually similar to  the approach of 
\cite{wakeley2012gene} and \cite{Diamantidis2024},   who  consider gene
genealogies within a fixed pedigree of diploid individuals. 
  We consider  a haploid population and  work with population ancestries,  read one sample  tree from each
ancestry,  and average over
population ancestries.      We generate a
complete sample tree forward-in-time    before  reading  statistics
off the tree.  In our construction each individual occupies a level,
and records the level occupied by its' immediate ancestor.  For
example,
\begin{equation}
\label{eq:44}
\begin{split}
& \qquad \text{levels} \\
\text{time} &  \qquad 1 \qquad \ldots \qquad \ell \qquad \ldots \qquad N \\
g& \qquad i \qquad \ldots \qquad j \qquad \ldots \qquad k
\end{split}
\end{equation}
where the  individual occupying  level $1$ resp.\ $\ell$ resp.\ $N$ at
time $g$  derives from the individual occupying level $i$ resp.\ 
$j$ resp.\ $k$ in the immediately previous generation.    If the  individual  occupying
level $\ell$ at time $g$  contributes  $x$ surviving  offspring, then
$x$ individuals at time $g+1$   will be labelled (or ``point to'')
 $\ell$.  In   this way one can record the ancestry  of the entire
population.  Given  a complete
sample tree   we read statistics (branch lengths) off the
tree, and  repeat the procedure, starting from scratch with a new
population ancestry.   }

\breyta{
\begin{equation}
\label{eq:42}
\xymatrix @R=0.1pt @C=0.1pt @M=0.1pt {
\customlabel{pp:a}{{\sf a}} \\
 \text{past} &  \vdots \\
&\circ & \ldots & \circ \\
&  \circ & \ldots & \circ \\
\text{present}\,\, & \circ & \ldots & \circ \\
}\qquad
\customlabel{pp:b}{{\sf b}}\,\xymatrix @R=0.1pt @C=0.1pt @M=0.1pt {
\vdots \\
\circ & \ldots & \circ \\
\circ & \ldots & \circ \\
\bullet & \bullet & \bullet  \\
}\qquad
\customlabel{pp:c}{{\sf c}}\,\xymatrix @R=0.1pt @C=0.1pt @M=0.1pt { 
 &  & \bullet \\
\bullet & \bullet &  \\
\vdots \\
\bullet & \bullet & \bullet  \\
}\qquad
\end{equation}
The process of sampling one quenched gene genealogy is illustrated in
\eqref{eq:42}. We start by generating a population ancestry
(\ref{eq:42}\ref{pp:a}; recalling \eqref{eq:44}  a population ancestry is a record of
the ancestral relations of all the gene copies $\circ$   in the
population),  at some arbitrary time  called ``present'' we draw a
random sample (the sampled gene copies are shown as $\bullet$ in
(\ref{eq:42}\ref{pp:b})),   once the identity of the sampled gene
copies in the given population ancestry  is known the sample tree
(\ref{eq:42}\ref{pp:c})  is fixed, and one then simply reads the
required statistic off the fixed tree;  in (\ref{eq:42}\ref{pp:c})  the
sampled gene copies and the gene copies ancestral to the sampled ones
are shown as  $\bullet$.   The quantity  $\overline \rho_{i}^{N}(n)$
\eqref{eq:functionals} is then obtained by  averaging the statistics
obtained by   repeating  the process illustrated in
\eqref{eq:42}, each time starting  with a new  population
ancestry generated independently of other population ancestries.
Thus,  $\overline \rho_{i}^{N}(n)$ obtained in this way from a finite
number of independent population ancestries 
is an  approximation  of  $\rho_{i}^{N}(n)$.  Here, we are concerned
with a haploid panmictic population  of constant size evolving
according to a specific model of sweepstakes reproduction. 
 The population  could, however,  include complex demography
such as  diploidy \citep{eldon24:_genediploid},  bottlenecks, population structure,  selfing,
dormancy (seedbanks), etc.
Clearly,  what the   limit law is  governing quenched
gene  genealogies  as defined here, and what   $\rho_{i}^{N}(n)$ converges to as
$N\to \infty$, are open questions.}

\breyta{
In
the annealed (non-quenched) approach, we update the statistics
generation by generation  without  ever generating a tree. With
$\svigi{b_{1},\ldots,b_{m}}$ denoting block sizes, starting with
$\svigi{b_{1},\ldots,b_{n}} = (1,\ldots,1)$,  the annealed approach of
sampling ``branch lengths'' $L_{i}^{N}(n)$   consists basically  of
two steps (recall we are working in discrete time),  
\begin{enumerate}
\item  $L_{b}^{N}(n) \leftarrow 1 + L_{b}^{N}(n)$    for   $b =  b_{1},\ldots, b_{m}$ 
\item sample $X_{1}, \ldots, X_{N}$ according to a specific example of
\eqref{eq:PXiJ},  merge blocks according to a
multivariate hypergeometric  using the  sampled offspring numbers, and
update block sizes
\end{enumerate}
}
 It is not {\it a priori} clear if
$\rho_{i}^{N}(n)$ and $\varrho_{i}^{N}(n)$ are different or not.     Here
we use simulations to compare $\overline \rho_{i}^{N}(n)$ and
$\overline \varrho_{i}^{N}(n)$ (recall \eqref{eq:functionals}).  We
leave a mathematical investigation of gene genealogies conditional on
the population ancestry in haploid populations evolving according to
sweepstakes reproduction to future work.
\breyta{Section \S~\ref{sec:appr-qeewid-r_inn}} contains a brief description of
the algorithm used to compute $\overline \rho_{i}^{N}(n)$.

In Figure~\ref{fig:quenchedannealedkpdbounded} and Figure~\ref{fig:QAaddDD} in Appendix~\ref{sec:estimatequenched}  we compare $\overline\rho_{i}^{N}(n) $ to   $\overline\varrho_{i}^{N}(n)$ (blue lines)  
 with $\alpha$ as shown and with the
population evolving according to Definition~\ref{def:haplrandomalpha}
(Figure~\ref{fig:qai},  Figures~\ref{fig:qaiaddDD}--\ref{fig:qaiv}) and
Definition~\ref{def:alpha-random-one}
(Figure~\ref{fig:qav},
Figures~\ref{fig:qavaddDD}--\ref{fig:qavi}). \breyta{Even though   $\overline \rho_{i}^{N}(n)$ and $\overline
\varrho_{i}^{N}(n)$   broadly agree, 
there is a visible  difference. }  In
Figure~\ref{fig:quenchedannealedkpdbounded} and
Figure~\ref{fig:QAaddDD}   we take
$\varepsilon_{N} = \overline\varepsilon_{N} = 0.1$ throughout since we
are comparing predictions of the pre-limiting model and the focus is
on seeing effects of conditioning on the population ancestry 
while keeping everything else the same between the two methods.

It has been pointed out that increasing sample size may affect the
site-frequency spectrum when the sample is from a finite population
evolving according to the Wright-Fisher model (and so cause deviations
from the one predicted by the Kingman coalescent), in particular when
the sample size exceeds the effective size  $1/c_{N}$ 
\citep{WT2003,BCS2014,Melfi2018b,Melfi2018}.  We compare $\overline
\varrho_{i}^{N}(n)$ and $\overline \varrho_{i}(n)$
 for different sample sizes to
investigate if there is any noticeable effect of increasing sample
size.  Recalling that
$ C_{\kappa}^{N} c_{N} \overset{c}{\sim} 1 $ (recall \eqref{eq:CNmap})
we conjecture it is not surprising to see little effect of increasing
sample size when   $\overline \varrho_{i}^{N}(n)$ and $\overline \varrho_{i}(n)$
 otherwise broadly agree (e.g.\
Figures~\ref{fig:strjalNe3KBC} and \ref{fig:strjalNe3KBF}).  Further
graphs are in   Figure~\ref{fig:deltanullBetaERiNaddA} in Appendix~\ref{sec:sampl-discr-trees}.

\begin{figure}[htp]
\centering
\captionsetup[subfloat]{labelfont={scriptsize,sf,md,up},textfont={scriptsize,sf}}
\subfloat[$\alpha = 0.25$,Definition~\ref{def:haplrandomalpha}]{\label{fig:qai}\includegraphics[scale=0.6]{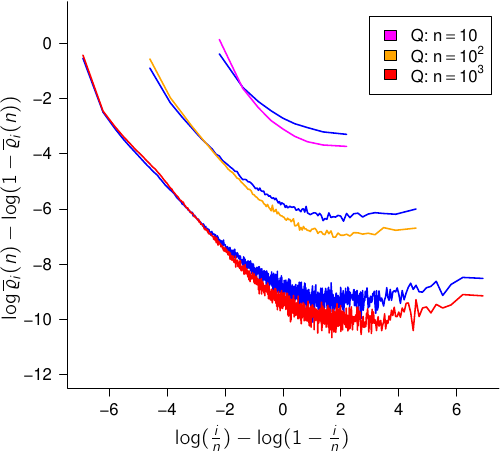}}
\subfloat[$\alpha = 0.01$,  Definition~\ref{def:alpha-random-one}]{\label{fig:qav}\includegraphics[scale=0.6]{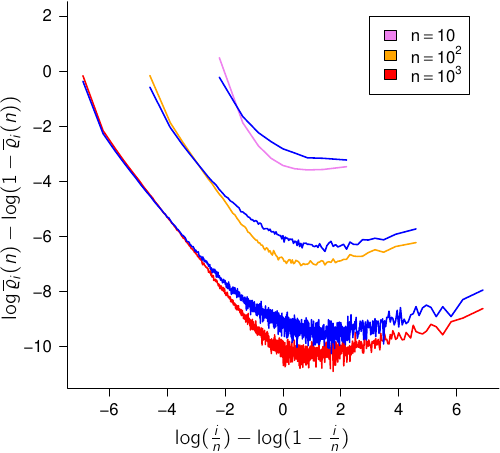}}
\caption{Quenched vs.\ annealed.  Comparing 
$\overline\rho_{i}(n)$ \breyta{(estimates of mean relative
branch lengths  when conditioning on the population ancestry)}   and $\overline \varrho_{i}^{N}(n)$
(blue lines; recall \eqref{eq:functionals})  when the population evolves according to
Definition~\ref{hschwpop} and Definition~\ref{def:haplrandomalpha} 
\breyta{(type $A$ random environment)}  a;   and
Definition~\ref{def:alpha-random-one} \breyta{(type $B$
random environment)}    b; $N=10^{3}$,
$\alpha$ as shown, $\kappa = 2$, $\uN = N^{1/\alpha}\log N$ (a) and
$\uN = N$ (b);
$\varepsilon_{N} = \overline\varepsilon_{N} = 0.1$; the approximations
are from $10^{4}$ experiments.  Sections~\ref{sec:appr-eriNn}
($\overline \varrho_{i}^{N}(n)$)  and
\ref{sec:appr-qeewid-r_inn}  ($\overline\rho_{i}^{N}(n)$) contain brief descriptions of the sampling
algorithms.  \breyta{Further graphs are in
Figure~\ref{fig:QAaddDD} in Appendix~\ref{sec:estimatequenched}.  The scale of the ordinate
 (vertical axis) may vary between graphs.  } }
\label{fig:quenchedannealedkpdbounded}
\end{figure}


\subsection{Time-changed  Beta- and Poisson-Dirichlet-coalescents}
\label{sec:lambda-coal-with}

The population size of natural populations may vary  over time. Here we
briefly  consider  fluctuations (primarily with population growth in
mind)    in the spirit of
\citep{freund2020cannings}.    Freund   studies Cannings models  with
 large enough fluctuations in population size 
 to lead to  time-changed $\Lambda$-coalescents.   See   \citep{BBMST09} for an
alternative approach with a focus on  
 recurrent bottlenecks.       Suppose $(N_{r})_{r\in \N_{0}}$ is a
sequence of population sizes where $N_{r}$ is the population size $r$
generations into the past with  $N_{0} \equiv N$, and such that
$v(s) \equiv \lim_{N\to\infty} N_{\lfloor s/c_{N} \rfloor}/N$ where $s
> 0$ and 
$c_{N}$ is the coalescence probability (Definition~\ref{def:cNhapl})  for the fixed-population size
case 
\citep[Equation~4]{freund2020cannings}.  In particular, when the
underlying population model is the one of \citep{schweinsberg03} as
giving rise to the Beta$(2-\alpha,\alpha)$-coalescent with
$1 \le \alpha < 2$, the time-change function $G(t)$ is of the form
$G(t) = \int_{0}^{t}\svigi{v(t)}^{1-\alpha}dt$
\citep[Theorem~3]{freund2020cannings}.  It follows that  gene
genealogies will be  essentially  unaffected by variations  in
population size  whenever the fluctuations   satisfy
\citep[Equation~4]{freund2020cannings} and  $\alpha$ is close to 1.

A time-changed $\delta_{0}$-Beta$(\gamma,2-\alpha,\alpha$)-coalescent
with $0 < \alpha < 2$ and whose time-change $G(t)$ is independent of
$\alpha$ follows almost  immediately from 
Theorems~\ref{thm:haplrandomalpha} and \ref{thm:hapl-alpha-random-one}
and  \citep[Lemma~4]{freund2020cannings}.
\S~\ref{sec:prooftimechangeLambda} contains a proof of
Theorem~\ref{thm:convtimechangedLambdacoal}.  

\label{sec:population-growth-}

\begin{shaded*}
\begin{thm}[A time-changed
$\delta_{0}$-Beta$(\gamma,2-\alpha,\alpha)$-coalescent]
\label{thm:convtimechangedLambdacoal}
Suppose a haploid population evolves according to
Definitions~\ref{hschwpop} and Definition~\ref{def:haplrandomalpha}
when $1 \le \alpha < 2$, and Definition~\ref{def:alpha-random-one} when $0 < \alpha
< 1$.  Suppose $\svigi{N_{r}}_{r\in \N}$ is a sequence of population
sizes with $N_{r}$ the population size $r$ generations into the past,
and where
\begin{equation}
\label{eq:7}
\begin{split}
& 0 < h_{1}(t)N \le N_{r} \le h_{2}(t)N \\
& N_{\lfloor t/c_{N} \rfloor}/N \to v(t) \quad \text{as $N\to \infty$}
\end{split}
\end{equation}
for bounded positive  functions $h_{1},h_{2},v: [0,\infty) \to
(0,\infty)$. 
Then
$\set{\xi^{n,N}(\lfloor t/c_{N}\rfloor ); t \ge 0}$ converges to
$\set{\xi^{n}(G(t)); t \ge 0}$ where  $\set{\xi^{n}(t);t \ge 0}$ is the
$\delta_{0}$-Beta$(\gamma,2-\alpha,\alpha)$-coalescent with $0 <
\gamma \le 1$ and  $0 <
\alpha < 2$, and  $G(t) = \int_{0}^{t}\svigi{v(s)}^{-1}{\rm
d}s$.   
\end{thm}
\end{shaded*}

A growing population  is sometimes modelled as increasing
exponentially in size  \citep{Donnelly1995}. 
Figure~\ref{fig:deltanullincbetaexpgrowth}  (see also Figure~\ref{fig:addbetaEE}
in Appendix~\ref{sec:further-graphs})    holds   examples of
$\overline \varrho_{i}(n)$    (recall \eqref{eq:functionals}) \breyta{  for the time-changed
$\delta_{0}$-Beta$(\gamma,2-\alpha,\alpha)$-coalescent}      when
$N_{r} = \lfloor N_{r-1}\svigi{1 - \rho c_{N} } \rfloor$, so that
$ v(s) = \lim_{N\to \infty} N_{\lfloor s/c_{N} \rfloor}/N = e^{-\rho
s}$ for $s > 0$ and  $c_{N}$ being the coalescent probability for the
fixed-size case.   Write $\tau_{m}$ for the time when
$\set{\xi^{n}\svigi{G(t)}; t \ge 0 }$ reached $m$ blocks (with
$\tau_{n} = 0$).  The
distribution of the time during which there are $m$ lineages is
\begin{displaymath}
\prob{T > \tau_{m} + t_{m} | \tau_{m}} = \exp\svigi{-\lambda_{m}\svigi{G(t_{m} + \tau_{m}) - G(\tau_{m})}}
\end{displaymath}
where $\lambda_{m}$ is the total merging rate of  $\set{\xi^{n}}$ when
$m$ blocks \citep[Remark~7, Equation~23]{freund2020cannings}.  It follows that the
time  $t_{m}$ during which  $\set{\xi^{n}\svigi{G(t)}; t \ge 0}$  has
$m$ lineages   is given by  (for $\rho > 0$)
\begin{equation}
\label{eq:8}
t_{m} =  \frac{1}{\rho}\log \svigi{1- \frac{1}{\phi_{m}}\log(1-U) }
\end{equation}
where $U$ is a random uniform on $(0,1)$,   $\phi_{m} =
(1/\rho)\lambda_{m}\exp(\rho\tau_{m})$,  and
$\tau_{n} = 0$.  As for the case when $\set{\xi^{n}}$ is the  Kingman
coalescent \citep{Donnelly1995}, the time-change acts to   increases  $\overline\varrho_{1}(n)$
   when  $\set{\xi^{n}}$ is
a $\delta_{0}$-Beta$(\gamma,2-\alpha,\alpha)$-coalescent.  

A time-changed $\delta_{0}$-Poisson-Dirichlet$(\alpha,0)$-coalescent
follows  almost immediately  from   Case~\ref{item:4} of
Theorem~\ref{thm:haplrandomalpha} and
\citep[Lemma~2]{freund2020cannings}. \S~\ref{sec:proof-thm-timechanged-poissondri}
contains a proof of Theorem~\ref{thm:time-changed-deltanullpoidr}.  

\begin{shaded*}
\begin{thm}[Time-changed
$\delta_{0}$-Poisson-Dirichlet$(\alpha,0)$-coalescent]
\label{thm:time-changed-deltanullpoidr}
Suppose a haploid  population evolves according to
Definition~\ref{hschwpop} and \breyta{type $A$ random environment} (Definition~\ref{def:haplrandomalpha})
with $0 < \alpha < 1$.  Suppose $\svigi{N_{r}}_{r\in \N}$ is as in
Theorem~\ref{thm:convtimechangedLambdacoal} (recall \eqref{eq:7}) and
such  that  $\zeta\svigi{N_{r}}/N_{r}^{1/\alpha} \to \infty$ as
$N_{r}\to \infty$ for all $r\in \N$.    
Then $\set{\xi^{n,N}\svigi{\lfloor t/c_{N} \rfloor }; t \ge 0 }$
converges to  $\set{\xi^{n}\svigi{G(t)};t \ge 0 }$ where
$\set{\xi^{n}\svigi{t};t \ge 0 }$ is a
$\delta_{0}$-Poisson-Dirichlet$(\alpha,0)$-coalescent, and  $G(t) =
\int_{0}^{t} \svigi{v(s)}^{-1}{\rm d}s$ with  $v$ as in \eqref{eq:7}. 
\end{thm}
\end{shaded*}
Figure~\ref{fig:time-changed-deltanull-poissondiri} (see also
Figure~\ref{fig:adddpdFF} in Appendix~\ref{sec:further-graphs})   contains 
examples  of $\overline \varrho_{i}(n)$ (recall
\eqref{eq:functionals})  for the  time-changed
$\delta_{0}$-Poisson-Dirichlet$(\alpha,0)$-coalescent.  As for the
time-changed   $\delta_{0}$-Beta$(\gamma,2-\alpha,\alpha)$-coalescent,
the time-change applied to the
$\delta_{0}$-Poisson-Dirichlet$(\alpha,0)$-coalescent, corresponding
to exponential population growth with rate $\rho$,   acts to  extend
$\overline\varrho_{1}(n)$.

\begin{figure}[htp]
\centering
\captionsetup[subfloat]{labelfont={scriptsize,sf,md,up},textfont={scriptsize,sf}}
\subfloat[$\alpha = 0.01$, $\gamma = 1$]{\includegraphics[scale=0.5]{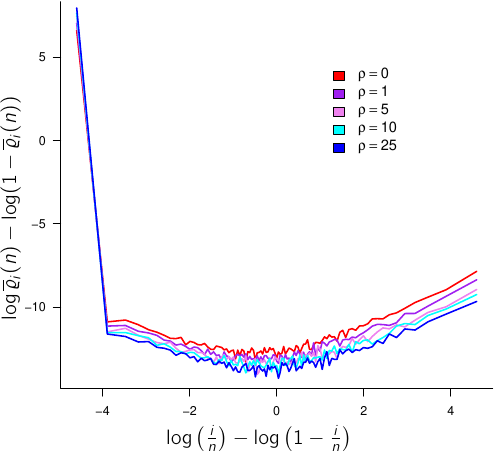}}
\subfloat[$\alpha = 1.5$, $\gamma = 0.1$]{\includegraphics[scale=0.5]{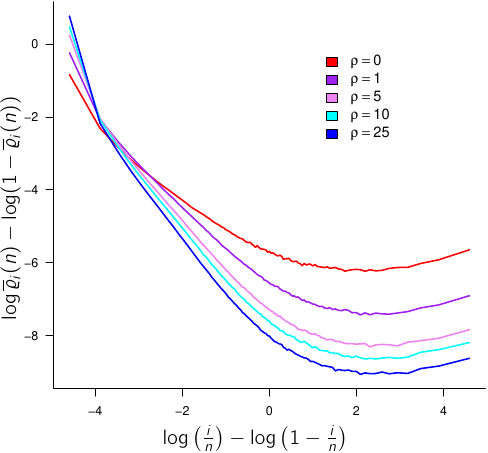}}
\caption{Approximations  $\overline \varrho_{i}(n)$
(\breyta{estimates of mean relative branch lengths,} recall
\eqref{eq:functionals})  \breyta{predicted by}    the time-changed   $\delta_{0}$-Beta$(\gamma,2-\alpha,\alpha)$-coalescent with
time-change $G(t) = \int_{0}^{t}e^{\rho s}{\rm d}s =  \one{\rho > 0}\svigi{e^{\rho
t} - 1}/\rho + \one{\rho = 0}t$  with  $\rho$ as shown,   $\kappa
= 2$, $c=1$, $n=100$.   The scale of the ordinate ($y$-axis) may vary between the
graphs; results from $10^{5}$ experiments. \breyta{Figure~\ref{fig:addbetaEE}
in Appendix~\ref{sec:further-graphs}  contains further graphs.    The scale of the ordinate
 (vertical axis) may vary between graphs.  }}
\label{fig:deltanullincbetaexpgrowth}
\end{figure}

\begin{figure}[htp]
\centering
\captionsetup[subfloat]{labelfont={scriptsize,sf,md,up},textfont={scriptsize,sf}}
\subfloat[$\alpha = 0.01$, $c=1$]{\includegraphics[scale=0.5]{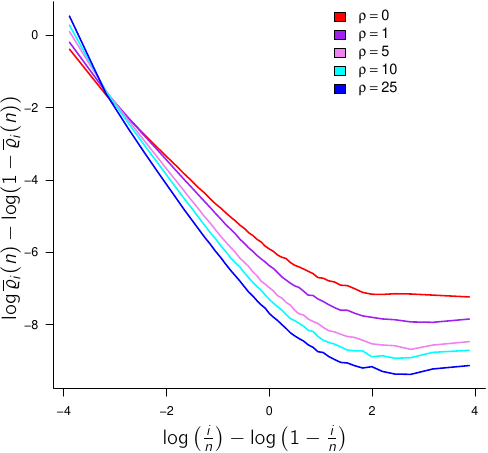}}
\subfloat[$\alpha = 0.99$, $c=100$]{\includegraphics[scale=0.5]{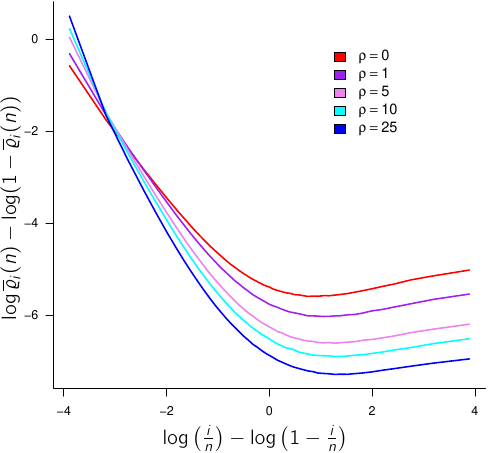}}
\caption{Approximations $\overline \varrho_{i}(n)$
(\breyta{estimates of mean relative branch lengths,}  recall
\eqref{eq:functionals})  \breyta{predicted by}   the
time-changed  $\delta_{0}$-Poisson-Dirichlet$(\alpha,0)$-coalescent
with time-change   $G(t) = \int_{0}^{t}e^{\rho s}{\rm d}s =  \one{\rho > 0}\svigi{e^{\rho
t} - 1}/\rho + \one{\rho = 0}t$  with  $\rho$ as shown,  
 $n=50$, $\kappa = 2$; 
results from   $10^{6}$ experiments.  The scale of the ordinate
 (vertical axis) may vary between graphs.  \breyta{Further graphs
are in Figure~\ref{fig:adddpdFF} in Appendix~\ref{sec:further-graphs}.}}
\label{fig:time-changed-deltanull-poissondiri}
\end{figure}

\section{Conclusion}
\label{sec:conclusion}

Our main results are {\it (i)} continuous-time coalescents derived  from
population genetic  models of sweepstakes reproduction; the coalescents
are either the Kingman coalescent, or specific families of the
$\delta_{0}$-Beta- or
$\delta_{0}$-Poisson-Dirichlet-coalescents (see Appendix~
~\ref{sec:betapoissondirichlet} for an example of a $\Xi$-coalescent
without an atom at zero); {\it (ii)} in contrast
to \citep{schweinsberg03} (recall \eqref{eq:14} in
Theorem~\ref{thm:Schwthm4}), time in (almost)  all our coalescents is measured
in units proportional to either  $N/\log N$ or $N$  generations (see
Appendix~\ref{sec:betapoissondirichlet} for an exception);  {\it (iii)}
time-changed  $\delta_{0}$-Beta- or
$\delta_{0}$-Poisson-Dirichlet-coalescents where the time-change is
independent of $\alpha$; {\it (iv)}
when $\set{\xi^{n,N}}$  is in the domain of attraction of the
$\delta_{0}$-Beta$(\gamma, 2-\alpha,\alpha)$-coalescent the  $\overline \varrho_{i}^{N}(n) $  are not matched by 
$\overline \varrho_{i}(n)$ (recall \eqref{eq:functionals}),   in particular for small values of $\alpha$; {\it
(v)} \breyta{even though  $\overline \rho_{i}^{N}(n)$ are broadly similar to 
$\overline \varrho_{i}^{N}(n) $ for the models considered here, the
difference  between them  may indicate that  $\lim_{N\to \infty}
 \rho_{i}^{(N)} \neq \lim_{N\to \infty} \varrho_{i}^{N}(n)$.}

We investigate  multiple-merger coalescents  where the driving measure
$\Xi$  is of the form
\begin{equation}
\label{eq:35}
\Xi = \begin{cases} \delta_{0} +    \one{\frac{ \uN }{N} \gneqq 0}\svigi{
\one{1 \le \alpha < 2} \Lambda_{+} +  \one{0 < \alpha < 1}\Xi_{\alpha}}  & \text{Definition~\ref{def:haplrandomalpha} holds} \\
\delta_{0} +    \one{\frac{ \uN }{N} \gneqq 0}
\Lambda_{+} & \text{Definition~\ref{def:alpha-random-one} holds}
\end{cases}
\end{equation}
where $\Lambda_{+}$ is a finite measure on $(0,1]$ without an atom at
zero.  In our formulation  $\Xi$  always retains
an atom at zero.  The coalescents derived here are based on a
construction where most of the time small (relative to the population
size) families occurs, but occasionally the environment favours the
generation of (at least one) large family.  Sweepstakes reproduction  
 is incorporated as in Definition~\ref{def:haplrandomalpha}
 \breyta{(type $A$ random environment)} 
and Definition~\ref{def:alpha-random-one} \breyta{(type $B$
random environment)}   with an upper bound on the
number of potential offspring any one individual may produce.
 \breyta{Sweepstakes reproduction,  for example  in broadcast
spawning marine organisms,  may occur through     chance matching of
favorable environmental conditions and  reproduction
\citep{HP11,Li1998,H82,B94}.   Random  environments similar to the
ones  in Definitions~\ref{def:haplrandomalpha}   and ~\ref{def:alpha-random-one}   
seem  natural  for  modeling  sweepstakes reproduction.     }

 In our formulation (recall \eqref{eq:PXiJ})   the upper bound $\uN$ on the number of
potential offspring  determines if the limiting coalescent admits 
multiple-mergers (recall~\eqref{eq:35}).    For comparison,
$\Lambda = \one{\alpha \ge 2}\delta_{0} + \one{1 \le \alpha <
2}\Lambda_{+}$ absent a bound on the number of potential offspring
\citep{schweinsberg03}. Moreover, we show that restricting $\uN$ to
the order of the population size is sufficient to obtain a
multiple-merger coalescent. In contrast, convergence to the original
(complete) Beta$(2-\alpha,\alpha)$-coalescent requires the (slightly)
stronger assumption of $\uN/N \to \infty$.  Applying an upper bound
appears reasonable for at least some highly fecund populations.
Considering broadcast spawners, we see potential offspring as
representing fertilised eggs (or some later stage in the development
of the organisms before maturity is reached).   It is plausible that the number of fertilised eggs
produced by an arbitrary individual for (at least some) broadcast
spawners may be at most a fraction of the population size.  Our
numerical results show that an upper bound on the number of potential
offspring markedly affects the predicted site-frequency spectrum (see
Appendix~\ref{sec:relat-expect-sfs} for examples comparing
$\varphi_{i}(n)$ from \eqref{eq:varphi} between
$\delta_{0}$-Beta$(\gamma,2-\alpha,\alpha)$-coalescents).  Thus, even
though our extensions (\eqref{eq:PXiJ},
Theorems~\ref{thm:haplrandomalpha} and  \ref{thm:hapl-alpha-random-one})
 of \eqref{eq:4} are relatively straightforward,
the implications  for inference  of our results   are significant.

Theorems~\ref{thm:haplrandomalpha} and \ref{thm:hapl-alpha-random-one}
combine to yield a continuous-time
$\delta_{0}$-Beta$(\gamma, 2-\alpha,\alpha)$-coalescent where
$0 < \gamma \le 1$ and $0 < \alpha < 2$, with the understanding that
the population model when $1\le \alpha < 2$  is different from the one
when $0 < \alpha < 1$.  When $1 \le \alpha < 2$ the population evolves
as in Definition~\ref{def:haplrandomalpha} \breyta{(type $A$
random environment)}, and according to
Definition~\ref{def:alpha-random-one} \breyta{(type $B$
random environment)}   when $0 < \alpha < 1$.  Thus, we
have extended the Beta-coalescent \citep{schweinsberg03} to include an 
atom at zero,  a truncation point ($\gamma$; see also
\cite{Eldon2026biorxiv}),   and with the range
of $\alpha$ extended from $[1,2)$ to $(0,2)$.  However, the number of
parameters involved has now increased from one $(\alpha)$ to four
$(c,\gamma, \alpha,\kappa)$, with the mean $m_{\infty}$ being a
function of $\kappa$ (recall \eqref{eq:57}).   Moreover, when in the domain of attraction of
the $\delta_{0}$-Beta$(\gamma,2-\alpha,\alpha)$-coalescent
$\EE{X_{1}} < \infty$ regardless of the value of $\alpha$; under
\eqref{eq:4} $\EE{X_{1}} = \infty$  when $ 0< \alpha \le 1$.  For the
$\delta_{0}$-Beta$(\gamma, 2-\alpha,\alpha)$ coalescent the unit of
time is proportional to $N/\log N$ generations when $\kappa = 2$, and
to $N$ generations when $\kappa > 2$.    For
comparison, time in the original complete Beta-coalescent is measured
as shown in \eqref{eq:14} \citep[Lemma~16; Lemma~13]{schweinsberg03}.
Effective population size $(N_{e})$ is often taken as $1/c_{N}$, and
so would be proportional to at least $N/\log N$ for our models. 
One can ask what assumptions on the population size $(N)$ and the
mutation rate $(\mu)$ are required to recover the mutations in a
sample of DNA sequences used to obtain an estimate of $\alpha$.  Using
the $N^{\alpha - 1}$ timescale \breyta{($1 < \alpha < 2$)},  it can be argued  that recovering
observed mutations requires strong assumptions on $N$ (and/or $\mu$ if
unknown) when  the estimate of $\alpha$ is near 1 \citep{Eldon2026biorxiv}.   Moreover, our
results  show that sweepstakes reproduction  does not
necessarily imply a small $N_{e}/N$ ratio \citep{JFB:JFB13143,A04}.

We obtain continuous-time
$\delta_{0}$-Poisson-Dirichlet$(\alpha,0)$-coalescents with time
measured in units proportional to at least $N/\log N$
generations. \breyta{Sagitov \cite[\S~3]{S03} discusses an    
example of a Poisson-Dirichlet-coalescent deduced from a
 compound multinomial distribution.}
  Combining Definition \ref{hschwpop} and
\eqref{eq:4} when $0 < \alpha < 1$ yields a
Poisson-Dirichlet$(\alpha,0)$-coalescent without rescaling time
($c_{N}\overset {c}{\sim} 1$) \citep{schweinsberg03}. However,
convergence to the
$\delta_{0}$-Poisson-Dirichlet$(\alpha,0)$-coalescent requires that
$\uN/N^{1/\alpha}\to \infty$, imposing a strong assumption on the
distribution of potential offspring \breyta{(recall
$0 < \alpha < 1$)}.  We leave it for later to investigate the
coalescents that may arise when $\uN/N^{1/\alpha}$
\breyta{converges as $N\to \infty$  to 
some constant $K > 0$;}  we
conjecture (see Remark~\ref{rm:incompletePD}) that the resulting
coalescent may be a form of an `incomplete Poisson-Dirichlet'
coalescent, similar to the case $\uN/N \to K$ required for convergence
to an incomplete Beta-coalescent.  Other scalings of $\uN$ may also
lead to new non-trivial coalescents.

The introduction of the coalescent, a sample-based approach for
investigating evolutionary histories of natural populations moved
population genetics forward. The coalescent   enables efficient generation of gene
genealogies, in particular when applying recent advances in algorithm
design \citep{Baumdicker2021,Kelleher2016}.  However, a
coalescent-based inference is only helpful when the coalescent
correctly approximates the trees generated by the corresponding
ancestral process.  Here, we check this approximation by comparing
functionals (recall~\eqref{eq:54} and~\eqref{eq:functionals}) of $\set{\xi^{n}}$ and
$\set{\xi^{n,N}}$ when the population is haploid and panmictic of
constant size and evolving under sweepstakes reproduction.  Given the
poor agreement between $\overline \varrho_{i}(n)$ and
$\overline \varrho_{i}^{N}(n)$ in some cases, in particular when the
law for the number of potential offspring is strongly skewed, the
usefulness of coalescent-based inference is in doubt in this case.
 \breyta{Multiple-merger coalescents
have been applied to various   broadcast-spawning marine organisms
\citep{AH2015,Arnasonsweepstakes2022,Vendrami2021,Niwa2016,https://doi.org/10.1111/mec.16774},
 outbreaks of  tuberculosis \citep{10.1093/molbev/msaa179},  and even to
  the propagation of cancer cells\citep{10.1098/rsos.171060}.
  Moreover, multiple-merger coalescents may be suitable for
  investigating  plants and (crop-infesting)  fungi that reproduce by distributing
 huge numbers of seeds or  spores  \citep{Minadakis2025,Jigisha2025}.
 Thus, the models studied here, and the resulting  
 multiple-merger coalescents may be relevant for quite
 a broad range of study systems.  Moreover,   }

One might want to extend the models applied here to diploid
populations (where gene copies always occur in pairs regardless of the
reproduction mechanism),    and to include complex demography such as
recurrent bottlenecks \citep{BBMST09}.   Multiple-merger coalescents
derived from population models of sweepstakes reproduction  in
diploid populations tend to be $\Xi$-coalescents
\citep{MS03,BLS15,BBE13}.  Site-frequency spectra predicted by
$\Xi$-coalescents can  be different from the ones predicted by
$\Lambda$-coalescents \citep{Blath2016,BBE2013a}.  As we have seen
(Theorem~\ref{thm:hapl-alpha-random-one}), coalescents derived for
haploid panmictic populations of constant size and evolving according
to sweepstakes   may also be continuous-time
$\Xi$-coalescents (the
$\delta_{0}$-Poisson-Dirichlet$(\alpha,0)$-coalescent).  \breyta{
In any extensions of the models considered here, one should compare  $\overline \varrho _{i}^{N}(n)$ to $\overline \varrho_{i}(n)$
, and to $\overline \rho_{i}^{N}(n)$.  Should
any    observed difference 
between $\overline \varrho_{i}^{N}(n)$ and  $\overline
\rho_{i}^{N}(n)$ not be an artifact of the finite   population
size used in  simulations such that  $\lim_{N\to\infty}
\rho_{i}^{N}(n) \neq \lim_{N\to \infty}\varrho_{i}^{N}(n) $  then we
may have to rethink gene genealogies from scratch. }

\section{Proofs}\label{sec:proofs-1}

\breyta{
This section contains  proofs of Theorem
~\ref{thm:haplrandomalpha} in
\S~\ref{sec:proof-thm-randomallhaploidnewmodel}, of
Theorem~\ref{thm:hapl-alpha-random-one} in
\S~\ref{sec:proof-thm-random-alpha-one}, Theorem~\ref{thm:convtimechangedLambdacoal} in
\S~\ref{sec:prooftimechangeLambda}, and of
Theorem~\ref{thm:time-changed-deltanullpoidr} in
\S~\ref{sec:proof-thm-timechanged-poissondri}.  First, 
in  \S~\ref{sec:verifying-PDtransition-probability},  we verify~\eqref{eq:34}.}
\subsection{Verifying~\eqref{eq:34}}\label{sec:verifying-PDtransition-probability}

In this section we check that~\eqref{eq:34} holds.  For  $x,a\in \mathds{R}$
and   $m\in \{0,1,2, \ldots \}$  define  ${[x]}_{0;a} = 1$, and  
\begin{displaymath}
{[x]}_{m;a} = x(x+a)\cdots(x+(m-1)a), \quad m \ge 1
\end{displaymath}
so that ${[\alpha]}_{r-1;\alpha} = \alpha^{r-1}(r-1)!$ and
${[1]}_{m-1;1} = (m-1)!$.  Then, with $b \ge 2$, $0 < \alpha < 1$, and
$k_{1} + \cdots + k_{r} = b$ where $\{k_{1}, \ldots, k_{r} \}$ is a
partition of $b$, the exchangeable probability function is 
\begin{displaymath}
\breyta{p(k_{1},\ldots,k_{r}) = }  \frac{[\alpha]_{r-1;\alpha}}{[1]_{b-1;1} }\prod_{i=1}^{r}[1-\alpha]_{k_{i}-1;1} =    \frac{\alpha^{r-1}(r-1)! }{(b-1)!}\prod_{i=1}^{r}[1-\alpha]_{k_{i}-1;1}
\end{displaymath}
 \citep[Proposition~9,  Equation~16]{pitman1995exchangeable},   where,
 recalling \eqref{eq:1},   
\begin{displaymath}
[1-\alpha]_{k_{i}-1;1} = (1-\alpha)(2-\alpha)\cdots (k_{i}-\alpha - 1) = (k_{i}-\alpha - 1)_{k_{i}-1}
\end{displaymath}
Following \citep{schweinsberg03} write 
\begin{displaymath}
  p_{s}(k_{1}, \ldots, k_{r}) = p(k_{1},\ldots, k_{r}, \underset{s
\text{ times}}{\underbrace{1, \ldots, 1}}) 
\end{displaymath}
where   
\begin{displaymath}
   p_{b;k_{1}, \ldots, k_{r};0} =   p_{0}(k_{1}, \ldots, k_{r}) =  \frac{\alpha^{r-1}(r-1)!}{(b-1)!}\prod_{i=1}^{r}[1-\alpha]_{k_{i}-1;1}
\end{displaymath}
and 
\begin{displaymath}
p_{0}(k_{1}, \ldots, k_{j-1},k_{j}+1, k_{j+1},\ldots, k_{r})  = \frac{ \alpha^{r-1}(r-1)!  }{b!} \breyta{[1-\alpha]_{k_{j};1}}   \prod_{ \breyta{ i\in [r]\setminus\set{j}} }[1-\alpha]_{k_{i}-1;1}
\end{displaymath}
Then, by \citep[page 137]{schweinsberg03},
\begin{equation}
\label{eq:29}
\begin{split}
p_{s+1}(k_{1},\ldots, k_{r} ) &  = p_{s}(k_{1}, \ldots, k_{r}) - \sum_{j=1}^{r}  p_{s}(k_{1}, \ldots, k_{j-1},k_{j}+1, k_{j+1}, \ldots, k_{r}) \\
& - sp_{s-1}(k_{1}, \ldots, k_{r},2 )
\end{split}
\end{equation}
By \eqref{eq:29} with $s = 0$, noting that  $[1-\alpha]_{k;1} =
[1-\alpha]_{k-1;1}(k-\alpha)$, 
\begin{displaymath}
\begin{split}
p_{1}(k_{1}, \ldots, k_{r}) & = p_{0}(k_{1},\ldots, k_{r}) - \sum_{j=1}^{r}p_{0}(k_{1}, \ldots, k_{j-1},k_{j}+1, k_{j+1}, \ldots , k_{r}) \\
& =  \frac{\alpha^{r-1}(r-1)! }{(b-1)!}\left( \prod_{i=1}^{r}[1-\alpha]_{k_{i}-1;1}  - \frac{1}{b} \sum_{j =1}^{r} [1-\alpha]_{k_{j};1} \prod_{\substack{i=1 \\ i \neq j}}^{r}[1-\alpha]_{k_{i}-1;1} \right) \\
& =  \frac{\alpha^{r-1}(r-1)! }{(b-1)!}\left( \prod_{i=1}^{r}[1-\alpha]_{k_{i}-1;1} - \frac{1}{b} \sum_{j=1}^{r}(k_{j}-\alpha)\prod_{i=1}^{r}[1-\alpha]_{k_{i}-1;1} \right) \\
& =   \frac{\alpha^{r-1}(r-1)! }{(b-1)!}\prod_{i=1}^{r}[1-\alpha]_{k_{i}-1;1} \left( 1-  \frac{1}{b}\sum_{j=1}^{r}(k_{j}-\alpha)  \right)
\end{split}
\end{displaymath}
We then have,    recalling $k_{1} + \cdots  + k_{r} = b$,  
\begin{displaymath}
\begin{split}
p_{1}(k_{1}, \ldots, k_{r}) & = \frac{\alpha^{r-1}(r-1)!}{(b-1)!}\prod_{i=1}^{r}[1-\alpha]_{k_{i}-1;1}\left( 1- \sum_{j=1}^{r}\frac{k_{j} - \alpha}{b}  \right) \\
& =   \frac{\alpha^{r-1}(r-1)!}{(b-1)!}\prod_{i=1}^{r}[1-\alpha]_{k_{i}-1;1} \left( 1  - \frac{\sum_{j=1}^{r}k_{j}}{b}  + \frac{r\alpha}{b} \right) \\
& =  \frac{\alpha^{r} r!}{b!}\prod_{i=1}^{r}[1-\alpha]_{k_{i}-1;1} =    \frac{\alpha^{r} r!}{b!}\prod_{i=1}^{r}(k_{i} - \alpha - 1)_{k_{i} - 1}
\end{split}
\end{displaymath}
Similarly one uses~\eqref{eq:29} to check~\eqref{eq:34} for $s=1$, hence to~\eqref{eq:34} by induction over $s$.  

\subsection{Proof of Theorem~\ref{thm:haplrandomalpha}}\label{sec:proof-thm-randomallhaploidnewmodel}

In this section we give a proof of Theorem~\ref{thm:haplrandomalpha}.
Therefore, we consider a haploid population evolving according to
Definitions~\ref{hschwpop} and~\ref{def:haplrandomalpha}
\breyta{(type $A$ random environment)}    with law
~\eqref{eq:PXiJ} on the number of potential offspring produced by each
individual.    \breyta{To briefly  orient the reader
(the proofs can always be skipped on first reading to get an overview
of the lemmas leading to  the theorem),
Proposition~\ref{pr:GHbounds} and Lemmas~\ref{lm:boundssumrka} and
~\ref{lm:ineqs} collect useful results for deriving the limits of
$C_{\kappa}^{N}c_{N}$ (recall $C_{\kappa}^{N}$ from~\eqref{eq:CNmap})
in Lemma~\ref{lm:cNhaplrandomall}.  Lemma~\ref{lm:convKingmanrall}
proves  Case~\ref{item:1} of
Theorem~\ref{thm:haplrandomalpha}  by checking the condition  given
in~\eqref{condi}.        The proof of Case~\ref{item:2} of 
Theorem~\ref{thm:haplrandomalpha} involves  checking the conditions of
 Proposition~\ref{prp:conditionsconvergenceLambda},
 Lemma~\ref{lm:tailprobability} identifies the $\Lambda_{+}$-measure
 (Condition~\eqref{eq:existenceLambda} of  Proposition~\ref{prp:conditionsconvergenceLambda}), 
 and   Lemma~\ref{lm:simmergersvanisharandall}  checks  that
 simultaneous  mergers vanish as $N\to \infty$
 (Condition~\ref{eq:simtozero} of
 Proposition~\ref{prp:conditionsconvergenceLambda}).        The proof of
 Case~\ref{item:4} of  Theorem~\ref{thm:haplrandomalpha}  follows the
 proof of \citep[Theorem~4d]{schweinsberg03}.
 Lemmas~\ref{lm:randomalphamNbounded}  $(m_{\infty} < \infty)$  and~\ref{lm:almostsureconv}
  ($S_{N}/(Nm_{N})\to 1$ almost surely) are ancillary  results
  required for  checking the convergence of  $C_{\kappa}^{N}c_{N}$.
  Note in particular that, since $m_{\infty} < \infty$ in our
  framework,   we can use the same method
  for  $\alpha = 1$ as for  $1 < \alpha < 2$.  }  

\begin{propn}[\citep{Eldon2026biorxiv}; Lemma~7.2]\label{pr:GHbounds}
Suppose $G$ and $H$ are positive  functions on $[1,\infty)$, 
$H$ \breyta{be} monotone  increasing,   $G$ monotone
decreasing,  and $\int HG^{\prime}$ exists.  
Then for  $\ell,m \in \IN$ 
\begin{displaymath}
-\int_{\ell}^{m+1} H(x-1)G^{\prime}(x)dx \le  \sum_{k=\ell}^{m} H(k)(G(k) - G(k+1)) \le -\int_{\ell}^{m+1} H(x)G^{\prime}(x)dx
\end{displaymath}
\end{propn}

We will  use Proposition~\ref{pr:GHbounds} to obtain bounds on
$\sum_{k=\ell}^{m} k(k-1) {(k+M)}^{-2}\svigi{k^{-a} - {(1+k)}^{-a}}$; when
we do so we intend to replace $M$, $\ell$, and $m$ as later specified
to identify the limit  of $C_{\kappa}^{N}c_{N}$ as $N\to \infty$.

\begin{lemma}[Bounds on $\sum_{k}k(k-1){(k +M)}^{-2}(k^{-a} - {(k+1)}^{-a})$]\label{lm:boundssumrka} 
For $k\in \IN$, with $M \gg 1$,  $a >  0$, $1 \le \ell \ll m $ all  fixed; write $r_{k}(a) = k(k-1){(k+M)}^{-2}(k^{-a} - {(k+1)}^{-a})$.
\begin{enumerate}
\item~\label{item:15} When $0 < a < 1$ we have
\begin{equation}
\label{eq:boundsrka}
\begin{split}
 &  \frac{a}{(M-1)^{a}} \int_{\frac{\ell}{\ell + M - 1 } }^{\frac{m+1 }{M + m  } } u^{1-a}(1-u)^{a-1}du  + O(\tfrac{1}{M^{2}})  \le  \sum_{k=\ell}^{m} r_{k}(a) \\
 & \le \frac{a}{M^{a}} \int_{\frac{\ell}{\ell + M } }^{\frac{m+1 }{M + m + 1 } } u^{1-a}(1-u)^{a-1}du  + O(\tfrac{1}{M^{2}} )
\end{split}
\end{equation}
\item \label{item:10}  When $a = 1$ we have
\begin{displaymath}
\begin{split}
& \frac{m-\ell + 1}{(M+\ell - 1)(M+m) } + O\left( \frac{ \log M}{M^{2}} \right) \le    \sum_{k=\ell}^{m} r_{k}(a) \\
&  \le \frac{m-\ell + 1}{(M+\ell)(M+m+1) } +  O\left( \frac{ \log M}{M^{2}} \right)
\end{split}
\end{displaymath}
\item \label{item:12} When $1 < a < 2$ the same  bounds as  in~\eqref{eq:boundsrka} for Case~\ref{item:15}  hold.
\item \label{item:13}  When $a = 2$
\begin{displaymath}
   \sum_{k=\ell}^{m} r_{k}(a) =     2M^{-2} \log M  +  O(M^{-2})
\end{displaymath}
\item~\label{item:14}  When $2 < a $
\begin{equation}
\label{eq:kappatwoplus}
\begin{split}
& \frac{a}{a-2}\frac{\ell^{2-a} }{(\ell + M - 1)^{2} } - \frac{3a}{a-1 }\frac{\ell^{1-a} }{(\ell + M - 1)^{2} } +  \frac{2\ell^{-a} }{(\ell + M-1)^{2}} \\
& + O(m^{1-a}(m + M)^{-2}) + O(M^{-3}) \le  \sum_{k=\ell}^{m}r_{k}(a) \\
& \le \frac{a}{a - 2}\frac{ \ell^{2-a}}{ (\ell + M)^{2} }  -  \frac{a}{a-1}\frac{\ell^{1-a}}{(\ell + M)^{2} }   + O(m^{2-a}(m+ 1 + M)^{-2})   + O(M^{-a})
\end{split}
\end{equation}
\end{enumerate}
\end{lemma}
\begin{proof}[Proof of Lemma~\ref{lm:boundssumrka}]
Take   $ H(x) =  x(x-1)(x+M)^{-2}$ and $G(x) = x^{-a}$ in Proposition~\ref{pr:GHbounds}   to obtain 
\begin{multline}
\label{eq:boundsGSints}
 \int_{\ell}^{m+1} \frac{(x-1)(x-2)}{(x+M-1)^{2} }\frac{a }{x^{a+1}} dx \le  \sum_{k=\ell}^{m} \frac{k(k-1) }{(k + M)^{2} } \left(\frac{1}{k^{\alpha}} - \frac{1}{(1 + k)^{\alpha}} \right) \\   \le  \int_{\ell}^{m+1} \frac{x(x-1)}{(x+M)^{2} }\frac{a }{x^{a+1}} dx
\end{multline}
Now the task  is  to evaluate the integrals in \eqref{eq:boundsGSints}. One can write
\begin{multline}
\label{eq:32}
\int \frac{(x-1)(x-2)}{(x+M-1)^{2}}\frac{a}{x^{1+a}}dx   =  \int \frac{ax^{1-a} }{(x+M-1)^{2}}dx -  \int \frac{3a }{x^{a}(x+M-1)^{2} }dx \\ 
    +  \int \frac{2a}{x^{1+a}(x+M-1)^{2} }dx  
  \end{multline}
  and 
  \begin{align}
\label{eq:31}
\int \frac{x(x-1) }{(x+M)^{2}}\frac{a}{x^{1+a}}dx &  =  \int \frac{ax^{1-a} }{(x+M)^{2}}dx  -  \int \frac{a }{x^{a}(x+M)^{2}}dx 
\end{align}
and evaluate each integral on the right in \eqref{eq:32} and
\eqref{eq:31} separately using standard integration techniques.

When $a = 1$  using  \eqref{eq:31} and integration by partial fractions we see 
\begin{displaymath}
\begin{split}
\int_{\ell}^{m+1}\frac{x(x-1)}{(x+M)^{2}x^{2}}dx &  =  \frac{m-\ell + 1}{(M + \ell)(M + m + 1) } \\
& -   \frac{1}{M^{2}} \left( \log\frac{ m + 1  }{M + m + 1 } + \log\frac{\ell + M }{\ell}  + M\frac{\ell - m - 1  }{(m+ M + 1)(M + \ell ) } \right)
\end{split}
\end{displaymath}
The lower bound in Case~\ref{item:10} is obtained in the same way using  \eqref{eq:32}.

When $0 < a < 1$ (Case~\ref{item:15}) or   $1 < a < 2$  (Case~\ref{item:12})  the substitution
$y = M/(x + M)$ (so that $x = My^{-1} - M$ and $dx = -My^{-2}dy$) gives 
\begin{displaymath}
a\int_{\ell}^{m+1} \frac{x^{1-a}}{(x+M)^{2} }dx =  \frac{a}{M^{a}} \int_{\frac{M }{M + m + 1 } }^{\frac{M}{M+\ell} }(1-y)^{1-a}y^{a-1}dy =   \frac{a}{M^{a}} \int_{\frac{\ell}{\ell + M } }^{ \frac{m+1 }{M + m + 1 } }u^{1-a}(1-u)^{a-1}du  
\end{displaymath}
The same substitution and integration by parts gives
\begin{displaymath}
a\int_{\ell}^{m + 1} \frac{x^{-a}}{(x+M)^{2} }dx = \left[ \frac{a}{1-a}\frac{x^{1-a}}{(x+M)^{2} }   \right]_{\ell}^{m+1} + \frac{2a}{1-a} \int_{\ell}^{m+1} \frac{x^{1-a} }{(x+M)^{3} }dx \in O(M^{-2}) 
\end{displaymath}
A similar calculation for the lower bound in \eqref{eq:boundsGSints} gives Case~\ref{item:15}  and   \eqref{eq:boundsrka}  for $1 < a < 2$ (Case~\ref{item:12}). 

When $a = 2$ (Case~\ref{item:13})   \eqref{eq:31} and  integration by partial fractions gives
\begin{displaymath}
\begin{split}
\int_{\ell}^{m+1}\frac{x(x-1)}{{(x+M)}^{2}}\frac{2}{x^{3}}dx & =  {2M^{-2}\log(1 + M/\ell)} + O\svigi{M^{-2}} + O\svigi{M^{-3}\log M}
\end{split}
\end{displaymath}

When $a > 2$ (Case~\ref{item:14})  integration by parts and the substitution $y = M/(x + M)$
as in the case $1 < a < 2$ give the result.   When $2 < a < 3$  integration by parts gives 
\begin{displaymath}
\begin{split}
\int_{\ell}^{m+1} \frac{ax^{1-a}}{(x+M)^{2}}dx & =  \frac{a}{2-a}\left[ \frac{x^{2-a}}{(x+M)^{2}}  \right]_{\ell}^{m+1} - \frac{2a}{a-2}\int_{\ell}^{m+1}\frac{x^{2-a}}{(x+M)^{3}}dx  \\
&  = \frac{a}{a-2}\frac{\ell^{2-a}}{(\ell+M)^{2}  } +  O\svigi{ \frac{m^{2-a}}{ (m+M)^{2} } } + O(M^{-a}) 
\end{split}
\end{displaymath}
using the substitution $y = M/(x+M)$ for the error term $O(M^{-a})$.
When $a \ge 3$ we simply repeat the calculation on
$\int x^{2-a}(x+M)^{-3}dx$ to see that we retain the same leading
term.   Considering  $\int x^{-a}(x+M)^{-2}dx$  we again  apply integration by parts to see
\begin{displaymath}
\begin{split}
\int_{\ell}^{m + 1} \frac{a}{x^{a}(x+M)^{2} }dx &  =  \frac{a}{1-a}\left[ \frac{x^{1-a}}{(x+M)^{2}}  \right]_{\ell}^{m +1}  - \frac{2a}{a-1} \int_{\ell}^{m + 1}\frac{x^{1-a} }{(x+M)^{3}}dx \\
& =  \frac{a}{a-1}\frac{\ell^{1-a}}{(\ell+M)^{2}} + O\svigi{\frac{ m^{1-a}}{ (m + M)^{2}  }  } + O\svigi{M^{-3} }
\end{split}
\end{displaymath}
The remaining terms of Case~\ref{item:14} are identified in a similar way (here we omit the details).
\end{proof}

\begin{lemma}[Bounds on $k^{-a} - (k+1)^{-a}$]
\label{lm:ineqs}
  For  any  $0 < a \le  1$ and $k\in \IN$  we have 
\begin{equation}
\label{eq:20}
    \frac{a}{(k+1)^{1+a} } \le   \frac{1}{k^a} - \frac{1}{(1 + k)^a} \le \frac{1}{k^{1 + a}}
\end{equation}
When $a \ge  1$ and $k\ge 2$  we have 
\begin{equation}
\label{eq:21}
     \frac{1}{k^{1 + a}} \le  \frac{1}{(k-1)^a} - \frac{1}{k^a}  \le  \frac{a}{(k-1)^{a + 1}}
\end{equation}
\end{lemma}
\begin{proof}[Proof of Lemma~\ref{lm:ineqs}]
 Suppose $0 < a \le  1$.   Then
\begin{displaymath}
1 - \frac{1}{k} =  \frac{k-1}{k} \le  \frac{k}{k+1} \le  \left( \frac{k}{k+1} \right)^{a}
\end{displaymath}
so that  $k^{-a} -  k^{-1-a} \le (k+1)^{-a} $.  
Rearranging terms   gives the upper   bound in   \eqref{eq:20}.

Considering the lower bound in  \eqref{eq:21} we see, for any $a \ge  1$, 
\begin{displaymath}
1 + \frac{1}{k} = \frac{k+1}{k} \le \frac{k }{k-1} \le \left( \frac{k}{k-1} \right)^{a} 
\end{displaymath}
Then $k^{-a} + k^{-1-a} \le (k-1)^{-a}$ and rearranging terms gives the bound.

 For the upper bound in   \eqref{eq:21} we recall the inequality
\begin{displaymath}
x^{r} - y^{r} < r x^{r - 1}(x - y)
\end{displaymath}
for any $x > y > 0$ and any rational  $r > 1$   \citep[\S 74, Equation~9]{hardy02}. Suppose $x = k+1$ and $y = k$.  We see
\begin{displaymath}
\frac{1}{k^{r}} - \frac{1}{(1+k)^{r}}  =  \frac{ (k+1)^{r} - k^{r}}{k^{r}(k + 1)^{r} } < \frac{r(k+1)^{r-1} }{k^{r}(k+1)^{r} }  = \frac{r}{k^{r}(1+k)} < \frac{r}{k^{r+1}}
\end{displaymath}
Since for any fixed $x >  0$ the function $r \mapsto x^{-r} $ is
continuous in $r$, and the rationals are dense in $\IR$, the upper
bound in \eqref{eq:21} follows.

The lower bound in \eqref{eq:20} follows by similar  calculations,
using the inequality $sx^{s-1}(x-y) < x^{s} - y^{s}$ for 
$0 < s < 1$ rational  and any $x,y \in \IR$ with   $x > y > 0$ \citep[\S 74, Equation~10,  pp.\ 144]{hardy02}.
\end{proof}
The following lemma gives a condition on
\breyta{$\svigi{\varepsilon_{N}}_{N}$ (recall
$\varepsilon_{N}$ is the probability of event $E$ in   type $A$ random
environment defined in  
Definition~\ref{def:haplrandomalpha})}    for the
mean  \breyta{$\EE{X_{1}} =  m_{N}$}  from \eqref{eq:5} to be finite for all $N$, and for
\breyta{ $m_\infty  =  \lim_{N\to \infty}m_{N}$}    to be finite.
\begin{lemma}[Finite $m_{\infty}$]
\label{lm:randomalphamNbounded}
Under the conditions of Theorem~\ref{thm:haplrandomalpha},      if
\begin{displaymath}
\varepsilon_{N} \in   
\begin{cases}
 O(\uN^{\alpha -1 }) & \text{when $0 < \alpha < 1$} \\
  O( 1/ \log \uN )    & \text{when $ \alpha = 1$} \\
 O(1) &   \text{when $\alpha > 1$} \\
\end{cases}
\end{displaymath}
then, with $X$ the  random number of potential offspring of an arbitrary individual, 
\begin{displaymath}
  \limsup_{N\to \infty}\EE{X } < \infty. 
\end{displaymath}
\end{lemma}%
\begin{proof}[Proof of Lemma~\ref{lm:randomalphamNbounded}]
Recalling event $E$ from Definition~\ref{def:haplrandomalpha} it is
straightforward to check that
$\limsup_{N\to \infty}\EE{X |E } < \infty$ whenever $\alpha > 1$.
Therefore, it is sufficient to consider $\EE{X |E}$ for the case
$0 < \alpha \le 1$.
First with  $0 < \alpha < 1$ and  using the upper bound in \eqref{eq:20} in    Lemma~\ref{lm:ineqs} and recalling $\overline{f_{\alpha}}$ from \eqref{eq:isfg}, 
\begin{equation}
\label{eq:43}
\begin{split}
\EE{ X | E }   \le &   \overline{f_{\alpha}}\sum_{k=1}^{\uN} k\left( \frac{1}{k^\alpha} - \frac{1}{(1+k)^\alpha}  \right)
   \le  \overline{f_{\alpha}} \sum_{k=1}^{\uN} {k^{-\alpha}}  \le \overline{f_{\alpha}} + \overline{f_{\alpha}} \int_{1}^{\uN}{x^{-\alpha}}dx \\
 = &   \one{ 0 < \alpha < 1}\frac{\overline{f_{\alpha}}}{1-\alpha} (\uN ^{1-\alpha} - \alpha)  + \one{\alpha = 1}\overline{f_{\alpha}} (1 +   \log \uN ) \\
\end{split}
\end{equation}
With  event $E^{\sf c}$ as in   Definition~\ref{def:haplrandomalpha} and   $\EE{X}  =  \EE{X | E }\varepsilon_{N} + \EE{X |
E^{\sf c} }(1-\varepsilon_{N}) $ the lemma  follows.
\end{proof}%
\begin{remark}[Lemma~\ref{lm:randomalphamNbounded} using tail probability]
We see, with $0 < \alpha < 1$
\begin{displaymath}
\EE{X | E} \le \overline{f_{\alpha}}\sum_{k=0}^{\uN}\prob{X > k | E}  = \overline{f_{\alpha}} \sum_{k=1}^{\uN} \left( \frac{1}{k^{\alpha}} - \frac{1}{(\uN + 1)^{\alpha} } \right) \sim   \frac{\alpha\overline{f_{\alpha}}}{1-\alpha}      \left( \uN^{1-\alpha} - 1 \right) 
\end{displaymath}
in agreement with the choice of $\varepsilon_{N}$ in
Lemma~\ref{lm:randomalphamNbounded}.
\end{remark}
\begin{lemma}[$S_{N}/(Nm_{N}) \to 1$ almost surely]
\label{lm:almostsureconv}
Suppose $ X_{1}, \ldots, X_{N}$ are independent positive
integer-valued random variables with
$X_{1}, \ldots, X_{N} \vartriangleright  \mathds L(\alpha, \uN)$ (recall
\eqref{eq:28}) with probability $ \varepsilon_{N}$, and
$X_{1},\ldots, X_{N} \vartriangleright  \mathds L(\at, \uN)$ with probability
$1 - \varepsilon_{N}$, where $1 \le \alpha < 2 $ and $\at \ge 2 $
fixed.  Recall \breyta{$S_{N} = \sum_{i=1}^{N}X_{i}$} from
\eqref{SN} and  \breyta{$\EE{X_{1}} =  m_{N}$ from
\eqref{eq:5} in }
Definition~\ref{def:not}, and choose $\varepsilon_{N}$   as in
Lemma~\ref{lm:randomalphamNbounded} such that
$\limsup_{N\to \infty}m_{N} < \infty$. Then $S_{N}/(Nm_{N}) \to 1$
almost surely as $N\to \infty$.
\end{lemma}
\begin{proof}[Proof of Lemma~\ref{lm:almostsureconv}]
Let $\overline X_{i} = (X_{i} - m_{N})/(Nm_{N})$ for $i \in [N]$.  The
$\overline X_{1}, \ldots, \overline X_{N}$ are i.i.d.\ and
$\EE{\overline X_{1}} = 0$.  Then
$\sum_{i=1}^{N} \overline X_{i} \to 0$ almost surely
\citep{Etemadi1981}.  Hence,
$S_{N}/(Nm_{N}) = 1 + \sum_{i=1}^{N}\overline{X}_{i} \to 1$ almost
surely as $N\to \infty$.
\end{proof}
We verify \eqref{eq:cNthmrandall} under the conditions of
Case~\ref{item:2} of Theorem ~\ref{thm:haplrandomalpha}. \breyta{  Recall
$c_{N}$ from  Definition~\ref{def:cNhapl} and \eqref{eq:6}.}
\begin{lemma}[$c_{N}$ under  Theorem ~\ref{thm:haplrandomalpha}]
\label{lm:cNhaplrandomall}
Suppose the conditions of Theorem~\ref{thm:haplrandomalpha} hold with
$1 \le \alpha < 2$, $\at \ge 2$, and $c > 0$ all fixed, and that
$\uN/N \gneqq 0$.  Let $L \equiv L(N)$ be a function of $N$ such that
$L/N \to 0$, and
\begin{equation}
\label{eq:33}
\begin{split}
\varepsilon_{N}  =  cN^{\alpha - 2} \left( \one{\kappa > 2} + \one{\kappa = 2}\log N \right)   
\end{split}
\end{equation}
Then, with $\kappa + 2 < c_{\kappa} < \kappa^{2}$ when $\kappa > 2$
and with $C_{\kappa}^{N} $ as in \eqref{eq:CNmap}, $\gamma$ as in
\eqref{eq:27}, $B(\gamma,2-\alpha,\alpha)$ as in
Definition~\ref{def:betacoal}, $f_{\alpha}^{(\infty)} = g_{\alpha}^{(\infty)}$ and
$\underline {g_{\kappa}}(2) = \overline {f_{\kappa}}(2)$, we have
$ \lim_{N\to \infty} C_{\kappa}^{N} c_{N} = C_{\kappa,\alpha,\gamma}$ where
$C_{\kappa,\alpha,\gamma}$ as in \eqref{eq:23} and $m_{\infty}$ as in
\eqref{eq:49}, i.e.\
\begin{equation}
\label{eq:boundscNrandall}
\begin{split}
&  C_{\kappa,\alpha,\gamma}  =   \frac{  2 \overline {f_{\kappa}}(2) }{m_{\infty}^{2}}\left( \one{\kappa = 2}  + \one{\kappa > 2}  \frac{ c_{\kappa} }{2^{\kappa}(\kappa -2)(\kappa - 1) } \right)  +       \frac{ c\alpha f_{\alpha}^{(\infty)} }{m_{\infty}^{\alpha}}  B(\gamma, 2-\alpha,\alpha) 
\end{split}
\end{equation}
\end{lemma}
Choosing  $\varepsilon_{N}$ \breyta{(recall
Definition~\ref{def:haplrandomalpha} of type $A$ random environment)}  as  in Lemma~\ref{lm:cNhaplrandomall}
gives  $\varepsilon_{N}/c_{N} \gg 0$ in all cases and so we 
expect with high probability to see at least one event $E$ 
\breyta{(the event $X_{1},\ldots,X_{N} \vartriangleright
\mathds L(\alpha,\zeta(N))$)}  over a period of time proportional to  $1/c_{N}$
generations. 
\begin{proof}[Proof of Lemma~\ref{lm:cNhaplrandomall}]
 Choosing $\varepsilon_{N}$
as in \eqref{eq:33} yields a finite $m_{\infty}$ by
Lemma~\ref{lm:randomalphamNbounded} and almost sure convergence of
$S_{N}/(Nm_{N})$ to 1   by Lemma~\ref{lm:almostsureconv}.
Let, for $0 < \delta < 1$ fixed and  such that $(1-\delta)m_{N} > 1$,     
\begin{equation}
\label{eq:41}
\begin{split}
M_{+} &  :=  (1 + \delta)Nm_{N}, \quad M_{-}   :=  (1 -\delta)Nm_{N}, \\
Y_{+} & :=  \frac{(X_{1})_{2}}{(X_{1} + M_{+})^{2}}, \quad  Y_{-}  :=  \frac{(X_{1})_{2}}{(X_{1} + M_{-})^{2}}, \\
R_{N} &  :=  \frac{(X_{1} )_{2}}{S_{N}^{2}}\one{S_{N}\ge N}
\end{split}
\end{equation}
where \breyta{$\EE{X_{1}} = m_{N}$ and $m_{\infty} =
\lim_{N\to \infty}m_{N}$} as in \eqref{eq:5} in
Definition~\ref{def:not}.  By Theorem~\ref{thm:Schwthm4} (recall
$c_{N}$ from  Definition~\ref{def:cNhapl})
\begin{displaymath}
c_{N} \sim N\EE{R_{N}}
\end{displaymath}
Adapting the arguments of \citep[Lemma~13]{schweinsberg03} and using
Lemma~\ref{lm:almostsureconv} we have for any  fixed  $\epsilon > 0$
(and assuming $N$ is large enough that $N-1\approx N$)
\begin{equation}
\label{eq:46}
(1-\epsilon)\EE{Y_{+} } \le  \EE{R_{N}} \le  \epsilon \EE{\frac{(X_{1})_{2} }{\max\{ X_{1}^{2},N^{2} \} } } + \EE{Y_{-}} 
\end{equation}
With \breyta{ $C_{\kappa}^{N} =   \one{\kappa > 2}N +
\one{\kappa = 2}N/\log N $}  as in \eqref{eq:CNmap}  we now  see that  it   suffices to consider lower  bounds on    $C_{\kappa}^{N}N\EE{Y_{+}}$ and upper
bounds on $C_{\kappa}^{N} N \EE{Y_{-}}$ 
 and to check that
\begin{equation}
\label{eq:EXovermaxfin}
\limsup_{N\to \infty} C_{\kappa}^{N}  N \EE{ \frac{(X_{1})_{2}}{ \max\{X_{1}^{2},N^{2} \} } } < \infty
\end{equation}
We  check that  \eqref{eq:EXovermaxfin} holds.  Since $(X_{1} + N)^{2} \le 4\max\{X_{1}^{2}, N^{2}\}$ we see 
\begin{displaymath}
N\EE{ \frac{(X_{1})_{2}}{ \max\{X_{1}^{2}, N^{2} \} } } \le  4N \EE{\frac{ (X_{1})_{2}}{ (X_{1} + N)^{2} } }
\end{displaymath}
and so to   \eqref{eq:EXovermaxfin}  combining  Lemma~\ref{lm:boundssumrka} and \eqref{eq:33}.

Write $r_{k}(a,M) \equiv  (k)_{2}(k+M)^{-2}(k^{-a} - (k + 1)^{-a})$.
Equation~\eqref{eq:boundscNrandall} will follow from
Lemma~\ref{lm:boundssumrka} after substituting for $\ell$ and $m$.
When $\alpha = 1$
\begin{displaymath}
\begin{split}
   \EE{Y_{-} | E } &  = \sum_{k=1}^{\uN}\frac{ k(k-1) }{(k + M_{-})^{2} }\prob{X_{1} =k | E}  \le   \overline{f_{\alpha}}(1)\sum_{k=1}^{L-1}r_{k}(\alpha,M_{-} ) +  \overline{f_{\alpha}}(L)\sum_{k=L}^{\uN}r_{k}(\alpha,M_{-} ) \\
   & \le \overline{f_{\alpha}}(L)  \frac{\uN + 1}{(M_{-} + L)(M_{-} + \uN +1)}       + O(LM_{-}^{-2})
   \end{split} 
\end{displaymath}
after substituting for  $\ell$ and $m$ in Case~\ref{item:10}  of  Lemma~\ref{lm:boundssumrka} so that  by \eqref{eq:33}  when $\alpha = 1$ 
\begin{displaymath}
\limsup_{N\to \infty} C_{\kappa}^{N} N  \EE{Y_{-} | E}\varepsilon_{N} \le
\begin{cases}
cf_{\alpha}^{(\infty)} \frac{K }{(1-\delta)m_{\infty}(K + (1-\delta)m_{\infty}) } &  \text{when $\uN/N \to K$}  \\
cf_{\alpha}^{(\infty)} ((1-\delta)m_{\infty})^{-1} &  \text{when $\uN/N \to \infty$}   
\end{cases}
\end{displaymath}
In the same way we get
\begin{displaymath}
   \EE{Y_{+} | E }  = \sum_{k=1}^{\uN}\frac{ k(k-1) }{(k + M_{+})^{2} }\prob{X_{1} =k | E}  \ge   \underline{g_{\alpha}}(1)\sum_{k=1}^{L-1}r_{k}(\alpha,M_{+} ) +  \underline{g_{\alpha}}(L)\sum_{k=L}^{\uN}r_{k}(\alpha,M_{+} ) \\
\end{displaymath}
 and when $\alpha = 1$ we have   the lower bound 
\begin{displaymath}
\liminf_{N\to \infty} C_{\kappa}^{N} N  \EE{Y_{+} | E}\varepsilon_{N} \ge
\begin{cases}
cg_{\alpha}^{(\infty)} \frac{K }{(1+\delta)m_{\infty}(K + (1+\delta)m_{\infty}) }  &  \text{when $\uN/N \to K$}  \\
cg_{\alpha}^{(\infty)} ((1+\delta)m_{\infty})^{-1}  &  \text{when $\uN/N \to \infty$}   
\end{cases}
\end{displaymath}
Write
\begin{displaymath}
\begin{split}
\gamma_{-} &  :=  \one{ \frac{\uN}{N} \to K  } \frac{K }{(1-\delta)m_{\infty} + K} + \one{\frac{\uN }{N} \to \infty } \\
\gamma_{+} &  :=  \one{ \frac{\uN}{N} \to K  } \frac{K }{(1 + \delta)m_{\infty} + K} + \one{\frac{\uN }{N} \to \infty }
\end{split}
\end{displaymath}
When $1 < \alpha < 2$ combining   Case~\ref{item:12} of  Lemma~\ref{lm:boundssumrka} with $\ell = 1$ and $m = \uN$  and  \eqref{eq:33} we get   the bounds 
\begin{displaymath}
\begin{split}
& \limsup_{N\to \infty}C_{\kappa}^{N} N \EE{Y_{-}| E}\varepsilon_{N} \le \frac{ c\alpha f_{\alpha}^{(\infty)}}{ ((1-\delta)m_{\infty})^{\alpha} }\int_{0}^{1}\one{0 < u \le \gamma_{-}} u^{1-\alpha}(1-u)^{\alpha-1}du  \\
& \liminf_{N\to \infty}C_{\kappa}^{N} N  \EE{Y_{+}| E}\varepsilon_{N} \ge \frac{ c\alpha g_{\alpha}^{(\infty)}}{ ((1+\delta)m_{\infty})^{\alpha} }\int_{0}^{1}\one{0 < u \le \gamma_{+}} u^{1-\alpha}(1-u)^{\alpha-1}du 
\end{split}
\end{displaymath}
When $\kappa = 2$ combining  Case~\ref{item:13}  of  Lemma~\ref{lm:boundssumrka} with $\ell = 2$ and $m = \uN$  and \eqref{eq:33}  we get the bounds 
\begin{displaymath}
\begin{split}
& \limsup_{N\to \infty} C_{\kappa}^{N} N \EE{Y_{-} | E^{\sf c} }(1- \varepsilon_{N}) \le   2 c\overline {f_{\kappa}}(2)  ((1-\delta)m_{\infty})^{-2}  \\
& \liminf_{N\to \infty} C_{\kappa}^{N} N \EE{Y_{+} | E^{\sf c} }(1- \varepsilon_{N}) \ge   2 c\underline {g_{\kappa}}(2)  ((1+\delta)m_{\infty})^{-2}  
\end{split}
\end{displaymath}
When $\kappa > 2$ combining Case~\ref{item:14}    of  Lemma~\ref{lm:boundssumrka}   with $\ell = 2$ and $m = \uN$  and \eqref{eq:33}  we get
\begin{displaymath}
\begin{split}
\limsup C_{\kappa}^{N} N \EE{Y_{-} | E^{\sf c}}(1- \varepsilon_{N}) &  \le  \overline {f_{\kappa}}(2) \left( \frac{\kappa 2^{2-\kappa} }{\kappa - 2} - \frac{\kappa 2^{1-\kappa} }{\kappa - 1 }  \right) = \frac{ 2^{1-\kappa}\kappa ^{2}}{ (\kappa -2)(\kappa - 1)  } \overline {f_{\kappa}}(2)          \\
\liminf C_{\kappa}^{N} N \EE{Y_{+} | E^{\sf c}}(1- \varepsilon_{N}) &  \ge  \underline {g_{\kappa}}(2) \left( \frac{\kappa 2^{2-\kappa} }{\kappa - 2} - \frac{3\kappa 2^{1-\kappa} }{\kappa - 1 } +                 2^{1-\kappa} \right)  =   \frac{ 2^{1-\kappa}(\kappa + 2) }{(\kappa-2)(\kappa - 1) } \underline {g_{\kappa}}(2)   \\
\end{split}
\end{displaymath}
noting that, when $\ell = 2$, the leading terms  in the lower resp.\ upper bound  in    \eqref{eq:kappatwoplus} simplify  to
\begin{displaymath}
     \frac{1}{(1 + M)^{2} }\frac{2^{1-a} (a+2) }{(a-2)(a-1) } \quad \text{resp.\  } \quad  \frac{1}{(2 + M)^{2}}\frac{ 2^{1-a} a^{2}}{(a-2)(a-1) }
\end{displaymath}
The lemma   now follows from the assumption $f_{\alpha}^{(\infty)} =
g_{\alpha}^{(\infty)}$ and $\underline {g_{\kappa}}(2) = \overline
{f_{\kappa}}(2)$ and taking $\epsilon$ and $\delta$ to 0. The proof of
Lemma~\ref{lm:cNhaplrandomall} is complete.
 \end{proof}


\begin{lemma}[Identifying the $\Lambda_{+}$-measure]
\label{lm:tailprobability}
Under the conditions of part~\ref{item:2} of
Theorem~\ref{thm:haplrandomalpha} with $\gamma$ as in \eqref{eq:27}, 
$c$ from \eqref{eq:33}, \breyta{ $S_{N} =
\sum_{i=1}^{N}X_{i}$} as in \eqref{SN},  \breyta{ $\lim_{N\to
\infty} 
C_{\kappa}^{N} c_{N}$}  as in
\eqref{eq:boundscNrandall},  $f_{\alpha}^{(\infty)} =
g_{\alpha}^{(\infty)}$,   $C_{\kappa,\alpha,\gamma}$ as in
\eqref{eq:23},  \breyta{$m_{\infty} = \lim_{N\to
\infty}\EE{X_{1}}$ as in \eqref{eq:5},}
\begin{displaymath}
\begin{split}
     \lim_{N\to \infty}   \frac{N}{c_{N}} \prob{\frac{X_{1}}{S_{N}}\one{S_{N} \ge N}  \ge x} = \frac{\alpha c  f_{\alpha}^{(\infty)}}{C_{\kappa,\alpha,\gamma} m_{\infty}^{\alpha} }  \int_{x}^{1 }\one{0<x\le \gamma} u^{-1-\alpha} (1-u)^{\alpha-1} du
\end{split}
\end{displaymath}
\end{lemma}
\begin{proof}[Proof of Lemma~\ref{lm:tailprobability}]
The proof follows the proof of \citep[Lemma~14]{schweinsberg03}.   Let $\epsilon > 0$ and $\delta >
0$ both  be fixed with  $(1-\delta)m_{N} > 1$.     Then, by 
Lemma~\ref{lm:almostsureconv}  we can  take   $N$ large enough that  
\begin{displaymath}
\prob{  M_{-} \le  \widetilde{S}_{N} \le  M_{+}  } > 1 - \epsilon 
\end{displaymath}
holds  with  \breyta{$M_{-} =  (1-\delta)Nm_{N}$ and $M_{+} =
(1+\delta)Nm_{N}$}  from \eqref{eq:41}  and
\breyta{ $\widetilde{S}_{N} = \sum_{i=2}^{N}X_{i} $}  from
\eqref{SN2}  in  Definition~\ref{def:not}.   For any
\breyta{fixed}   $0 < x < 1$
\begin{multline}
\label{eq:boundstailprob}
 (1-\epsilon)\prob{\frac{(X_{1})_{2}}{(X_{1} + M_{+})^{2} } \ge x } \le   \prob{\frac{X_{1} }{S_{N}}\one{S_{N} \ge N} \ge x  } \le  \epsilon \prob{ X_{1} \ge Nx} \\  +  \prob{\frac{X_{1}}{X_{1} + M_{-} } \ge x }  
\end{multline}
 Recall  event $E$ \breyta{( $X_{1},\ldots,X_{N}
 \vartriangleright \mathds L(\alpha,\zeta(N))$)} from
Definition~\ref{def:haplrandomalpha} \breyta{(type $A$ random
environment)}  and that by
Lemma~\ref{lm:randomalphamNbounded} and Lemma~\ref{lm:cNhaplrandomall}
we can choose $(\varepsilon_{N})$ as in \eqref{eq:33} so that
$m_{\infty} < \infty$.  Using \eqref{eq:12} we see
\breyta{(recall $m_{N} = \EE{X_{1}}$ as in  \eqref{eq:5})}
\begin{displaymath}
\begin{split}
& \prob{\frac{X_{1}}{X_{1} + M_{-}} \ge x | E}   =  \prob{X_{1} \ge \frac{x}{1-x}M_{-} | E } \\
& \le  \left( \left( \frac{1-x}{x}\right)^{\alpha}\frac{1}{((1-\delta)Nm_{N})^{\alpha}}  -  \frac{1}{(\uN + 1)^{\alpha}}  \right)\overline{f_{\alpha}} \left( \frac{x}{1-x}M_{-}  \right)   \\
 &  =  \frac{ \overline{f_{\alpha}}  \left( \frac{x}{1-x}M_{-}  \right)   }{ ( (1-\delta)Nm_{N})^{\alpha}  }\left(  \left( \frac{1-x}{x}  \right)^{\alpha}   -  \left( \frac{  (1-\delta)Nm_{N}   }{ \uN +1}  \right) ^{\alpha}  \right) \\
\end{split}
\end{displaymath}
For  $0 < \alpha < 2$,   $w \ge 0$,  and  $0 < x \le  1/(1 + w)$ we
have \citep[Equation~10.2]{Eldon2026biorxiv}
\begin{equation}
\label{eq:idLambdaintegral}
\left( \frac{1-x}{x}  \right)^{\alpha} -  w^{\alpha} =  \alpha \int_{x}^{\frac{1}{1 + w}}u^{-1-\alpha}(1-u)^{\alpha - 1}du  
\end{equation}
using the substitution $y = (1-u)/u$.  Write  $K_{N}^{-} \equiv    M_{-} /(\uN +1) $ so that 
\begin{equation}
\label{eq:widecheckgamma}
\lim_{N\to \infty}\frac{1}{1 + K_{N}^{-}}  =   \one{ \frac{\uN}{N} \to K } \frac{ K}{K + (1-\delta)m_{\infty} } + \one{  \frac{\uN}{N} \to \infty  } \equiv  \overline\gamma
\end{equation}
By Lemma~\ref{lm:cNhaplrandomall},  recalling
$\varepsilon_{N}$ from   \eqref{eq:33},  with $1 \le \alpha < 2$, 
\begin{displaymath}
\begin{split}
& \limsup_{N\to \infty} \frac{N}{c_{N}} \prob{\frac{X_{1}}{X_{1} + M_{-}} \ge x | E}\varepsilon_{N} \\
& \le \limsup_{N\to \infty} \frac{ \alpha  \overline{f_{\alpha}}   \left( \frac{x}{1-x}M_{-}  \right)    }{ ((1-\delta)m_{N})^{\alpha}  } \frac{N^{1-\alpha}  \varepsilon_{N} }{c_{N}}   \int_{x}^{\frac{1}{1 + K_{N}^{-} } }  u^{-1-\alpha}(1-u)^{\alpha - 1}du  \\
& =      \frac{\alpha c f_{\alpha}^{(\infty)} }{C_{\kappa,\alpha,\gamma} ((1-\delta)m_{\infty})^{\alpha} } \int_{x}^{1}\one{0 < x \le \overline \gamma } u^{-1-\alpha}(1-u)^{\alpha - 1}du  
\end{split}
\end{displaymath}
using that $\lim_{k\to \infty}f_{\alpha}(k) = f_{\alpha}^{(\infty)}$ (recall
\eqref{eq:isfg}), and 
$C_{\kappa,\alpha,\gamma}$ is as in 
\eqref{eq:boundscNrandall} in Lemma~\ref{lm:cNhaplrandomall}.

We consider the  lower bound  in \eqref{eq:boundstailprob}.    
Write   $K_{N}^{+} \equiv   M_{+}   /( \uN +1) $ so that
\begin{equation}
\label{eq:gammahat}
\lim_{N\to \infty} \frac{1}{1 + K_{N}^{+}}  =  \one{ \frac{\uN}{N} \to K } \frac{ K}{K + (1+\delta)m_{\infty} } + \one{  \frac{\uN}{N} \to \infty  }  \equiv  \widehat \gamma
\end{equation}
We see  using \eqref{eq:12} and that $\lim_{k\to\infty}g_{\alpha}(k)=g_{\alpha}^{(\infty)}$ (recall \eqref{eq:isfg}) and with $1 \le \alpha < 2$   
\begin{displaymath}
\begin{split}
& \liminf_{N\to \infty}(1-\epsilon)\frac{N}{c_{N}} \prob{\frac{X_{1}}{X_{1}  + M_{+}}  \ge x | E } \varepsilon_{N} \\
\ge &  \liminf_{N\to \infty}\frac{(1-\epsilon)\underline {g_{\alpha}}   \left( \frac{x}{1-x}M_{+}  \right)  }{((1+\delta)m_{N})^{\alpha}}\frac{N^{1-\alpha}\varepsilon_{N} }{c_{N}}\int_{x}^{\frac{1}{1 + K_{N}^{+}} }u^{-1-\alpha}(1-u)^{\alpha - 1}du  \\
=  &  \frac{(1-\epsilon) \alpha c {g}_{\alpha}^{(\infty)}}{C_{\kappa,\alpha,\gamma} ((1+\delta)m_{\infty})^{\alpha} } \int_{x}^{1} \one{0<x\le \widehat \gamma  } u^{-1-\alpha}(1-u)^{\alpha - 1}du  
\end{split}
\end{displaymath}
The lemma follows after taking $\epsilon$ and $\delta$ to zero and
recalling that $f_{\alpha}^{(\infty)} = g_{\alpha}^{(\infty)}$ by
assumption. The proof of Lemma~\ref{lm:tailprobability} is complete.
\end{proof}

The following lemma shows that our choice \eqref{eq:33} of
$\varepsilon_{N}$ guarantees that simultaneous mergers vanish in the
limit when $1 \le \alpha < 2$. 
\begin{lemma}[Simultaneous mergers vanish]
\label{lm:simmergersvanisharandall}
\breyta{Condition}~\eqref{eq:simtozero} holds  under the conditions of Case~\ref{item:2} of
Theorem~\ref{thm:haplrandomalpha}. 
\end{lemma}
\begin{proof}[Proof of Lemma~\ref{lm:simmergersvanisharandall}]
We see, following the proof of \citep[Lemma~15]{schweinsberg03},
\breyta{recalling  $S_{N}=\sum_{i=1}^{N}X_{i}$ as in
\eqref{SN}, $c_{N}$ from Definition~\ref{def:cNhapl},
$\varepsilon_{N}$ from Definition~\ref{def:haplrandomalpha} (type $A$
random environment)} 
\begin{displaymath}
\begin{split}
&  \frac{N^{2}}{c_{N}}\EE{\frac{(X_{1})_{2}(X_{2})_{2} }{S_{N}^{4}}\one{S_{N} \ge N} } 
\le    \frac{N^{2}}{c_{N}}\EE{\frac{(X_{1})_{2}(X_{2})_{2} }{\max\{X_{1}^{2}, N^{2}\}  \max\{X_{2}^{2}, N^{2}\}   }\one{S_{N} \ge N} } \\
& \le  \frac{N^{2}}{c_{N}}\left( \EE{\frac{(X_{1})_{2}}{ \max\{X_{1}^{2}, N^{2}\} }  } \right)^{2} \le   16  \frac{N^{2}}{c_{N}}\left( \EE{\frac{(X_{1})_{2}}{ (X_{1} + N)^{2} }  } \right)^{2}  \\
\end{split}
\end{displaymath}
When $\alpha = 1$ we use Case~\ref{item:10} of
Lemma~\ref{lm:boundssumrka} to see that
$\EE{(X_{1})_{2}(X_{1}+N)^{-2} |E }\varepsilon_{N} \overset{c}{\sim}
N^{-2}$ and Case~\ref{item:13} of the same lemma to see that
$\EE{ (X_{1})_{2}(X_{1} + N)^{2} | E^{\sf c} } \overset c \sim
N^{-2}\log N $ when $\kappa = 2$  as $N\to \infty$ and  with $\varepsilon_{N}$ as in
\eqref{eq:33}; hence
$\limsup_{N\to \infty} \svigi{N^{2}/c_{N}} \svigi{\EE{ (X_{1})_{2}(X_{1}+N)^{-2}}
}^{2} = 0$ using \eqref{eq:cNthmrandall}; the case for
$1<\alpha < 2$ follows similarly, hence the lemma using
\eqref{eq:SchwLm6}.
\end{proof}

\begin{lemma}[Convergence to the  Kingman coalescent]
\label{lm:convKingmanrall}
Suppose the  conditions of Case~\ref{item:1} of Theorem~\ref{thm:haplrandomalpha} hold.  Then $Nc_{N}^{-1}\EE{(X_{1})_{3}S_{N}^{-3}\one{S_{N}\ge N} } \to 0$.   
\end{lemma}
\begin{proof}[Proof of Lemma~\ref{lm:convKingmanrall}]
 On $S_{N}\ge N$ we have  $S_{N}^{3} \ge \max\{X_{1}^{3}, N^{3} \}$ so that with $\uN < N$ we have 
 \begin{displaymath}
  \frac{N}{c_{N}}\EE{\frac{(X_{1})_{3} }{S_{N}^{3}} \one{S_{N} \ge N} } \le \frac{N}{c_{N}}\EE{\frac{(X_{1})_{3} }{\max\{X_{1}^{3}, N^{3} \} } }
\end{displaymath}
When  $\alpha > 1$  \eqref{eq:21} gives
\begin{displaymath}
\begin{split}
 \frac{N}{c_{N}}\EE{\frac{(X_{1})_{3} }{\max\{X_{1}^{3}, N^{3} \}} | E }\varepsilon_{N} &  \le \frac{\alpha  \overline {f_{\alpha}} N \varepsilon_{N} }{c_{N} N^{3} }\sum_{k=2}^{\uN} \frac{((k -1) + 1)^{3}}{(k-1)^{1+\alpha}} \le   \frac{4 \alpha  \overline {f_{\alpha}}  \varepsilon_{N} }{c_{N} N^{2} }\sum_{k=2}^{\uN} \frac{(k -1)^{3} + 1}{(k-1)^{1+\alpha}} \\
\end{split}
\end{displaymath}
Combining   \eqref{eq:33} and   Lemma~\ref{lm:cNhaplrandomall} we see with  $1 < \alpha < 2$ 
\begin{displaymath}
 \frac{4 \alpha  \overline {f_{\alpha}}  \varepsilon_{N} }{c_{N} N^{2} }\sum_{k=2}^{\uN} \frac{(k -1)^{3} + 1}{(k-1)^{1+\alpha}} \le  { 4\alpha c \overline {f_{\alpha}}  } \frac{  \uN^{3-\alpha} }{N^{3-\alpha}} +  O(N^{\alpha - 3})
\end{displaymath}
When  $ \alpha  =   1$ we have, by the upper bound in \eqref{eq:20} and, since $\uN/N \to 0$ by assumption,  
 \begin{displaymath}
 \frac{N}{c_{N}}\EE{\frac{(X_{1})_{3} }{\max\{X_{1}^{3}, N^{3} \} } | E }\varepsilon_{N}  \le  \frac{\overline {f_{\alpha}} N \varepsilon_{N} }{c_{N} N^{3} } \sum_{k=1}^{\uN }k \le  \frac{\overline {f_{\alpha}} \varepsilon_{N} }{c_{N}N^{2}}\uN(\uN + 1) + O(N^{-1})
\end{displaymath}
The lemma follows from Lemma~\ref{lm:cNhaplrandomall} with $\varepsilon_{N}$ as in \eqref{eq:33} so that $\limsup_{N\to \infty} \varepsilon_{N}/c_{N} < \infty$.
\end{proof}
\begin{remark}[Convergence to the Kingman coalescent when $0 < \alpha < 1$]
\label{rm:convKsmallalpha}
Suppose the conditions of Theorem~\ref{thm:haplrandomalpha} hold with
$\uN/N \to 0$ and $\kappa > 2$ and  $0 < \alpha < 1$.  Then (assuming $Nc_{N} \sim c$
and $N\varepsilon_{N} \sim c^{\prime}$ for $c,c^{\prime} > 0$ fixed)  similar calculations as in the
proof of Lemma~\ref{lm:convKingmanrall} give that
\begin{displaymath}
\frac{N}{c_{N}}\EE{\frac{(X_{1})_{3} }{\max\{X_{1}^{3}, N^{3} \} }|E }\varepsilon_{N} \le  c^{\prime\prime} N^{-2}\uN^{3-\alpha}  + O(N^{2-\kappa})
\end{displaymath}
for some fixed $c^{\prime\prime} > 0$.  Hence, by
\citep[Proposition~4]{schweinsberg03} convergence to the
Kingman coalescent only holds if $\uN^{3-\alpha}/N^{2} \to 0$ when
$0 < \alpha < 1$ under our assumptions on $c_{N}$ and
$\varepsilon_{N}$. \breyta{Taking $\zeta(N) = N^{x}$ shows
that convergence to the Kingman coalescent follows if  $x < 2/
(3-\alpha) < 1$, compare with Case~\ref{item:1} of Theorem~\ref{thm:haplrandomalpha}. }
 \end{remark}

\begin{proof}[Proof of parts \ref{item:1} and \ref{item:2} of Theorem~\ref{thm:haplrandomalpha}]
Part~\ref{item:1} of Theorem~\ref{thm:haplrandomalpha} is
Lemma~\ref{lm:convKingmanrall} when $1 \le \alpha < 2$ using
\eqref{eq:SchwLm6} and \citep[Proposition~2]{schweinsberg03} (recall
\eqref{condi}).

We turn to part \ref{item:2} of Theorem~\ref{thm:haplrandomalpha}.  We
need to check the three conditions of \citep[Proposition~3]{schweinsberg03};
these are \eqref{eq:cntozero}, \eqref{eq:simtozero}, and
\eqref{eq:existenceLambda} in
Proposition~\ref{prp:conditionsconvergenceLambda}.
Lemma~\ref{lm:cNhaplrandomall} gives \eqref{eq:cntozero}; choosing
$\varepsilon_{N}$ as in \eqref{eq:33} gives
Lemma~\ref{lm:randomalphamNbounded} so that
Lemma~\ref{lm:almostsureconv} holds and with that
Lemma~\ref{lm:cNhaplrandomall}.  Condition~\eqref{eq:simtozero} is Lemma~\ref{lm:simmergersvanisharandall}.  
 We turn to  
\eqref{eq:existenceLambda}.  We need to show that (recall $\gamma$
from \eqref{eq:27}, \breyta{and $\nu_{1}$ is the random
number of surviving offspring of (arbitrary) individual 1})
\begin{equation}
\label{eq:convtail}
\lim_{N\to \infty} \frac{N}{c_{N}} \prob{\nu_{1} > Nx } =  \int_{x}^{1} y^{-2}\Lambda_{+}dy =  C\int_{x}^{1} \one{0 < x \le \gamma}y^{-1-\alpha}y^{\alpha - 1}dy
\end{equation}
for $0 < x < 1$. By the arguments of the proof of 
part (c) of \breyta{Theorem}~4 of \citep{schweinsberg03} using bounds on the tail
probability of a hypergeometric \citep{Chvtal1979} we have
\begin{subequations}
\begin{align}
\label{eq:supinftails}
\limsup_{N\to\infty} \frac{N}{c_{N}} \prob{\nu_{1} > Nx} & \le \limsup_{N\to \infty} \frac{N}{c_{N}} \prob{ \frac{X_{1}}{S_{N}} \ge x - \epsilon } \\
\label{eq:22}
\liminf_{N\to \infty} \frac{N}{c_{N}} \prob{\nu_{1} > Nx} & \ge \liminf_{N\to \infty} \frac{N}{c_{N}} \prob{ \frac{X_{1}}{S_{N}} \ge x + \epsilon } 
\end{align}
\end{subequations}
for $0 < \epsilon \le x$.  Lemma~\ref{lm:tailprobability} together
with \eqref{eq:supinftails} and \eqref{eq:22} and taking $\epsilon$
to 0 gives \eqref{eq:convtail}.  
\end{proof}


\begin{proof}[Proof of part~\ref{item:4} of Theorem~\ref{thm:haplrandomalpha}]
\breyta{Finally, we check  part~\ref{item:4} of Theorem~\ref{thm:haplrandomalpha}, 
 the case  $0 < \alpha < 1$ and $\uN/N^{1/\alpha} \to \infty$.}

Following  \citep[\S~4]{schweinsberg03} let   $Y_{(1)} \ge Y_{(2)} \ge \ldots \ge Y_{(N)}$ be  the ranked values of  $N^{-1/\alpha}X_{1}, \ldots, N^{-1/\alpha} X_{N}$.
With $p_{i} =  \prob{X_{1} \ge N^{1/\alpha}x_{i} | E} - \prob{X_{1}
\ge N^{1/\alpha} x_{i-1} | E}$ using  \eqref{eq:12}  we see, 
\breyta{recalling   $E$ is the event
$X_{1},\ldots,X_{N}\vartriangleright \mathds L(\alpha,\zeta(N))$ from
Definition~\ref{def:haplrandomalpha} (type $A$ random environment),
and  $(x_{j})_{j\in \N}$ is a decreasing  sequence of positive numbers
with $x_{0} = \infty$,  }
\begin{displaymath}
\begin{split}
&   \frac{\underline {g_{\alpha}}(N^{1/\alpha}x_{i}  ) }{N}  \left( x_{i}^{-\alpha} -  (\tfrac{N^{1/\alpha} }{\uN })^{\alpha} \right) -     \frac{\overline {f_{\alpha}}(N^{1/\alpha}x_{i-1} ) }{N}  \left( x_{i-1}^{-\alpha} -  (\tfrac{N^{1/\alpha} }{\uN +1 })^{\alpha} \right) \\
& = \frac{1}{N}\svigi{ \underline {g_{\alpha}}(N^{1/\alpha}x_{i} )  x_{i}^{-\alpha} - \overline {f_{\alpha}} (N^{1/\alpha}x_{i-1} ) x_{i-1}^{-\alpha} } + O(\uN^{-\alpha})   \\
&   \le p_{i} \le  \frac{\overline {f_{\alpha}} (N^{1/\alpha}x_{i} )  }   {N}  \left( x_{i}^{-\alpha} -  (\tfrac{N^{1/\alpha} }{\uN + 1 })^{\alpha} \right) -     \frac{\underline {g_{\alpha}} (N^{1/\alpha}x_{i-1}) }{N}  \left( x_{i-1}^{-\alpha} -  (\tfrac{N^{1/\alpha} }{\uN  })^{\alpha} \right)  \\
& =    \frac{1}{N}\left( \overline {f_{\alpha}} (N^{1/\alpha}x_{i} ) x_{i}^{-\alpha} - \underline {g_{\alpha}} (N^{1/\alpha}x_{i-1}) x_{i-1}^{-\alpha} \right) +  O(\uN^{-\alpha}) 
\end{split}
\end{displaymath}
for $1 \le i \le j$ so that $p = 1 - p_{1} - \ldots - p_{j} = 1 - \prob{Y_{1} \ge x_{j}| E}$.  Moreover,  for non-negative integers $n_{1}, \ldots, n_{j}$ we have
\begin{displaymath}
p^{N - n_{1} - \cdots - n_{j}} = \left( 1 -  \prob{X_{1} \ge N^{1/\alpha}x_{j} | E}  \right)^{N - n_{1} - \cdots - n_{j}}
\end{displaymath}
Using \eqref{eq:12}  again we have
\begin{multline}
\left(1 -  \frac{\overline{f_{\alpha}}(N^{1/\alpha}x_{j} ) }{Nx_{j}^{\alpha}} \left(1 - O\left( (\tfrac{N^{1/\alpha} }{\uN + 1})^{\alpha} \right)    \right) \right)^{N - \sum_{i}n_{i} } \le p^{N - \sum_{i}n_{i} } \\  \le \left(1 -  \frac{\underline{g_{\alpha}}(N^{1/\alpha}x_{j}) }{Nx_{j}^{\alpha}}  \left(1 - O\left( (\tfrac{N^{1/\alpha} }{\uN })^{\alpha} \right)    \right) \right)^{N - \sum_{i}n_{i} }  
\end{multline}
so that 
\begin{displaymath}
\begin{split}
  \exp\svigi{- f_{\alpha}^{(\infty)} x_{j}^{-\alpha} } \le \lim_{N\to \infty}  p^{N - \sum_{i}n_{i} } \le   \exp \svigi{-  g_{\alpha}^{(\infty)} x_{j}^{-\alpha} } 
\end{split}
\end{displaymath}
Then
\begin{displaymath}
\begin{split}
&   \prod_{i=1}^{j} \frac{ \exp \svigi{-f_{\alpha}^{(\infty)} ( x_{i}^{-\alpha} -  x_{i-1}^{-\alpha}) } \svigi{ {g}_{\alpha}^{(\infty)}  x_{i}^{-\alpha} - {f}_{\alpha}^{(\infty)} x_{i-1}^{-\alpha} }^{n_{i}} }{n_{i}!}  \\
\le &     \lim_{N\to\infty} \frac{ (N)_{n_{1} + \cdots + n_{j}}}{n_{1}!  \cdots n_{j}!} p_{1}^{n_{1}}\cdots  p_{j}^{n_{j}}p^{N - \sum_{i}n_{i}} \\
\le &  \prod_{i=1}^{j} \frac{ \exp\svigi{- g_{\alpha}^{(\infty)}  ( x_{i}^{-\alpha} -  x_{i-1}^{-\alpha}) }\svigi{ {f}_{\alpha}^{(\infty)}  x_{i}^{-\alpha} - {g}_{\alpha}^{(\infty)}x_{i-1}^{-\alpha} }^{n_{i}}        }{n_{i}!}
\end{split}
\end{displaymath}
Choosing $g_{\alpha}$ and $f_{\alpha}$ in \eqref{eq:PXiJ} such  that
$g_{\alpha}^{(\infty)} =  f_{\alpha}^{(\infty)}$  we have
$\prob{L_{1}^{N} = n_{1}, \ldots, L_{j}^{N } = n_{j}} \to \prob{L_{1}
= n_{1}, \ldots, L_{j} = n_{j} }$ where
$(L_{1}^{N}, \ldots, L_{j}^{N})$ and $(L_{1}, \ldots, L_{j})$ are as
in the proof of   \citep[Lemma 20]{schweinsberg03}.


As briefly mentioned in \S~\ref{back} Poisson-Dirichlet distributions
can be constructed from the ranked points
$Z_{1} \ge Z_{2} \ge \ldots $ of a Poisson point process with
characteristic measure
$\nu_{\alpha}((x,\infty)) = \one{x > 0}Cx^{-\alpha}$.   The
law of $(Z_{j}/\sum_{i}Z_{i})_{j\in \IN}$ is the Poisson-Dirichlet distribution with
parameters $(\alpha ,0)$.

The upper bound in  \eqref{eq:20} gives (recalling $\uN/N^{1/\alpha} \to \infty$)
\begin{displaymath}
\begin{split}
\EE{\sum_{i=1}^{N}Y_{i}\one{Y_{i} \le \delta} | E} & =  N^{1-1/\alpha}\EE{X_{1}\one{X_{1} \le N^{1/\alpha}\delta} | E }   \le  \overline {f_{\alpha}} N^{1-1/\alpha}  \sum_{k=1}^{\lfloor N^{1/\alpha}\delta \rfloor }k^{-\alpha}    \\
& \le \overline {f_{\alpha}}  N^{1-1/\alpha}  \int_{0}^{N^{1/\alpha}\delta}x^{-\alpha}dx  =  \overline {f_{\alpha}} \frac{\delta^{1-\alpha}}{1-\alpha} 
\end{split}
\end{displaymath}
On $E$ we then have that
$\left( Y_{1}, \ldots, Y_{j}, \sum_{i=j+1}^{N} Y_{i} \right)$
converges weakly to
$(Z_{1}, \ldots, Z_{j}, \sum_{i=j+1}^{\infty}Z_{i})$ as $N\to \infty$
by the arguments for \citep[Lemma~21]{schweinsberg03}.  With
$(W_{j})_{j\in \IN} = (Z_{j}/\sum_{i}Z_{i})_{j\in \IN}$ it follows
that, provided  $ C_{\kappa}^{N}  c_{N} \overset{c}{\sim} 1$ and with
$\varepsilon_{N} = c/C_{\kappa}^{N}$ for some $c > 0$,
\begin{displaymath}
\begin{split}
&  \lim_{N\to \infty} \frac{N^{r}}{c_{N}} \EE{ \frac{(X_{1})_{k_{1}} \cdots (X_{r})_{k_{r}} }{S_{N}^{k_{1} + \cdots + k_{r}} } \one{S_{N} \ge N }  } \\
& =   \lim_{N\to \infty} \frac{N^{r}}{c_{N}} \EE{ \frac{(X_{1})_{k_{1}} \cdots (X_{r})_{k_{r}} }{S_{N}^{k_{1} + \cdots + k_{r}} } \one{S_{N} \ge N } |E }\varepsilon_{N}  \\
& +       \lim_{N\to \infty} \frac{N^{r}}{c_{N}} \EE{ \frac{(X_{1})_{k_{1}} \cdots (X_{r})_{k_{r}} }{S_{N}^{k_{1} + \cdots + k_{r}} } \one{S_{N} \ge N } |E^{\sf c} }(1 - \varepsilon_{N}) \\
& = \frac{c}{c(1-\alpha)  + C_{\kappa} } \sum_{\substack{  i_{1}, \ldots, i_{r} = 1 \\ \text{all distinct}}}^{N} \EE{ W_{i_{1}}^{k_{1}}\cdots W_{i_{r}}^{k_{r}}  }  +  \one{r=1, k_{1}= 2}\frac{C_{\kappa} }{c(1-\alpha) + C_{\kappa} } \\
& =  \frac{c }{c(1-\alpha) + C_{\kappa} } \int_{\Delta_{+}}  \sum_{\substack{  i_{1}, \ldots, i_{r} = 1 \\ \text{all distinct}}}^{\infty} x_{i_{1}}^{k_{1}}\cdots x_{i_{r}}^{k_{r}}\frac{1}{\sum_{j=1}^{\infty}x_{j}^{2} }\Xi_{\alpha} (dx)         +    \one{r=1, k_{1}= 2}\frac{C_{\kappa} }{c(1-\alpha) + C_{\kappa} }
\end{split}
\end{displaymath}
for all $k_{1}, \ldots, k_{r} \ge 2$ and $r\in \IN$ where
$\Xi_{\alpha}$ is as in Definition~\ref{def:poisson-dirichlet} with
$C_{\kappa}$ as in \eqref{eq:37} and $C_{\kappa}^{N}$ as in
\eqref{eq:CNmap}.

It remains to verify that $C_{\kappa}^{N} c_{N} \overset{c}{ \sim}
1$. We have
$N\EE{ (X_{1})_{2}S_{N}^{-2}\one{S_{N} \ge N} | E } \to 1 - \alpha$
\citep[pp.\ 137]{schweinsberg03}.  By Lemma~\ref{lm:cNhaplrandomall}
$C_{\kappa}^{N}\overset{c}{\sim} N\EE{ (X_{1})_{2}S_{N}^{-2}\one{S_{N}
\ge N} | E^{\sf c} }$.  Then
$C_{\kappa}^{N} c_{N} \to c(1-\alpha) + C_{\kappa}$.
Part~\ref{item:4} of Theorem~\ref{thm:haplrandomalpha}  now follows
from our choice of $\varepsilon_{N}$ and \citep[Proposition~1; Theorem
2.1]{schweinsberg03,MS01} where
\begin{displaymath}
\begin{split}
\lim_{N\to \infty} \frac{ \EE{ (\nu_{1})_{k_{1}} \cdots (\nu_{r})_{k_{r}} } }{N^{k_{1} + \cdots + k_{r} - r }c_{N} } &  =  \frac{C_{\kappa}}{c(1-\alpha)  + C_{\kappa}}\one{r=1, k_{1} = 2} \\
&  +   \frac{c }{c(1-\alpha) + C_{\kappa} }   \int_{\Delta_{+}}  \sum_{\substack{  i_{1}, \ldots, i_{r} = 1\\ \text{all distinct}}}^{\infty} x_{i_{1}}^{k_{1}}\cdots x_{i_{r}}^{k_{r}}\frac{1}{\sum_{j=1}^{\infty}x_{j}^{2} }\Xi_{\alpha} (dx) 
\end{split}
\end{displaymath}
The proof of  Theorem~\ref{thm:haplrandomalpha} is complete.
\end{proof}

\begin{remark}[The case $0< \alpha < 1$ and $\uN/N^{1/\alpha} \to K$]
\label{rm:incompletePD}
Suppose, under the conditions of Theorem~\ref{thm:haplrandomalpha}, we
have $0 < \alpha < 1$ and $\uN/N^{1/\alpha} \to K$
\breyta{($K > 0$ fixed).}

\breyta{Let     $K = x_{0} \ge x_{1} \ge x_{2} \ge \ldots
\ge x_{j} > 0$ all be fixed.}
   Then, with $f_{\alpha}^{(\infty)} = g_{\alpha}^{(\infty)}$
and    recalling \breyta{event}   $E$ from  Definition~\ref{def:haplrandomalpha}
\breyta{(type $A$ random environment)}, 
\begin{displaymath}
  \prob{X_{1} \ge N^{1/\alpha } x_{i} | E}  \overset{c}{\sim}   N^{-1}( x_{i}^{-\alpha} - K^{-\alpha})
\end{displaymath}
Taking   $p_{i} =  \prob{X_{1} \ge N^{1/\alpha} x_{i} | E } - \prob{X_{1} \ge N^{1/\alpha }x_{i-1} | E}$  for all  $1 \le i \le j $  we have 
$p_{i} \overset{c}{\sim}  N^{-1}( x_{i}^{-\alpha} -  x_{i-1}^{-\alpha})$.  Moreover,  taking  $p = 1 - p_{1} - \cdots - p_{j}$ we see
\begin{multline}
p^{N - \sum_{i}n_{i} } = \left(1 - \prob{X_{1} \ge N^{1/\alpha}x_{j} | E } \right)^{N - \sum_{i}n_{i}} \overset{c}{\sim} \left( 1 - N^{-1}(x_{j}^{-\alpha} - K^{-\alpha}) \right)^{N - \sum_{i}n_{i}} \\  \overset{c} \sim e^{-(x_{j}^{-\alpha} - K^{-\alpha}) }
\end{multline}
for $n_{1}, \ldots, n_{j}\in \IN\cup \{0\}$.   Then 
\begin{displaymath}
\begin{split}
\frac{(N)_{n_{1} + \cdots + n_{j}} }{n_{1}! \cdots n_{j}!} p_{1}^{n_{1}} \cdots p_{j}^{n_{j}}p^{N - \sum_{i}n_{i}} & \overset c  \sim \frac{N^{\sum_{i} n_{i}}}{n_{1}! \cdots n_{j}! }\prod_{i=1}^{j} N^{-n_{i}}(x_{i}^{-\alpha} - x_{i-1}^{-\alpha} )^{n_{i}} e^{-(x_{j}^{-\alpha} - K^{-\alpha} ) } \\
& =  \prod_{i=1}^{j} \frac{1}{n_{i}!} e^{-(x_{i}^{-\alpha} - x_{i-1}^{-\alpha}) }(x_{i}^{-\alpha} - x_{i-1}^{-\alpha})^{n_{i}}
\end{split}
\end{displaymath}
We conjecture that the last expression corresponds to the probability
$\prob{\bigcap_{i=1}^{j}\{ L_{i} = n_{i} \}}$ where $L_{i}$ is the
number of points of a Poisson point process
$\{\mathscr N_{\alpha,K} \}$ between $x_{i}$ and $x_{i-1}$ where
$\EE{L_{i}} = x_{i}^{-\alpha} - x_{i-1}^{-\alpha}$ with $x_{0} = K$
and $\{\mathscr N_{\alpha,K} \}$ has intensity measure
$\mathbb{\nu}_{\alpha,K}((x,\infty)) = \one{0 < x \le K} (x^{-\alpha}
- K^{-\alpha}) $.  Suppose $Z_{1} \ge Z_{2} \ge \ldots $ are the
ranked points of $\{\mathscr N_{\alpha,K} \}$.  Then, in contrast to
the case when $\uN/N^{1/\alpha} \to \infty$, it is not clear what the
law of $(Z_{j}/\sum_{i}Z_{i})_{j\in\IN }$ converges to.
\end{remark}

\subsection{Proof of Theorem~\ref{thm:hapl-alpha-random-one}}
\label{sec:proof-thm-random-alpha-one}
In this section we give a proof of
Theorem~\ref{thm:hapl-alpha-random-one} \breyta{(recall the
type $B$ random environment from 
Definition~\ref{def:alpha-random-one})}.
   \breyta{To briefly orient the reader through the proof of
Theorem~\ref{thm:hapl-alpha-random-one}  (the proofs can always be skipped on first reading to get an overview
of the lemmas leading to  the theorem),   we  use the  same approach
as for  Theorem~\ref{thm:haplrandomalpha}, namely checking the
conditions of  Proposition~\ref{prp:conditionsconvergenceLambda} for
Case~\ref{item:9} of Theorem~\ref{thm:hapl-alpha-random-one},  and
Condition~\eqref{condi}  for Case~\ref{item:7} of the theorem.        Recall from
Definition~\ref{def:alpha-random-one}  that  $\overline
\varepsilon_{N}$ is the probability of event $E_{1}$ (that one randomly
picked $X_{i} \vartriangleright \mathds L(\alpha,\zeta(N))$, and all
other $X_{j}\vartriangleright \mathds L(\kappa,\zeta(N))$),  
Lemma~\ref{lm:finitemean-alpha-random-one} identifies conditions on
$\overline\varepsilon_{N}$ such that   $m_{\infty} =
\lim_{N\to \infty}\EE{X} < \infty$, and in Lemma~\ref{lm:law-large-numbers-alpha-random-one}  further conditions
on  $\overline \varepsilon _{N}$ such that  $S_{N}/(Nm_{N})$ converges
almost surely to 1, in turn useful for  deriving the limits of
$C_{\kappa}^{N}c_{N}$ in Lemma~\ref{lm:cNrandomalphaone}.
Lemma~\ref{lm:Lambdameasurealphaone} identifies the
$\Lambda_{+}$-measure (recall Condition \eqref{eq:existenceLambda} in
Proposition~\ref{prp:conditionsconvergenceLambda}),
Lemma~\ref{lm:simmergersvanishrandone} checks that simultaneous
mergers vanish in the limit (recall Condition~\ref{eq:simtozero} of
Proposition~\ref{prp:conditionsconvergenceLambda}),  and
Lemma~\ref{lm:convkingmanralphaone} verifies  Case~\ref{item:7} of
Theorem~\ref{thm:hapl-alpha-random-one} by checking
Condition~\eqref{condi}.       }

First we identify conditions
for \breyta{ $m_{\infty} = \lim_{N\to \infty}\EE{X}$}  from
\eqref{eq:5} to be finite, \breyta{with $X$ denoting the random number of potential offspring of an arbitrary
individual.}
\begin{lemma}[Finite $m_\infty$]
\label{lm:finitemean-alpha-random-one}
Under the conditions of Theorem~\ref{thm:hapl-alpha-random-one}, with $\uN/N \gneqq 0$ and
$\log \uN \overset{c}{\sim} \log N$.  Recalling
$(\overline \varepsilon_{N})_{N}$ from
Definition~\ref{def:alpha-random-one} \breyta{($\overline \varepsilon_{N}$
is the probability of  event $E_{1}$ in type $B$ random environment)}
suppose
\begin{equation}
\label{eq:epsilonfinitemNalpharandone}
\overline \varepsilon_{N} \in 
\begin{cases}
O(1)  & \text{when $\alpha = 1$} \\
O( N\uN^{ \alpha - 1 })     & \text{when $0 < \alpha < 1$} \\
\end{cases}
\end{equation}
Then we have   $\limsup_{N\to \infty}\EE{X} < \infty$.
\end{lemma}
\begin{proof}[Proof of Lemma~\ref{lm:finitemean-alpha-random-one}]
We see, recalling the events $E_{1}$ and $E_{1}^{\sf c}$ from
Definition~\ref{def:alpha-random-one} \breyta{(type $B$
random environment;  event $E_{1}$
is when $X_{i}\vartriangleright \mathds L(\alpha,\zeta(N))$, and
$X_{j}\vartriangleright \mathds L(\kappa,\zeta(N))$ for all $j\in
[N]\setminus\{i\}$ with index $i$  picked uniformly at random; event
$E_{1}^{\sf c}$ is when  $X_{1},\ldots, X_{N} \vartriangleright
\mathds L(\kappa,\zeta(N))$)}, 
\begin{displaymath}
\EE{X} =  \EE{X | E_{1}}\overline \varepsilon_{N} +   \EE{X | E_{1}^{\sf c}}(1- \overline \varepsilon_{N}  )
\end{displaymath}
When  $\alpha = 1$   and $\at \ge   2$ using the upper bound in \eqref{eq:21} in Lemma~\ref{lm:ineqs}  
\begin{displaymath}
\begin{split}
\EE{X | E_{1}} &  \le  \frac{\overline{f_{\alpha}}}{N}\sum_{k=1}^{\uN} \frac{1}{k+1} +    \frac{N-1}{N}  \overline{f_{\kappa}}\sum_{k=1}^{\uN} k\left( \frac{1}{k^{\at } } -   \frac{1}{(1+k)^{\at } }  \right) \\
& \le       \frac{\overline f_{\alpha}  \log(\uN + 1)}{N} +  \kappa  \overline{f_{\kappa}} \sum_{k=1}^{\uN} k^{-\kappa} 
\end{split}
\end{displaymath}
When $\alpha  \in (0,1)$,  we have $\EE{X | E_{1}}
\le  c \uN^{1-\alpha}/N$ by \eqref{eq:43} and so the lemma by \eqref{eq:epsilonfinitemNalpharandone}.
\end{proof}

We  verify a convergence  result \breyta{ on $S_{N}/(Nm_{N})$  similar to
Lemma~\ref{lm:almostsureconv};  the result will be useful for deriving 
 limits of  $C_{\kappa}^{N}c_{N}$. }
\begin{lemma}[$S_{N}/(Nm_{N}) \to 1$ almost surely]
\label{lm:law-large-numbers-alpha-random-one}
Suppose the conditions of Theorem~\ref{thm:hapl-alpha-random-one} hold
with $m_{N} = \EE{X_{1}}$ and  $\alpha \in  (0,1]$. Let   $\overline \varepsilon _{N}$ from
Definition~\ref{def:alpha-random-one} \breyta{(type $B$
random environment)}  fulfill the condition from 
Lemma~\ref{lm:finitemean-alpha-random-one} and, with   $1 < r < 2$ fixed, 
\begin{equation}
\label{eq:50}
\overline \varepsilon_{N} \in  
\begin{cases}
O(1) &  \text{when $\uN/N \to K$} \\
o\left(1\wedge  \frac{N^{r} }{\uN^{r-\alpha}}\right) &  \text{when $\uN/N \to \infty$} \\
\end{cases}
\end{equation}
Then $S_{N}/(Nm_{N})$  converges almost surely  to 1  as $N\to \infty$
with  \breyta{$ S_{N} = \sum_{i=1}^{N}X_{i}$} as in  \eqref{SN}.
\end{lemma}
\begin{proof}[Proof of Lemma~\ref{lm:law-large-numbers-alpha-random-one}]
We follow the approach of \citep{PANOV2017379}.  Our choice  of
$\overline \varepsilon_{N}$  guarantees that
$\limsup_{N\to \infty}m_{N} < \infty$ by
Lemma~\ref{lm:finitemean-alpha-random-one}.
Write $\overline X_{i} =  (X_{i} - m_{N})/(Nm_{N})$, and 
\begin{displaymath}
\overline S_{N}  =  \sum_{i=1}^{N} \overline X_{i}  = \frac{S_{N} - Nm_{N} }{Nm_{N} } = \frac{S_{N}}{Nm_{N}} - 1
\end{displaymath}
Fix  $1 < r < 2$.  
 Provided      $\limsup_{N\to \infty} \EE{|\overline X_{1} |^{r} }  < \infty$ 
we have   by \citep[Theorem~2]{10.1214/aoms/1177700291} 
\begin{displaymath}
\EE{|\overline S_{N}|^{r}}  \le   \frac{2}{N^{r}m_{N}^{r} } \sum_{i=1}^{N}  \EE{ |X_{i} - m_{N}|^{r} } =   \frac{2}{m_{N}^{r}} N^{1-r} \EE{ |X_{1} - m_{N}|^{r} }
\end{displaymath}
 Then, using Lemma~\ref{lm:ineqs} with 
  $\alpha  \in (0,1]$,  and with $m_{N} < \infty$ for all $N$  by Lemma~\ref{lm:finitemean-alpha-random-one} 
\begin{displaymath}
\begin{split}
\EE{|X_{1} - m_{N}|^{r} | E_{1}} & \le   |m_{N}|^{r}\prob{X_{1} = 0} +     \frac{ \overline{f_{\alpha}}}{N} \sum_{k= 1}^{\uN}|k - m_{N}|^{r} \frac{1}{k^{1 + \alpha}}    +  \kappa  \overline{f_{\kappa}} \sum_{k= 1}^{\uN} |k-m_{N}|^{r}\frac{1}{k^{1+\kappa}} \\
& \le  |m_{N}|^{r} +    \frac{   \overline{f_{\alpha}} 2^{r-1} }{N} \sum_{k= 1}^{\uN}\frac{k^{r} + m_{N}^{r} }{k^{1+\alpha}}  + 2^{r-1} \kappa \overline {f_{\kappa}}  \sum_{k= 1}^{\uN} \frac{k^{r} + m_{N}^{r}}{k^{1+\kappa}}  \\ 
\end{split}
\end{displaymath}
where we have used that $|x+ y|^{r} \le 2^{r-1}(|x|^{r} + |y|^{r})$
for any $x,y$ real and $r \ge 1$ \citep[Equation~(1.12),
Proposition~3.1.10{\it (iii)} pp.\ 87]{athreya06:_measur}
\breyta{(follows easily from Jensen's inequality since
$\mathds R \ni  x\mapsto |x|^{r}$ is convex for $r \ge 1$)}.  Recalling
$0 < \alpha \le 1 < r < 2 \le \kappa $ we see
\begin{displaymath}
\begin{split}
 &  \frac{1 }{N} \sum_{k= 1}^{\uN} k^{r - \alpha -1}  \le \frac{1 }{N}  +    \one{r \le  \alpha + 1} \frac{1 }{N}   \int_{1}^{\uN} x^{r-\alpha-1 } dx  + \one{r > \alpha + 1}\frac{1 }{N} \int_{1}^{\uN + 1} x^{r-\alpha - 1} dx  \\
 &  =  O(  N^{-1}\uN^{r-\alpha})  
\end{split}
\end{displaymath}
Then $\EE{|X_{1} - m_{N}|^{r} | E_{1} } \in O(1 +  N^{-1}\uN^{r-\alpha})$.  
With $\overline \varepsilon_{N}$ as in \eqref{eq:50} we then have
\begin{displaymath}
\limsup_{N\to \infty} N^{1-r}\EE{|X_{1} - m_{N}|^{r} | E_{1} }\overline \varepsilon_{N}  = 0
\end{displaymath}
for any $1 < r < 2$.  By similar
calculations
$\limsup_{N\to \infty} N^{1-r}\EE{|X_{1} - m_{N}|^{r} | E_{1}^{\sf c} } (1
- \overline \varepsilon_{N}) = 0$; there is a $r > 1$ such that
$\limsup_{N\to \infty} \EE{ |\overline{S}_{N}|^{r} } = 0$.  Hence
$\overline S_{N}$ converges in probability to $0$, and since the
$\overline X_{i}$ are independent, $\overline S_{N}$ converges almost
surely to 0 (see e.g.\  \citep[Theorem~8.3.3 pp~251]{athreya06:_measur}). Hence, $S_{N}/(Nm_{N}) = 1 + \sum_{i=1}^{N} \overline{X}_{i}  \to 1$ almost surely.  
\end{proof}
\begin{remark}[Convergence of $S_{N}/(Nm_{N})$]
Using the upper bounds in \eqref{eq:20} and \eqref{eq:21} in
Lemma~\ref{lm:ineqs} it is straightforward to check that
$\EE{\overline X_{1}^{2}} = O \svigi{ N^{-3}\uN^{2-\alpha}\overline
\varepsilon_{N} + N^{-2}\log \uN } $. Then
$N\EE{\overline X_{1}^{2}} = O\svigi{ N^{-2}\uN^{2-\alpha}\overline
\varepsilon_{N} + N^{-1}\log \uN} $.  \breyta{ Taking}
\begin{displaymath}
\overline \varepsilon_{N} = O \svigi{ 1 \wedge N^{2}\uN^{\alpha - 2} \wedge
N\uN^{\alpha - 1}}
\end{displaymath}
then gives $\limsup_{N\to \infty} m_{N} < \infty$
by Lemma~\ref{lm:finitemean-alpha-random-one} and
$\limsup_{N\to \infty}N\EE{\overline X_{1}^{2}} < \infty $.  By
Lemma~\ref{lm:law-large-numbers-alpha-random-one}
$\overline S_{N} \to 0$ in probability, $\EE{\overline X_{1}} = 0$,  so  by
Khinchine-Kolmogorov's 1-series theorem $\overline S_{N} \to 0$ almost
surely. Then $S_{N}/(Nm_{N}) = 1 +  \overline S_{N}  \to 1$ almost surely. 
\end{remark}

We verify the scaling of $c_{N}$ (recall Definition~\ref{def:cNhapl}). 
\begin{lemma}[Scaling of $c_N$]
\label{lm:cNrandomalphaone}
Suppose the conditions of Theorem~\ref{thm:hapl-alpha-random-one} hold
with $0 < \alpha \le 1$, $\at \ge 2$, $c > 0$ and $0 < \varepsilon  < 1$
all fixed.  \breyta{Recall $C_{\kappa}^{N}
= \one{\kappa > 2}N + \one{\kappa = 2}N/\log N$ as in \eqref{eq:CNmap}.}    Let $(\overline \varepsilon_{N})_{N}$ from
Definition~\ref{def:alpha-random-one} \breyta{(type $B$
random environment;  $\overline \varepsilon_{N} $ is the probability
of event $E_{1}$ that a randomly picked $X_{i} \vartriangleright
\mathds L(\alpha,\zeta(N)$ and all other $X_{j}\vartriangleright
\mathds L(\kappa,\zeta(N))$)}  take the values
\begin{equation}
\label{eq:51}
\overline \varepsilon_{N} = \begin{cases}   cN^{\alpha - 1} & 
\text{when $0 < \alpha < 1$ and    $\at > 2$} \\
 c N^{\alpha - 1} \log N &
\text{when  $0 < \alpha < 1$ and   $\at = 2$} \\
\varepsilon   & \text{when $\alpha =  1$ and $\at > 2$} \\
\end{cases}
\end{equation}
and suppose $f_{\alpha }^{(\infty)} = g_{\alpha}^{(\infty)}$.  Then $C_{\kappa}^{N}c_{N} \to C_{\kappa,\alpha,\gamma}$ with $C_{\kappa}^{N}$
as in \eqref{eq:CNmap} and   $C_{\kappa,\alpha,\gamma}$ as in  
\eqref{eq:boundscNrandall} except    $c_{\alpha} = \one{0 < \alpha < 1}c + \one{\alpha = 1} \varepsilon$ replaces  $c$. 
\end{lemma}

\begin{remark}[Condition on $\overline \varepsilon_{N}$ in Lemma~\ref{lm:cNrandomalphaone}]
When $\uN/N\to K$ Lemma~\ref{lm:law-large-numbers-alpha-random-one}
imposes no new restrictions on $\overline \varepsilon_{N}$, as well as
when $\uN \overset{c}{\sim} N^{\gamma}$ with
$1 < \gamma < r/(r-\alpha)$; since $r/(r-\alpha) \le r + 1 - \alpha$
so that $N^{\alpha - 1} =  O(N^{r}\uN^{\alpha - r})$, and
$\limsup_{N\to \infty} m_{N} < \infty$ by
Lemma~\ref{lm:finitemean-alpha-random-one}. Further,
$N^{\alpha - 1 } = O(N^{1 - \gamma (1-\alpha)})$ for any
$\gamma > 1$. Thus, a mild condition on $\uN$ when $\uN/N \to \infty$
allows to choose $\overline \varepsilon_{N}$ as in \eqref{eq:51} such
that $S_{N}/(Nm_{N}) \to 1$ almost surely by
Lemma~\ref{lm:law-large-numbers-alpha-random-one} and  $\limsup_{N\to \infty} m_{N} < \infty$ by Lemma~
\ref{lm:finitemean-alpha-random-one}.  
\end{remark}

\begin{proof}[Proof of Lemma~\ref{lm:cNrandomalphaone}]
Recall $Y_{-}$ and $Y_{+}$ from \eqref{eq:41}. Similar to the proof of
Lemma~\ref{lm:cNhaplrandomall} we use
Lemma~\ref{lm:law-large-numbers-alpha-random-one} to get
\eqref{eq:46}. It then suffices to consider lower bounds on
$C_{\kappa}^{N} N\EE{Y_{+}}$    and upper bounds on
$C_{\kappa}^{N} N\EE{Y_{-}}$ and to check that \eqref{eq:EXovermaxfin}
holds.  These all follow as in the proof of
Lemma~\ref{lm:cNhaplrandomall}, noting that by
Definition~\ref{def:alpha-random-one}, writing $p_{k}(a) \equiv k^{-a} - (1+k)^{-a}$ for $k\in \N$ and $a > 0$, 
\begin{displaymath}
\begin{split}
  \EE{Y_{-} } &  \le  \left( \frac{\overline{f_{\alpha}}(1)}{N} \sum_{k=1}^{L-1}\frac{k(k-1)}{(k + M_{-})^{2}}p_{k}(\alpha) +  \frac{\overline{f_{\alpha}}(L)}{N}  \sum_{k=L}^{\uN}\frac{k(k-1)}{(k + M_{-})^{2}}p_{k}(\alpha)    \right) \overline \varepsilon_{N}   \\
   & +   \left( \overline{f_{\kappa}}(1) \frac{N- 1 }{N} \sum_{k=1}^{L-1}\frac{k(k-1)}{(k + M_{-})^{2}} p_{k}(\kappa) +  \overline{f_{\kappa}}(L) \frac{N-1 }{N}  \sum_{k=L}^{\uN}\frac{k(k-1)}{(k + M_{-})^{2}}p_{k}(\kappa) \right) \overline \varepsilon_{N}     \\
& +  \left( {\overline{f_{\kappa}}(1)} \sum_{k=1}^{L-1}\frac{k(k-1)}{(k + M_{-})^{2}}  p_{k}(\kappa)  +  {\overline{f_{\kappa}}(L)}  \sum_{k=L}^{\uN}\frac{k(k-1)}{(k + M_{-})^{2}} p_{k}(\kappa)   \right)(1 -  \overline \varepsilon_{N} )
\end{split}
\end{displaymath}
A lower bound on $\EE{Y_{+}}$ is of a similar form with
$\underline {g_{\alpha}}, \underline{g_{\kappa}}$ replacing
$\overline {f_{\alpha}}, \overline{f_{\kappa}}$, and the lemma follows
as in the proof of Lemma~\ref{lm:cNhaplrandomall} upon substituting
for $\ell $ and $m$ in Lemma~\ref{lm:boundssumrka} and taking limits.
\end{proof}


\begin{lemma}[Identifying the $\Lambda_{+}$-measure]
\label{lm:Lambdameasurealphaone}
Under the conditions of Theorem~\ref{thm:hapl-alpha-random-one} for all
$0 < x < 1$ with \breyta{ $m_{\infty} = \lim_{N\to
\infty}\EE{X_{1}} $}  as in \eqref{eq:5} and  \breyta{$S_{N}
= \sum_{i=1}^{N}X_{N}$ as in \eqref{SN}}  we have 
\begin{equation}
\label{eq:Lambdameasurebounds}
\begin{split}
   \lim_{N\to \infty}  \frac{N}{c_{N}}\prob{ \frac{X_{1}}{S_{N}}\one{S_{N} \ge N} \ge x  } =  \frac{\alpha c_{\alpha} f_{\alpha}^{(\infty)} }{ C_{\kappa,\alpha,\gamma}   m_{\infty}^{\alpha}}  \int_{x}^{1} \one{0 < x \le \gamma} y^{-1-\alpha}(1-y)^{\alpha -1 }dy
\end{split}
\end{equation}
\breyta{ with $\gamma =   \one{\zeta(N)/N\to \infty} +     \one{\zeta(N)/N \to K}K/(m_{\infty} + K)$ as in
\eqref{eq:27},} 
 $ C_{\kappa,\alpha, \gamma}$ as in  \eqref{eq:23} with
$c_{\alpha} = \one{0 < \alpha < 1}c +   \one{\alpha = 1}\varepsilon $
replacing $c$.
\end{lemma}
\begin{proof}[Proof of Lemma~\ref{lm:Lambdameasurealphaone}]
The proof follows the one of Lemma~\ref{lm:tailprobability}. When
$\alpha = 1$ and $\kappa > 2$ with $\overline\varepsilon_{N}$ as in
\eqref{eq:51} so that $ m_{\infty} < \infty$ by
Lemma~\ref{lm:finitemean-alpha-random-one} (recall \breyta{$M_{+} = (1+\delta)Nm_{N} $
from \eqref{eq:41}})
\begin{displaymath}
\begin{split}
& \prob{\frac{X_{1}}{X_{1} + M_{+} } \ge x | E_{1} }\overline \varepsilon_{N}  =  \prob{X_{1} \ge  \frac{x}{1-x} M_{+} | E_{1} }\overline \varepsilon_{N} \\
& \ge \frac{\underline {g_{\alpha}} ( \frac{x}{1-x}M_{+}  )   \varepsilon }{N}\left( \frac{1}{M_{+}} \frac{1-x}{x}  - \frac{1 }{\uN + 1} \right) \\
& +  \frac{\underline {g_{\kappa}} ( \frac{x}{1-x}M_{+}  ) \varepsilon }{M_{+}^{\kappa} } \frac{N-1}{N}  \left(  \left( \frac{1-x}{x} \right)^{\kappa} - \left(\frac{M_{+}}{\uN + 1} \right)^{\kappa} \right) 
\end{split}
\end{displaymath}
Using  Lemma~\ref{lm:cNrandomalphaone} (recall $\gamma$ from \eqref{eq:27}) and since (with $\alpha = 1$)
\begin{equation}
\label{eq:alphaoneint}
\frac{1-x}{x} - \one{\frac{\uN }{N}\to K }\frac{(1+\delta)m_{\infty} }{K} =   \int_{x}^{1}\one{0 < x \le \widehat\gamma} y^{-2}dy =   \int_{x}^{1}\one{0 < x \le \widehat\gamma} y^{-\alpha - 1}(1-y)^{\alpha-1}dy
\end{equation}
with $\widehat \gamma$ as in \eqref{eq:gammahat}  we get (after checking that $Nc_{N}^{-1}\prob{X_{1} \ge (X_{1} + M_{+})x | E_{1}^{\sf c}} \to 0 $)
\begin{displaymath}
\begin{split}
\liminf_{N\to \infty} \frac{N}{c_{N}} \prob{\frac{X_{1}}{X_{1} + M_{+} } \ge x | E_{1} }\overline \varepsilon_{N} & \ge   \frac{\varepsilon g_{\alpha}^{(\infty)}}{m_{\infty}}  \int_{x}^{1}\one{0 < x \le \widehat\gamma} y^{-2}dy
\end{split}
\end{displaymath}
after taking $\delta$ to 0.  For the upper bound we see (recall
\breyta{$M_{-} = (1-\delta)Nm_{N}$} from \eqref{eq:41}, and $\alpha = 1$)
\begin{displaymath}
\begin{split}
& \prob{\frac{X_{1} }{X_{1} + M_{-}} \ge x | E_{1} }\overline \varepsilon_{N} = \prob{X_{1}  \ge  \frac{x}{1-x}M_{-} | E_{1}}\overline \varepsilon_{N} \\
& \le  \frac{\overline{f_{\alpha}}( \frac{x }{1-x }M_{-} ) \varepsilon }{N}\left( \frac{1}{M_{-}} \frac{1-x}{x} -  \frac{1}{\uN + 1}  \right)  +  \frac{\overline{f_{\kappa}}( \frac{x }{1-x }M_{-} )  }{M_{-}^{\kappa} } \left( \left( \frac{1-x}{x}\right)^{\kappa} - \left(\frac{M_{-} }{\uN} \right)^{\kappa} \right) \\
\end{split}
\end{displaymath}
Analogously to the lower bound we get (when $\alpha = 1$ and $\kappa > 2$) with $\overline \gamma$ from \eqref{eq:widecheckgamma} 
\begin{displaymath}
\begin{split}
\limsup_{N\to \infty} \frac{N}{c_{N}} \prob{ \frac{X_{1}}{X_{1} + M_{-}} \ge x  } \le  \frac{\varepsilon f_{\alpha}^{(\infty)}}{m_{\infty}}  \int_{x}^{1}\one{0 < x \le \overline \gamma} y^{-2}dy
\end{split}
\end{displaymath}
The result  when $\alpha = 1$  follows after taking $\delta$ to 0. 
The case $0 < \alpha < 1$ proceeds in the same way where we make use of
\eqref{eq:idLambdaintegral}. 
\end{proof}
\begin{lemma}[Simultaneous mergers vanish]
\label{lm:simmergersvanishrandone}
Condition~\eqref{eq:simtozero} follows from the settings of
Theorem~\ref{thm:hapl-alpha-random-one}.
\end{lemma}
\begin{proof}[Proof of Lemma~\ref{lm:simmergersvanishrandone}]
The proof follows the arguments for
Lemma~\ref{lm:simmergersvanisharandall}.  Choosing
$\overline \varepsilon_{N}$ as in \eqref{eq:51} then by
Case~\ref{item:15} of Lemma~\ref{lm:boundssumrka} we have
$\EE{(X_{1})_{2}(X_{1} + N)^{-2} | E_{1} } \overline\varepsilon_{N}
\overset{c}{\sim} N^{-2}$ as $N\to \infty$ so that by
Lemma~\ref{lm:cNrandomalphaone}
$\limsup_{N\to \infty}N^{2}(\EE{ (X_{1})_{2}(X_{1} + N) })^{2}/c_{N} =
0$. Hence to \eqref{eq:simtozero} by the arguments for
Lemma~\ref{lm:simmergersvanisharandall}.
\end{proof}
\begin{remark}[Condition~\eqref{eq:simtozero} follows from  $m_{\infty} < \infty$]
Since $m_{\infty} < \infty$ under Theorem~\ref{thm:hapl-alpha-random-one}
(by Lemma~\ref{lm:finitemean-alpha-random-one} choosing
$\overline\varepsilon_{N}$ as in \eqref{eq:51}) and Case~\ref{item:2}
of Theorem~\ref{thm:haplrandomalpha} (by
Lemma~\ref{lm:randomalphamNbounded} choosing $\varepsilon_{N}$ as in
\eqref{eq:33}) condition \eqref{eq:simtozero} follows from Lemma~15 of \citep{schweinsberg03}.  
\end{remark}

\begin{lemma}[Convergence to the Kingman coalescent]
\label{lm:convkingmanralphaone}
When $\uN/N \to 0$ Case~\ref{item:7} of Theorem~\ref{thm:hapl-alpha-random-one} follows.
\end{lemma}
\begin{proof}[Proof of Lemma~\ref{lm:convkingmanralphaone}]
By \citep[Proposition~4]{schweinsberg03} (see also
\citep[\S~4]{mohle2000total}) it suffices to check that \eqref{condi}
holds.  By \eqref{eq:SchwLm6} in Theorem~\ref{thm:Schwthm4} it is
sufficient to show that
$Nc_{N}^{-1}\EE{(X_{1})_{3}S_{N}^{-3}\one{S_{N}\ge N}} \to 0$.  On
$S_{N}\ge N$ so that $S_{N}^{3} \ge \max\{X_{1}^{3}, N^{3} \} $
\begin{displaymath}
\begin{split}
\frac{N}{c_{N}}\EE{ \frac{ (X_{1})_{3}}{S_{N}^{3}}\one{S_{N} \ge N}  } \le   \frac{N}{c_{N}}\EE{ \frac{ (X_{1})_{3}}{\max\{X_{1}^{3}, N^{3} \} }\one{S_{N} \ge N}  }
\end{split}
\end{displaymath}
We see, for $ 0 < \alpha \le 1$ and $\kappa > 2$ using the upper bound
in Eq~\eqref{eq:20} in Lemma~\ref{lm:ineqs} and \eqref{eq:51} for some
$C^{\prime\prime} > 0$ fixed
\begin{displaymath}
\begin{split}
 \frac{N}{c_{N}}\EE{ \frac{ (X_{1})_{3}}{\max\{X_{1}^{3}, N^{3} \} }\one{S_{N} \ge N} |E_{1} } \overline \varepsilon_{N}   \le C^{\prime\prime}  \frac{N^{2}}{N^{5-\alpha}}\sum_{k=1}^{\uN}k^{2-\alpha} &   \le \frac{C^{\prime\prime}}{3-\alpha }   \frac{(\uN + 1)^{3-\alpha}}{N^{3-\alpha}}  \\
& \le \frac{ C^{\prime\prime}}{3-\alpha}  \frac{\uN + 1 }{N} 
\end{split}
\end{displaymath}
for all $N \ge \uN + 1$ (recall $\uN/N \to 0$ by assumption).    Hence
\begin{displaymath}
\limsup_{N\to \infty} \frac{N}{c_{N}}\EE{ \frac{ (X_{1})_{3}}{S_{N}^{3}}\one{S_{N} \ge N}  }  \le   \frac{ C^{\prime\prime}}{3-\alpha}   \limsup_{N\to \infty} \frac{\uN + 1}{N} = 0
\end{displaymath}
and the lemma follows using   \eqref{eq:SchwLm6}. 
\end{proof}
\begin{remark}[The condition $\EE{X_{1}^{2}} < \infty$]
In the same manner as in the proof of
Lemma~\ref{lm:convkingmanralphaone} one checks that
$\EE{X_{1}^{2}} < \infty$ when $\kappa > 2$ and $\uN/N \to 0$; then by \citep[Proposition~7]{schweinsberg03} Case~ 
\ref{item:7} of Theorem~\ref{thm:hapl-alpha-random-one} holds.
\end{remark}
\begin{proof}[Proof of Case~\ref{item:7} of Theorem~\ref{thm:hapl-alpha-random-one}]
  Case~\ref{item:7} of Theorem~\ref{thm:hapl-alpha-random-one} is   Lemma~\ref{lm:convkingmanralphaone}.  
\end{proof}
\begin{proof}[Proof of Case~\ref{item:9} of Theorem~\ref{thm:hapl-alpha-random-one}]
Condition~\ref{eq:cntozero} is Lemma~\ref{lm:cNrandomalphaone}; condition   \ref{eq:simtozero}  is Lemma~\ref{lm:simmergersvanishrandone}; 
\ref{eq:existenceLambda} follows from
Lemma~\ref{lm:Lambdameasurealphaone} in the same way as for
Case~\ref{item:2} of Theorem~\ref{thm:haplrandomalpha}. Case~\ref{item:9}
of Theorem~\ref{thm:hapl-alpha-random-one} then follows by
Proposition~\ref{prp:conditionsconvergenceLambda}.  
\end{proof}

\subsection{Proof of Theorem~\ref{thm:convtimechangedLambdacoal}}
\label{sec:prooftimechangeLambda}

\begin{proof}[Proof of Theorem~\ref{thm:convtimechangedLambdacoal}]
Theorem~\ref{thm:convtimechangedLambdacoal}
follows from   Theorem~3 of
\citep{freund2020cannings}.  Write $X_{i}(r)$ for the number of
potential offspring  of  individual $i$ (arbitrarily labelled)  for
all $i \in [N_{r}]$.  Suppose $\zeta(N_{r})/N \gneqq 0  $ for all $r$,
and that either $\zeta(N_{r})/N \to \infty$ for all $r$, or
$\zeta(N_{r})/N \to K$ for all $r$.    Choosing $f$ and $g$ in \eqref{eq:PXiJ} 
such that  $\EE{X_{i}(r)} > 1$ for all $i \in [N_{r}] $ and $r \in
\N_{0}$  gives \citep[Lemma~1]{freund2020cannings},
\begin{displaymath}
\prob{\sum_{i=1}^{N_{r}}X_{i}(r) < N_{r-1}} < c^{N}
\end{displaymath}
for some constant $0 < c < 1$. By  Theorems~\ref{thm:haplrandomalpha}
and \ref{thm:hapl-alpha-random-one} (recall \eqref{eq:cNthmrandall})    it holds that  $c_{N}
\overset{c}{\sim} \one{\kappa > 2} N^{-1} + \one{\kappa = 2}N^{-1}\log
N$.
Then, by  \citep[Lemma~4]{freund2020cannings}, $\set{\xi^{n,N}\svigi{
\lfloor t/c_{N} \rfloor }; t \ge 0}$ converges weakly to
$\set{\xi^{n}(G(t)); t \ge 0}$, and $G(t) =
\int_{0}^{t}(v(s))^{-1}{\rm d}s$.

Define  $1_{A}(s) \equiv 1$ when $s\in A$ for some set $A$,  take
$1_{A}(s) \equiv 0$ otherwise.      We  need to check that, for the case $\kappa = 2$, it holds that
\begin{displaymath}
\begin{split}
& \sum_{r=1}^{\lfloor t/c_{N}\rfloor } \svigi{ \frac{N_{r}}{N}  }^{-1} \frac{\log N_{r} }{\log N}1_{[rc_{N}, (r+1)c_{N}) }(s)  = \sum_{r=1}^{\lfloor t/c_{N}\rfloor }\svigi{\frac{N_{r} }{N} }^{-1}\frac{\log k_{r} + \log N }{\log N}1_{ [rc_{N}, (r+1)c_{N})}(s)  \\
=  & \sum_{r}^{\lfloor t/c_{N}\rfloor }\svigi{ \frac{N_{r}}{N} }^{-1} \frac{\log k_{r}}{\log N} 1_{ [rc_{N}, (r+1)c_{N}) }(s) +  \sum_{r=1}^{\lfloor t/c_{N}\rfloor }\svigi{\frac{N_{r}}{N}}^{-1} 1_{[rc_{N}, (r+1)c_{N}) }(s) 
\end{split}
\end{displaymath}
Using \eqref{eq:7} it follows that, as $N\to \infty$,  
\begin{displaymath}
 \sum_{r=1}^{\lfloor t/c_{N}\rfloor }\svigi{\frac{N_{r}}{N}}^{-1} 1_{[rc_{N}, (r+1)c_{N}) }(s) = \svigi{ \frac{N_{\lfloor s/c_{N}\rfloor }  }{N } }^{-1} \to \svigi{v(s)}^{-1}
\end{displaymath}
Recalling  \eqref{eq:7} again  take   $h_{1}(t) \le k_{r} \le
h_{2}(t)$  with 
$h_{2}$
bounded.  It remains to check that  
\begin{displaymath}
\sum_{r=1}^{\lfloor t/c_{N}\rfloor } \svigi{ \frac{N_{r}}{N} }^{-1} \frac{\log k_{r}}{\log N}1_{[rc_{N}, (r+1)c_{N}) }(s)  \le  \frac{\log h_{2}(t) }{\log N} \sum_{r=1}^{\lfloor t/c_{N}\rfloor } \svigi{\frac{N_{r}}{N}}^{-1}1_{[rc_{N}, (r+1)c_{N})}(s) \to 0
\end{displaymath}
as $N\to \infty$
\end{proof}

\begin{remark}[Exponential growth gives a well-defined population sequence]

\breyta{When $N_{r+1} = \lfloor N_{r}(1 + \rho c_{N})
               \rfloor$ as in exponential growth it holds
               \begin{displaymath}
               \frac{N_{r+1} - N_{r}}{N} = \frac{ N_{r}(1 + \rho c_{N}) - N_{r}}{N} = \frac{N_{r}}{N}\rho c_{N} 
\end{displaymath}
Since $N_{\lfloor s/c_{N} \rfloor }/N \to v(s)$ and $c_{N}\to 0$ it
holds $\svigi{\svigi{N_{r+1} - N_{r}}/N}_{r\in \N}$ is a   null-sequence
such that  \citep[Lemma~1]{freund2020cannings} holds for any fixed
$\rho > 0$. 
}
\end{remark}

\subsection{Proof of Theorem~\ref{thm:time-changed-deltanullpoidr}}
\label{sec:proof-thm-timechanged-poissondri}

In this section we give a proof of
Theorem~\ref{thm:time-changed-deltanullpoidr}. A version of Lemma~2 of
\citep[Lemma~2]{freund2020cannings} holds for $\Xi$-coalescents (when
started from a finite number $n$ of blocks; $\Xi$-$n$-coalescents)
provided {\it (i)}  we  have a sequence $\svigi{N_{r}}_{r\in \N_{0}}$ of
population sizes where $N_{r}$ is the population size $r$ generations
into the past such that  \eqref{eq:7} holds; {\it (ii)} since
$C_{\kappa}^{N}c_{N}\overset{c}\sim 1$ it holds that  $c_{N}$ is
regularly varying so that \citep[Equation~13]{freund2020cannings}
\begin{displaymath}
M_{1}(t) \le  \frac{ c_{N_{r}}}{ c_{N}}  \le M_{2}(t)  
\end{displaymath}
for all $r \le \lfloor t/c_{N} \rfloor$ for some bounded positive
functions $M_{1},M_{2} : [0,\infty) \to (0,\infty)$; {\it (iii)}
$c_{N} \to 0$ as $N\to \infty$ by \eqref{eq:cNthmrandall}; {\it (iv)}  
 the condition
given in \citep[Equation~5]{freund2020cannings} follows as in the proof
of \citep[Theorem~3]{freund2020cannings}; {\it (v)}
$\set{\xi^{n,N}\svigi{ \lfloor t/c_{N} \rfloor}; t \ge 0}$ converges
to  a continuous-time   $\Xi$-coalescent $\set{\xi^{n}}$  as $N\to \infty$ by
Case~\ref{item:4} of  Theorem~\ref{thm:haplrandomalpha}.  
Then, we have by
\citep[Lemma~2]{freund2020cannings} that 
$\set{\xi^{n,N}\svigi{\lfloor G_{N}^{-1}(t) \rfloor }; t \ge 0 }$ converges
weakly to  $\set{\xi^{n}}$  as $N\to \infty$, where $G_{N}$ is as in
\citep[Equation~1]{freund2020cannings}.

Since   \eqref{eq:7}  and     \eqref{eq:cNthmrandall} hold, 
\citep[Lemma~4]{freund2020cannings} establishes that  the convergence
of  $\set{\xi^{n,N}\svigi{\lfloor G_{N}^{-1}(t) \rfloor }; t \ge 0 }$
to   $\set{\xi^{n}}$  is equivalent to  the convergence of
$\set{\xi^{n,N} \svigi{\lfloor t/c_{N} \rfloor }; t \ge 0  }$ to
$\set{\xi^{n}\svigi{G(t) }; t \ge 0 }$, where $G(t) =\int_{0}^{t}\svigi{v(s)}^{-1}{\rm d}s$.


\section*{Declarations}

\subsection*{Ethical Statement} The author has no conflict of interest to
declare that are relevant to the content of this article

\subsection*{Funding}

Funded in part  by DFG Priority Programme SPP 1819 `Rapid Evolutionary
Adaptation' DFG grant Projektnummer 273887127 through SPP 1819 grant
STE 325/17 to Wolfgang Stephan;  Icelandic Centre of Research (Rann\'is)
through the Icelandic Research Fund (Ranns\'oknasj\'o{\dh}ur) Grant of
Excellence no.\ 185151-051 with Einar \'Arnason, Katr\'in
Halld\'orsd\'ottir, Alison Etheridge, and Wolfgang Stephan, and DFG
SPP1819 Start-up module grants with Jere Koskela and Maite Wilke
Berenguer, and with Iulia Dahmer.

\subsection*{Data availability} The documented software (C/C++ code) developed for the numerical
results  is freely  available at \\   \url{https://github.com/eldonb/gene_genealogies_haploid_populations_sweepstakes}.

\breyta{\subsection*{Acknowledgements} The author would like to thank two
anonymous reviewers for  useful comments.}


\begin{appendices}

\counterwithin*{equation}{section}
\counterwithin*{figure}{section}
\renewcommand\theequation{\thesection\arabic{equation}}
\renewcommand\thefigure{\thesection\arabic{figure}}

\section{The Beta$(2-\beta,\beta)$-Poisson-Dirichlet$(\alpha,0)$-coalescent}\label{sec:betapoissondirichlet}

\begin{figure}[htp]
\centering
\captionsetup[subfloat]{labelfont={scriptsize,sf,md,up},textfont={scriptsize,sf}}
\subfloat[$c=1$, $\beta = 1.01$]{\includegraphics[scale=0.5]{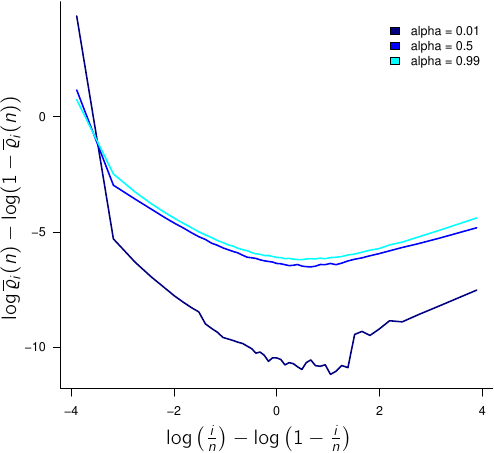}}
\subfloat[$c=100$, $\beta = 1.01$]{\label{fig:betapoissonB}\includegraphics[scale=0.5]{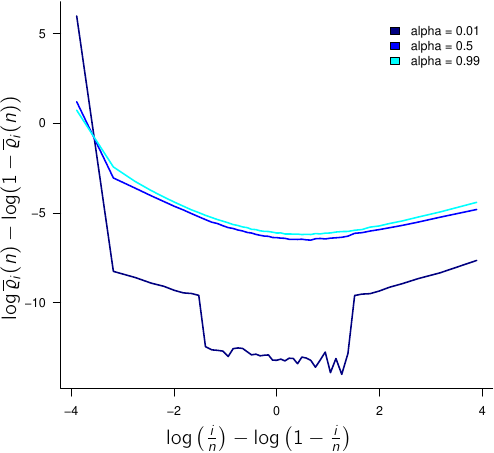}}
\subfloat[$c=10^{4}$, $\beta = 1.01$]{\label{fig:betapoissonC}\includegraphics[scale=0.5]{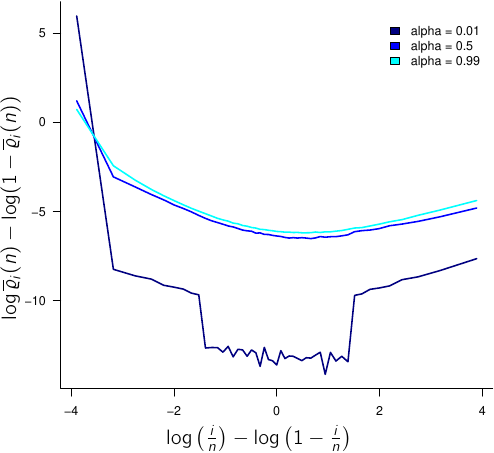}}\\
\subfloat[$c=1$, $\beta = 1.5$]{\includegraphics[scale=0.5]{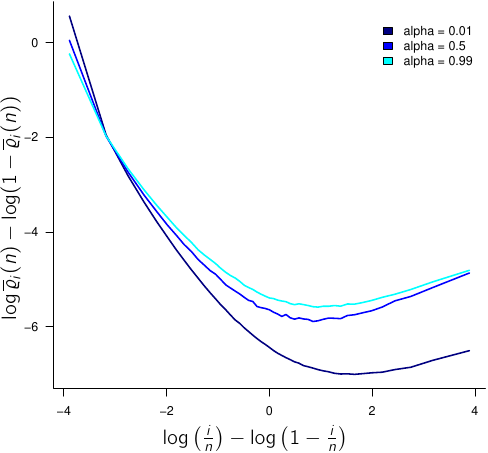}}
\subfloat[$c=100$, $\beta = 1.5$]{\label{fig:betapoissonE}\includegraphics[scale=0.5]{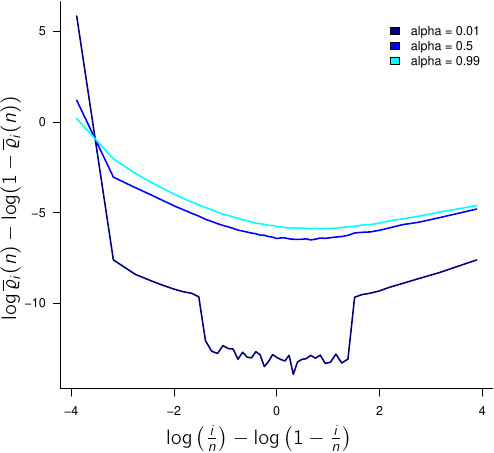}}
\subfloat[$c=10^{4}$, $\beta = 1.5$]{\label{fig:betapoissonF}\includegraphics[scale=0.5]{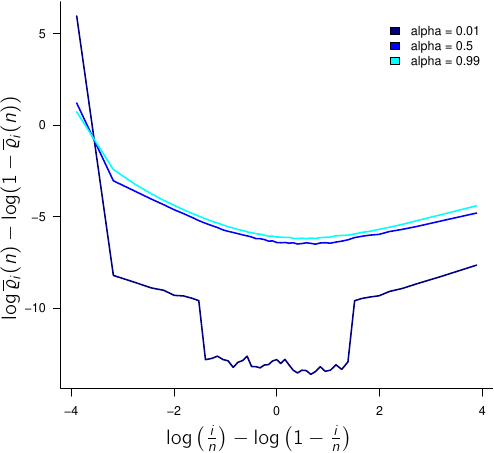}}
\caption{Examples of $\overline \varrho_{i}(n)$
\breyta{(estimates of mean relative branch lengths, recall
\eqref{eq:functionals})}  \breyta{predicted by }    the
Beta$(2-\beta,\beta)$-Poisson-Dirichlet$(\alpha,0)$-coalescent with
$n=50$, $\beta$, $\alpha$, $c$ as shown, and $m_{\infty}$
\eqref{eq:11}  approximated
with   $\svigi{2 + \svigi{1 + 2^{1-\beta}}/(\beta - 1)}/2$.   The scale of the
ordinate (y-axis) may vary between graphs (results from $10^{5}$
experiments).  \breyta{  The scale of the ordinate
 (vertical axis) may vary between graphs.} }
\label{fig:betapoissondiriA}
\end{figure}

\breyta{In this section we  briefly study the  Beta$(2-\alpha,\alpha)$-Poisson-Dirichlet$(\alpha,0)$-coalescent.}
Fix $0 < \alpha < 1$ and $1 < \beta < 2$.  Consider a population
evolving as in Definition~\ref{def:haplrandomalpha} with $\beta$
replacing $\kappa$, and such that  $\zeta(N)/N^{1/\alpha} \to \infty$
as $N\to \infty$.   The calculations leading to  Cases~\ref{item:2}
and \ref{item:4} of Theorem~\ref{thm:haplrandomalpha}, assuming
$\varepsilon_{N} \overset{c}\sim c_{N}$ where $N^{\beta - 1}c_{N}
\overset c \sim 1$ as $N\to \infty$   lead to a
coalescent with transition rates 
\begin{subequations}
\begin{align}
C_{\beta} & =  \frac{\beta f_{\beta}^{(\infty)} }{m_{\infty}^{\beta}}B(2-\beta,\beta) \label{eq:9} \\
& \underline{g_{\beta}} +  \underline{g_{\beta}}2^{1-\beta}/(\beta - 1) < m_{\infty} < \overline{f_{\beta}} + \overline{f_{\beta}}/(\beta - 1) \label{eq:11}  \\
\lambda_{n;k_{1}, \ldots, k_{r}; s} & =  \one{r=1}\frac{ C_{\beta}}{C_{\beta} + c(1-\alpha) }B(k-\beta,n-k+\beta) + \frac{cp_{n;k_{1},\ldots, k_{r};s} }{C_{\beta} + c(1-\alpha)}  \label{eq:13}
\end{align}
\end{subequations}
where $p_{n;k_{1},\ldots,k_{r};s}$ is as in \eqref{eq:34},
$B(a,b) = \int_{0}^{1}t^{a-1}(1-t)^{b-1}dt$ for $a,b > 0$, and $c > 0$
fixed (we take $\varepsilon_{N} = cN^{1-\beta}$).  We will refer to
the coalescent with transition rates as in \eqref{eq:13} as the
Beta$(2-\beta,\beta)$-Poisson-Dirichlet$(\alpha,0)$-coalescent.  It is
a multiple-merger coalescent (a $\Xi$-coalescent) without an atom at
zero. Moreover, one coalescent time unit is proportional to
$N^{\beta-1}$ generations, in contrast to the time scaling in
Theorems~\ref{thm:haplrandomalpha} and \ref{thm:hapl-alpha-random-one}
(recall \eqref{eq:cNthmrandall}). 
Figure~\ref{fig:betapoissondiriA} holds  examples of $\overline
\varrho_{i}(n)$ \eqref{eq:functionals} predicted by the
Beta$(2-\beta,\beta)$-Poisson-Dirichlet$(\alpha,0)$-coalescent. The
algorithm for sampling the coalescent is a straightforward adaption of
the one for sampling the
$\delta_{0}$-Poisson-Dirichlet$(\alpha,0)$-coalescent described in 
\S~\ref{sec:sampl-from-delta-pd}. 

\section{Examples of  $\varphi_{i}(n)$ for the
$\delta_{0}$-Beta$(\gamma,2-\alpha, \alpha)$-coalescent }
\label{sec:relat-expect-sfs}

\begin{figure}[htp]
\centering
\captionsetup[subfloat]{labelfont={scriptsize,sf,md,up},textfont={scriptsize,sf}}
\subfloat[ $\gamma = 0.05$]{\label{fig:graphexactELiELA}\includegraphics[scale=0.5]{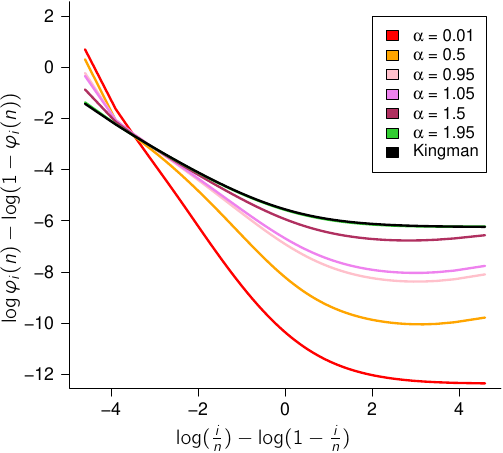}}
\subfloat[ $\gamma = 0.5$]{\label{fig:graphexactELiELB}\includegraphics[scale=0.5]{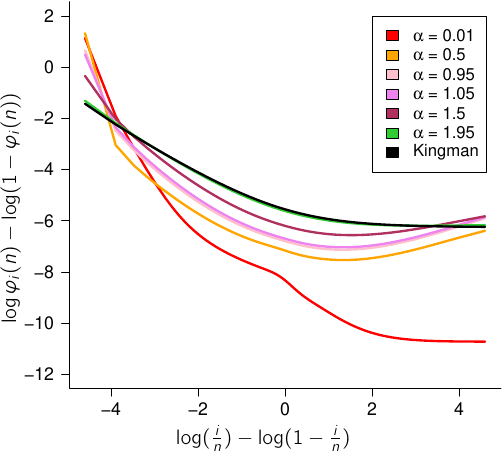}}
\subfloat[ $\gamma = 1$]{\label{fig:graphexactELiELC}\includegraphics[scale=0.5]{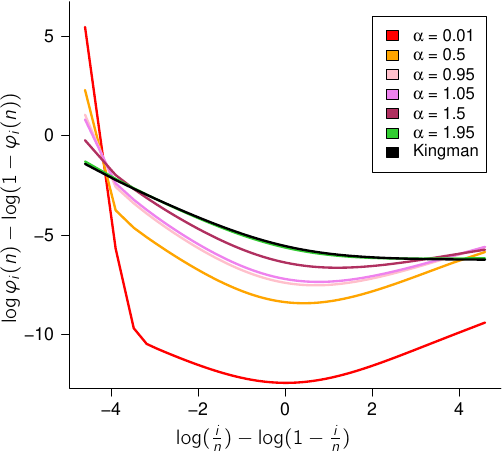}}\\
\subfloat[$\alpha = 0.01$]{\label{fig:graphexactELiELD}\includegraphics[scale=0.5]{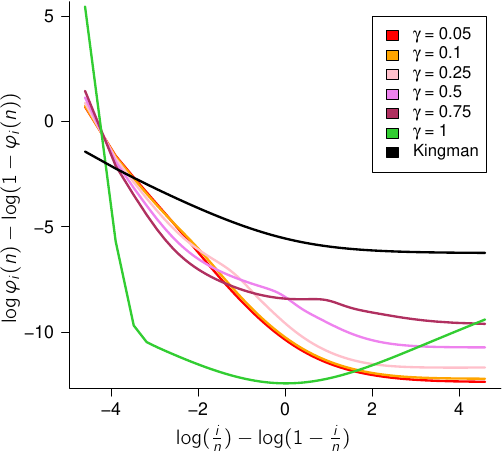}}
\subfloat[$\alpha = 0.5$]{\label{fig:graphexactELiELE}\includegraphics[scale=0.5]{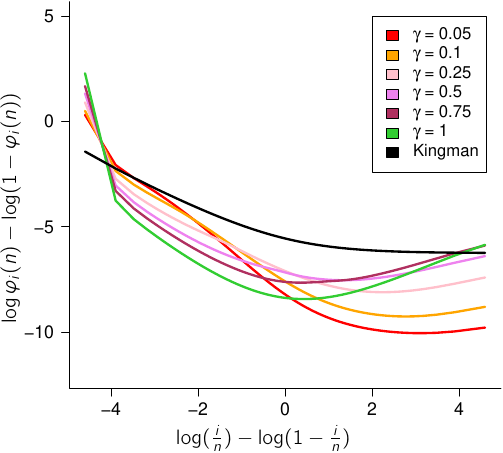}}
\subfloat[$\alpha = 0.95$]{\label{fig:graphexactELiELF}\includegraphics[scale=0.5]{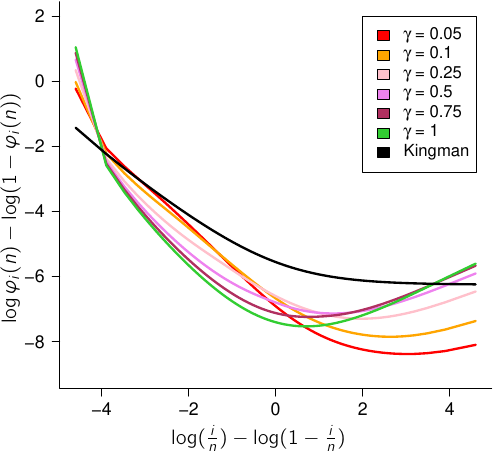}}\\
\subfloat[$\alpha = 1.05$]{\label{fig:graphexactELiELG}\includegraphics[scale=0.5]{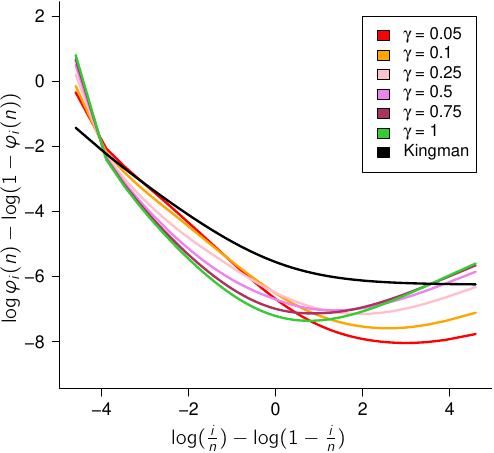}}
\subfloat[$\alpha = 1.5$]{\label{fig:graphexactELiELH}\includegraphics[scale=0.5]{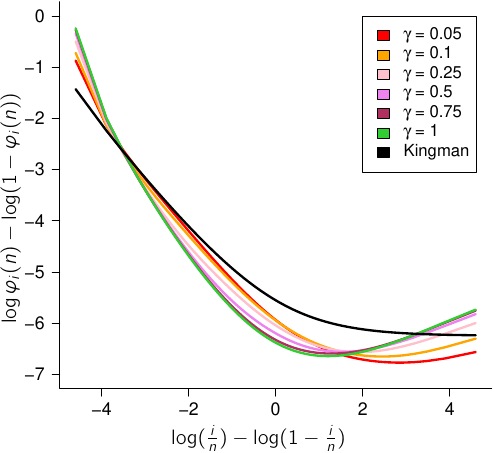}}
\subfloat[$\alpha = 1.95$]{\label{fig:graphexactELiELI}\includegraphics[scale=0.5]{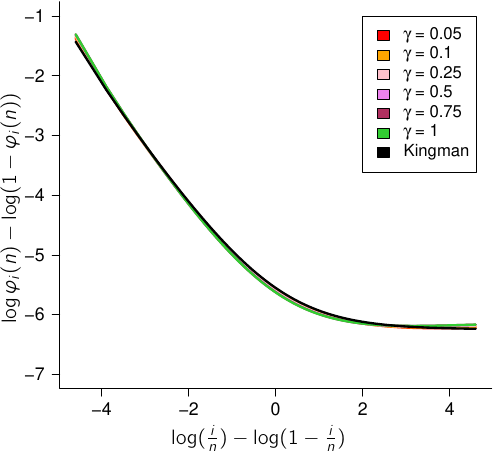}}
\caption{Relative expected branch length $\varphi_{i}(n)$
\eqref{eq:varphi} predicted by the
$\delta_{0}$-Beta$(\gamma,2-\alpha,\alpha)$-coalescent with transition rate
as in \eqref{eq:ratesarandall} for sample size $n=100$, $\alpha$ and
$\gamma$ as shown, and $c = 1$ and $\kappa = 2$.  The expected values
computed \breyta{using}  a recursion for $\EE{L_{i}(n)}$ derived for
$\Lambda$-coalescents \citep{BBE2013a}; for comparison we show the
logits of $\varphi_{i}(n) = i^{-1}/\sum_{j=1}^{n-1} j^{-1}$ predicted
by the Kingman coalescent \citep{F95}.  \breyta{  The scale of the ordinate
 (vertical axis) may vary between graphs.}  }
\label{fig:varphiA}
\end{figure}
\begin{figure}[htp]
\centering
\captionsetup[subfloat]{labelfont={scriptsize,sf,md,up},textfont={scriptsize,sf}}
\subfloat[$c = 1$, $\gamma = 0.1$]{\includegraphics[scale=0.6]{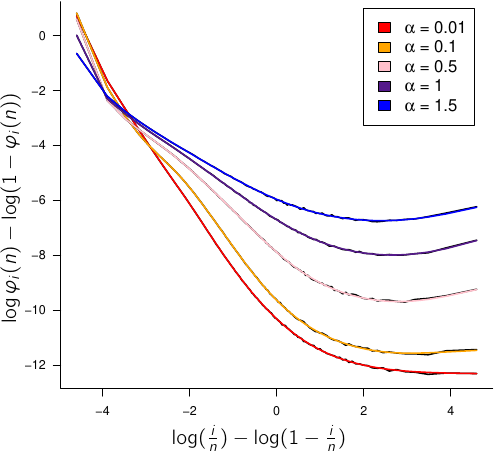}}
\subfloat[$c = \gamma = 1$]{\includegraphics[scale=0.6]{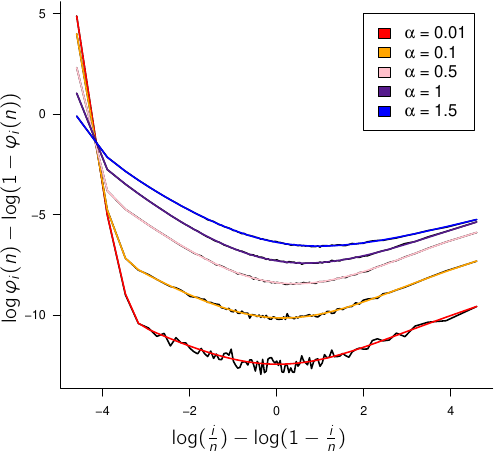}}\\
\subfloat[$c = 10$, $\gamma = 0.1$]{\includegraphics[scale=0.6]{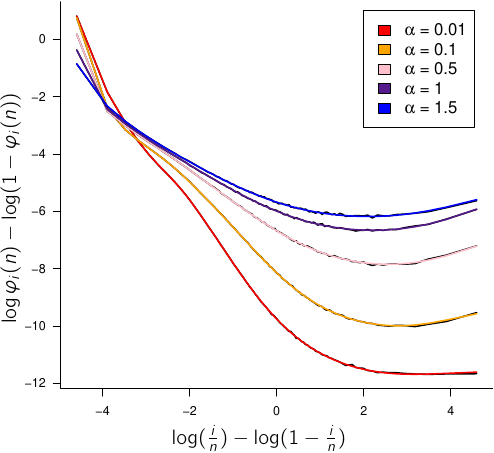}}
\subfloat[$c = 10$, $\gamma = 1$]{\includegraphics[scale=0.6]{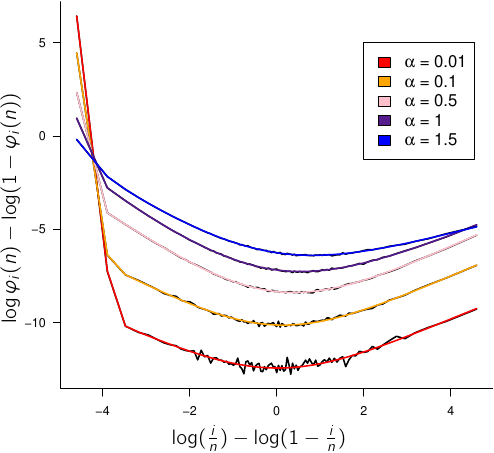}}
\caption{Comparing exact and approximated \breyta{relative
expected branch length}   $\varphi_{i}(n)$ \breyta{\eqref{eq:varphi}}  predicted by
the $\delta_{0}$-Beta$(\gamma,2-\alpha,\alpha)$ coalescent.  Checking
that approximations  of $\varphi_{i}(n)$ \eqref{eq:varphi} agree with exact
values of $\varphi_{i}(n)$ computed \breyta{using}  a recursion for
general $\Lambda$-coalescents \citep{BBE2013a}; sample size $n= 100$
with $\kappa = 2$ and $\alpha$ and $c$ and $\gamma$ as shown; the
coloured lines are exact logit values of $\varphi_{i}(n)$ and the
black lines the corresponding estimates for each value of
$\alpha$. The approximations  were obtained for $10^{5}$ experiments. \breyta{ The scale of the ordinate
 (vertical axis) may vary between graphs. } }
\label{fig:comparexactsim}
\end{figure}

In this section we give in Figure~\ref{fig:varphiA} examples of 
 $\varphi_{i}(n)$ as in
\eqref{eq:varphi}, \S~\ref{sec:numerics}, and predicted by the
$\delta_{0}$-Beta$(\gamma,2-\alpha,\alpha)$-coalescent.  In
Figures~\ref{fig:graphexactELiELA}--\ref{fig:graphexactELiELC}
$\gamma$ is fixed and we vary over $\alpha$ as shown, and in
Figures~\ref{fig:graphexactELiELD}--\ref{fig:graphexactELiELI}
$\alpha$ is fixed and we vary over $\gamma$ as shown; in each graph
the black line represents
$\varphi_{i}(n) = i^{-1}\sum_{j=1}^{n-1}j^{-1}$ as predicted by the
Kingman coalescent \citep{F95}. We see when both $\alpha$ and $\gamma$
are small enough $\varphi_{i}(n)$ is decreasing
(Figures~\ref{fig:graphexactELiELA}, \ref{fig:graphexactELiELB},
\ref{fig:graphexactELiELD}) and so will not predict a `U-shaped'
site-frequency spectrum seen as a characteristic of
$\Lambda$-coalescents \citep{freund2023interpreting}. When
$\alpha \le 1$ the distribution for the number of potential offspring
has a heavy right-tail (recall \eqref{eq:PXiJ}).  In
Figure~\ref{fig:comparexactsim} we check that simulated values of
$\varphi_{i}(n)$ predicted by the
$\delta_{0}$-Beta$(\gamma,2-\alpha,\alpha)$-coalescent match the exact
values of $\varphi_{i}(n)$ computed \breyta{using}  a recursion for
general $\Lambda$-coalescents \citep{BBE2013a}.

\section{\breyta{Further graphs comparing  $\overline \varrho_{i}^{N}(n)$ and
$\overline \varrho_{i}(n)$}}
\label{sec:sampl-discr-trees}

\begin{figure}[htp]
\centering
\captionsetup[subfloat]{labelfont={scriptsize,sf,md,up},textfont={scriptsize,sf}}
\subfloat[$\alpha = \gamma = 1$]{\label{fig:strjalNe3KBAaddA}\includegraphics[angle=0,scale=.5]{graphERnERN_KBA-crop}}
\subfloat[$\alpha = 1.25$, $\gamma = 1$]{\label{fig:strjalNe3KBB}\includegraphics[angle=0,scale=.5]{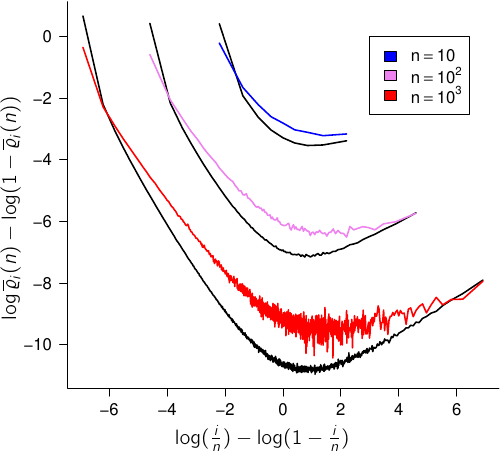}}
\subfloat[$\alpha = 1.5$, $\gamma = 1$]{\label{fig:strjalNe3KBC}\includegraphics[angle=0,scale=.5]{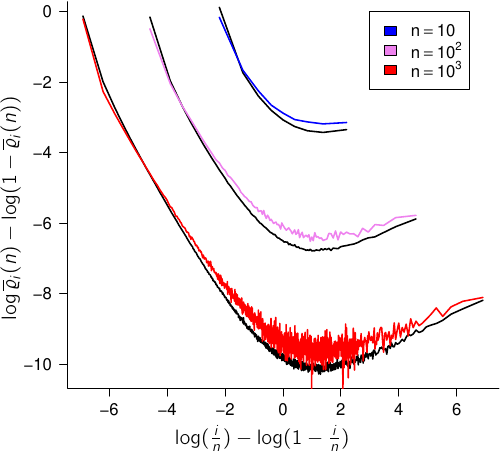}}\\
\subfloat[$\alpha = 1$, $\gamma =  \tfrac{1}{1 + \m} $]{\label{fig:strjalNe3KBD}\includegraphics[angle=0,scale=.5]{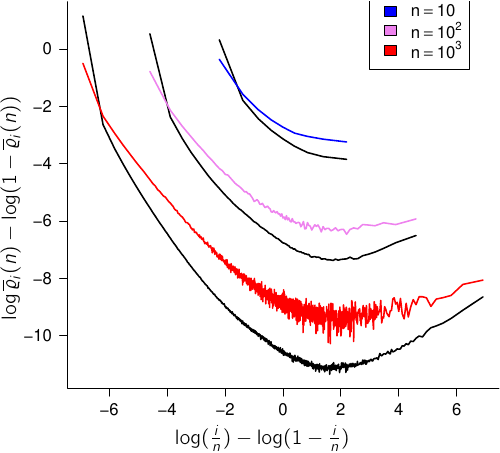}}
\subfloat[$\alpha = 1.25$, $\gamma =  \tfrac{1}{1 + \m} $]{\label{fig:strjalNe3KBE}\includegraphics[angle=0,scale=.5]{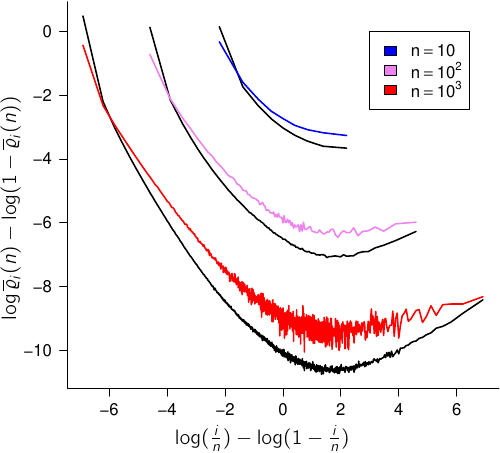}}
\subfloat[$\alpha = 1.5$, $\gamma =  \tfrac{1}{1 + \m} $]{\label{fig:strjalNe3KBFaddB}\includegraphics[angle=0,scale=.5]{graphERnERN_KBF-crop}}
\caption{\breyta{Type $A$ random
environment (Definition~\ref{def:haplrandomalpha}) }   and the 
$\delta_{0}$-Beta$(\gamma,2-\alpha,\alpha)$-coalescent.  Comparing
$\overline\varrho_{i}(n)$  (\breyta{estimates of  mean
relative branch lengths predicted by the  $\delta_{0}$-Beta$(\gamma,2-\alpha,\alpha)$-coalescent}, black lines)   and
$\overline\varrho_{i}^{N} $  (\breyta{estimates of mean
relative branch lengths predicted by $\set{\xi^{n,N}}$ when in
attraction of the   $\delta_{0}$-Beta$(\gamma,2-\alpha,\alpha)$-coalescent},  recall \eqref{eq:functionals})     when   $N = 10^{3}$,  
$\kappa = 2$, $c = 1$, and $\gamma$ as shown;     black lines are 
$\overline \varrho_{i}(n)$   for sample size $n$ as
shown with rates as in \eqref{eq:ratesarandall} in Theorem~\ref{thm:haplrandomalpha},  coloured lines  are
 $\overline \varrho_{i}^{N}(n)$ for a sample from a population of finite size $N$
evolving according to Definition~\ref{hschwpop} and Definition~\ref{def:haplrandomalpha}
with the potential offspring distributed as in \eqref{eq:40} and with
$\varepsilon_{N} = cN^{\alpha - 2}\log N$ as in \eqref{eq:33} in
Lemma~\ref{lm:cNhaplrandomall}; the case $\gamma = 1$ is compared to
$\uN = N\log N$, and $\gamma = 1/(1 + \m) $ to $\uN = N$ with
$\m$ as in \eqref{eq:57} approximating $m_{\infty}$;  $\overline \varrho_{i}^{N}(n)$
 from $10^{4}$ experiments.  The graphs  are an addition to the graphs
 in  Figure~\ref{fig:compareERiNe2A}.  \breyta{ The scale of the ordinate
 (vertical axis) may vary between graphs. }  }
\label{fig:deltanullBetaERiNaddA}
\end{figure}

\begin{figure}[htp]
\centering
\captionsetup[subfloat]{labelfont={scriptsize,sf,md,up},textfont={scriptsize,sf}}
\subfloat[$\alpha = 0.01$, $\gamma = 1$]{\includegraphics[angle=0,scale=.5]{graphERnERNarandone_KBA-crop}}
\subfloat[$\alpha = 0.25$, $\gamma = 1$]{\includegraphics[angle=0,scale=.5]{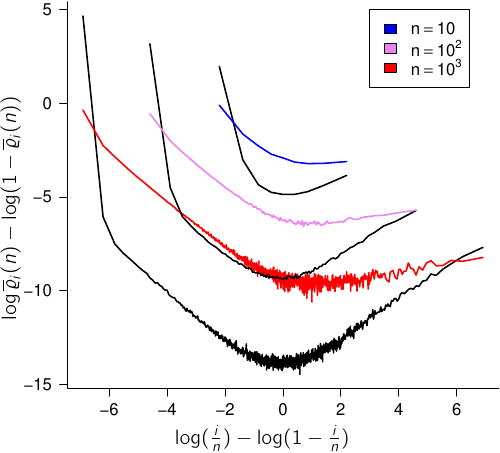}}
\subfloat[$\alpha = 0.5$, $\gamma = 1$]{\includegraphics[angle=0,scale=.5]{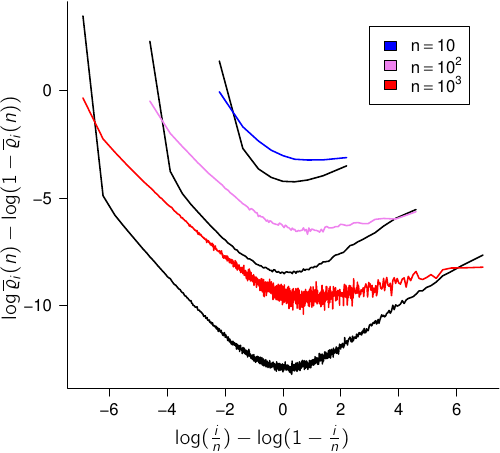}}\\
\subfloat[$\alpha = 0.01$, $\gamma = \tfrac{1}{1 + \m}$]{\includegraphics[angle=0,scale=.5]{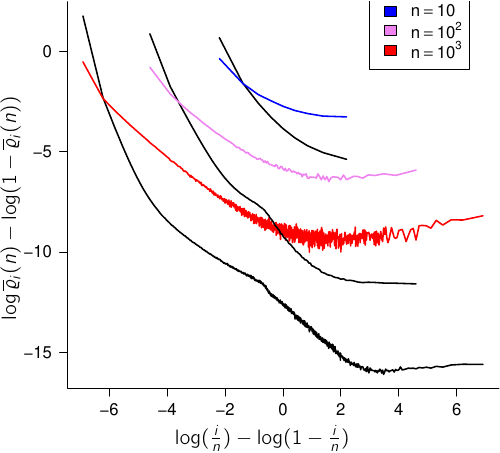}}
\subfloat[$\alpha = 0.25$, $\gamma = \tfrac{1}{1 + \m}$]{\includegraphics[angle=0,scale=.5]{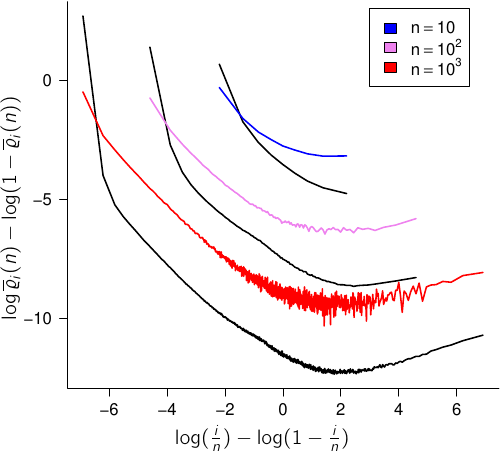}}
\subfloat[$\alpha = 0.5$, $\gamma = \tfrac{1}{1 + \m}$]{\includegraphics[angle=0,scale=.5]{graphERnERNarandone_KBF-crop}}
\caption{Definition~\ref{def:alpha-random-one} \breyta{(type
$B$ random environment)}  and the
$\delta_{0}$-Beta$(\gamma,2-\alpha,\alpha)$-coalescent. Comparing
$\overline\varrho_{i}(n)$  (\breyta{estimates of mean
relative branch lengths predicted by the  $\delta_{0}$-Beta$(\gamma,2-\alpha,\alpha)$-coalescent,} black lines)
and $\overline \varrho_{i}^{N}(n)$ (\breyta{estimates of mean
relative branch lengths predicted by $\set{\xi^{n,N}}$ when in
attraction of the
$\delta_{0}$-Beta$(\gamma,2-\alpha,\alpha)$-coalescent,}  recall \eqref{eq:functionals})  when $N=10^{3}$, $\kappa = 2$, $c = 1$, and
with $\alpha$, $\gamma$, and sampled size $n$ as shown; black lines
are  $\overline \varrho_{i}(n)$ with rates as in
\eqref{eq:ratesarandall}, coloured lines are estimates of
$\EE{R_{i}^{N}(n)}$ for a sample from a population evolving according
to Definition~\ref{hschwpop} and Definition~\ref{def:alpha-random-one}
with the potential offspring distributed as in \eqref{eq:40} and with
$\overline\varepsilon_{N} = cN^{\alpha - 1}\log N$ as in \eqref{eq:51}
in Lemma~\ref{lm:cNrandomalphaone}; the case $\gamma = 1$ is compared
to $\uN = N\log N$, and the case $\gamma = 1/(1 + \m)$
compared to $\uN = N$ with $\m$  as in
\eqref{eq:57} approximating $m_{\infty}$;  $\overline \varrho_{i}^{N}(n)$ from $10^{4}$
experiments.  The graphs are an addition to
Figure~\ref{fig:increasingnkbalessone}.  \breyta{ The scale of the ordinate
 (vertical axis) may vary between graphs. }  }
\label{fig:deltanullBetaERiNaddBB}
\end{figure}

\breyta{In this section we  compare  $\overline
\varrho_{i}^{N}(n)$ and $\varrho_{i}(n)$ (recall
\eqref{eq:functionals})    for the
$\delta_{0}$-Beta$(\gamma,2-\alpha,\alpha)$-coalescent; in
Figure~\ref{fig:deltanullBetaERiNaddA} (see also Figure~\ref{fig:compareERiNe2A})  the ancestral process $\set{\xi^{n,N}}$   is
associated with the type $A$ random environment
(Definition~\ref{def:haplrandomalpha}) such that $1 \le \alpha < 2$,
and in Figure~\ref{fig:deltanullBetaERiNaddBB} (see also Figure~\ref{fig:increasingnkbalessone})   the ancestral process
 is associated with  type $B$ random  environment
 (Definition~\ref{def:alpha-random-one}) such that $0 < \alpha < 1$.
 }

\section{\breyta{Further graphs  comparing  $\overline\varrho_{i}^{N}(n)$ and
$\overline\rho_{i}^{N}(n)$}}\label{sec:estimatequenched}

\begin{figure}[htp]
\centering
\captionsetup[subfloat]{labelfont={scriptsize,sf,md,up},textfont={scriptsize,sf}}
\subfloat[$\alpha = 0.25$,Definition~\ref{def:haplrandomalpha}]{\label{fig:qaiaddDD}\includegraphics[scale=0.6]{graphQAG-crop}}
\subfloat[$\alpha = 0.95$, Definition~\ref{def:haplrandomalpha}]{\label{fig:qaii}\includegraphics[scale=0.6]{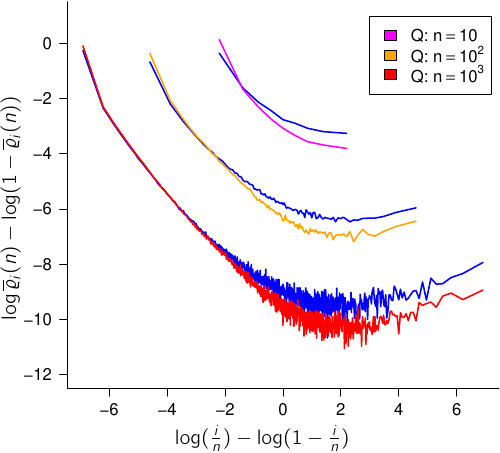}}\\
\subfloat[$\alpha = 1$, Definition~\ref{def:haplrandomalpha}]{\label{fig:qaiii}\includegraphics[scale=0.6]{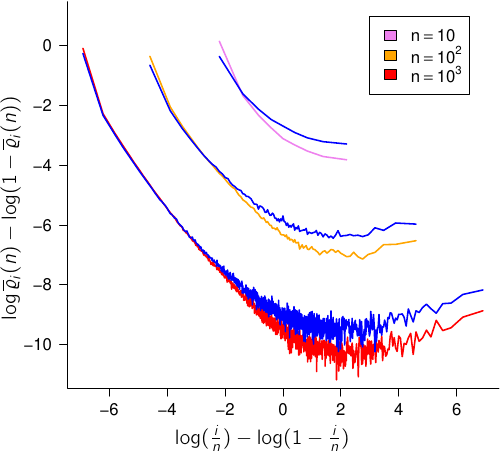}}
\subfloat[$\alpha = 1.5$, Definition~\ref{def:haplrandomalpha}]{\label{fig:qaiv}\includegraphics[scale=0.6]{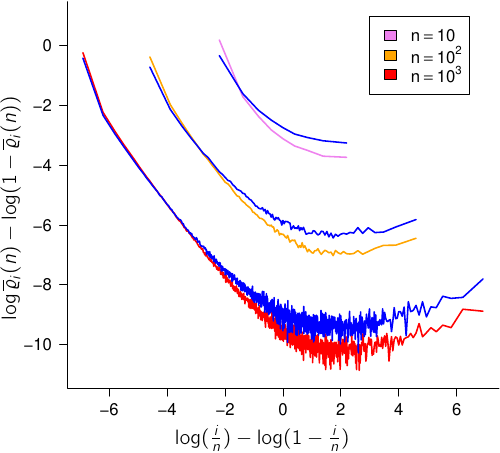}}\\
\subfloat[$\alpha = 0.01$,  Definition~\ref{def:alpha-random-one}]{\label{fig:qavaddDD}\includegraphics[scale=0.6]{graphquenchedannERiNEE-crop}}
\subfloat[$\alpha = 0.5$, Definition~\ref{def:alpha-random-one}]{\label{fig:qavi}\includegraphics[scale=0.6]{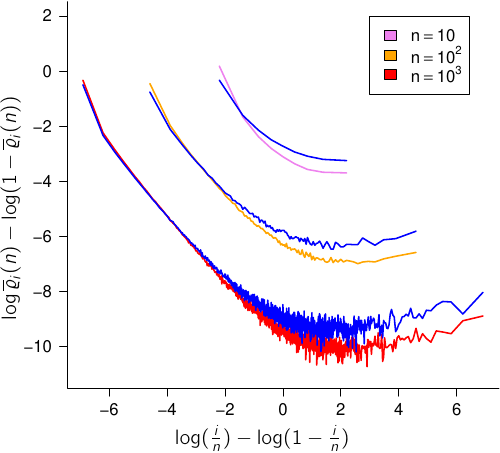}}
\caption{Quenched vs.\ annealed.  Comparing 
\breyta{ $\overline\rho_{i}^{N}(n)$}  \breyta{(estimates of mean relative
branch lengths  when conditioning on the population ancestry)}  and $\overline \varrho_{i}^{N}(n)$
(\breyta{estimates of mean relative branch lengths predicted
by  $\set{\xi^{n,N}}$,}  blue lines; recall \eqref{eq:functionals})  when the population evolves according to
Definition~\ref{hschwpop} and Definition~\ref{def:haplrandomalpha}
(\breyta{type $A$ random environment}; a,b,c,d)  and
Definition~\ref{def:alpha-random-one} (\breyta{type $B$
random environment;}   e,f); $N=10^{3}$,
$\alpha$ as shown, $\kappa = 2$, $\uN = N^{1/\alpha}\log N$ (a,b) and
$\uN = N$ (c,d,e,f);
$\varepsilon_{N} = \overline\varepsilon_{N} = 0.1$; the approximations
are from $10^{4}$ experiments.  Sections \S~\ref{sec:appr-eriNn} 
($\overline \varrho_{i}^{N}(n)$)  and
\S~\ref{sec:appr-qeewid-r_inn} ($\overline\rho_{i}^{N}(n)$) contain brief descriptions of the sampling
algorithms.  \breyta{  The scale of the ordinate
 (vertical axis) may vary between graphs.}    The graphs are an addition to
Figure~\ref{fig:quenchedannealedkpdbounded}.   }
\label{fig:QAaddDD}
\end{figure}

\breyta{In this section we compare in Figure~\ref{fig:QAaddDD}
$\overline \varrho_{i}^{N}(n)$ (the annealed mean relative branch lengths)   and
 $\overline \rho_{i}^{N}(n)$ (the quenched mean relative branch lengths,  recall~\eqref{eq:functionals}).
 Figure~\ref{fig:QAaddDD} is  an addition to
Figure~\ref{fig:quenchedannealedkpdbounded}.    }

\section{\breyta{Further graphs for the   $\delta_{0}$-Poisson-Dirichlet$(\alpha,
0)$-coalescent}}\label{sec:samplingKPD}

\begin{figure}[htp]
\centering
\subfloat[$c = 1$]{\includegraphics[angle=0,scale=.6]{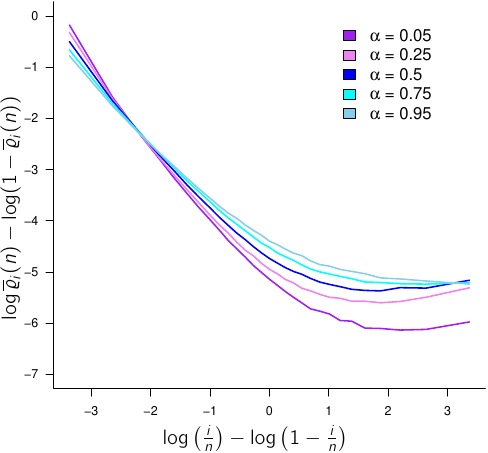}}
\subfloat[$c = 100$]{\includegraphics[angle=0,scale=.6]{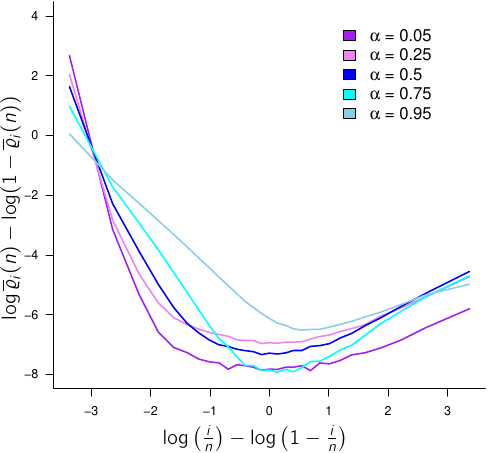}}\\
\subfloat[$c = 1$]{\label{fig:dpdERiNcompareA}\includegraphics[angle=0,scale=.6]{pdpERiNcompare_n30_c1-crop}}
\subfloat[$c = 100$]{\label{fig:dpdERiNcompareB}\includegraphics[angle=0,scale=.6]{pdpERiNcompare_n30_c100-crop}}
\caption{The $\delta_{0}$-Poisson-Dirichlet$(\alpha,0)$-coalescent.
Approximations  $(\overline \varrho_{i}(n))$ of $\EE{R_{i}(n)}$ 
(\breyta{approximations of mean relative branch lengths, 
 recall \eqref{eq:functionals};} lines) predicted by the
$\delta_{0}$-Poisson-Dirichlet$(\alpha,0)$-coalescent compared to
$\overline\varrho_{i}^{N}(n)$ (\breyta{estimates of mean
relative branch lengths predicted by $\set{\xi^{n,N}}$ when in
attraction of the  $\delta_{0}$-Poisson-Dirichlet$(\alpha,0)$-coalescent;}  circles)  when the population evolves according
to Definition~\ref{def:haplrandomalpha} \breyta{(type $A$
random environment)}   with
$\uN = \infty$, $\varepsilon = c(\log N)/N$ for $c$ as shown,
$N=3000$, $\alpha$ as shown and $\kappa = 2$. The scale of the
ordinate (y-axis) may vary between the graphs.   
See \S~\ref{sec:sampl-from-delta-pd}  for an algorithm
for sampling from the
$\delta_{0}$-Poisson-Dirichlet$(\alpha,0)$-coalescent.   \breyta{ The scale of the ordinate
 (vertical axis) may vary between graphs. }   The graphs are an addition to   Figure~\ref{fig:kpdbounded}.  }
\label{fig:dpdaddCC}
\end{figure}

\begin{figure}[htp]
\centering
\captionsetup[subfloat]{labelfont={scriptsize,sf,md,up},textfont={scriptsize,sf}}
\subfloat[ $c = 1$]{\includegraphics[scale=0.6]{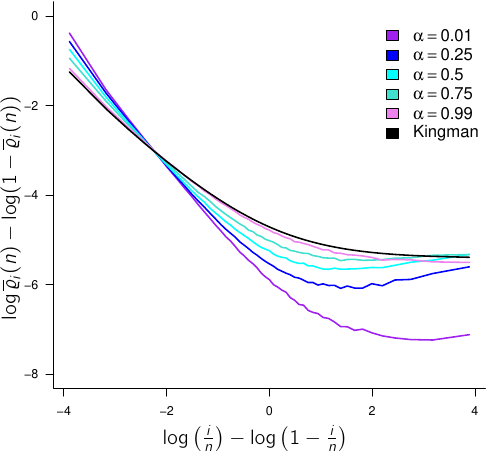}}
\subfloat[ $c = 10^{1}$]{\includegraphics[scale=0.6]{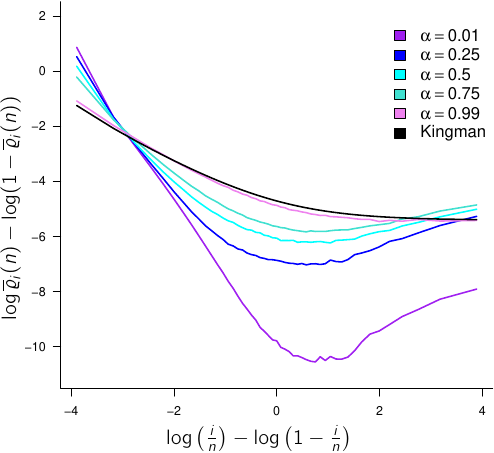}}
\subfloat[ $c = 10^{2}$]{\label{fig:dpdexamplesc}\includegraphics[scale=0.6]{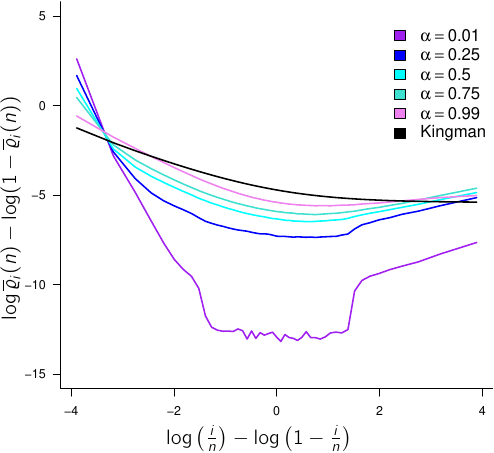}}
\caption{The-$\delta_{0}$-Poisson-Dirichlet$(\alpha,0)$ coalescent.     Examples of   \breyta{$\overline
\varrho_{i}(n)$ (approximations  of mean relative branch lengths)}  \eqref{eq:functionals}    of $\EE{R_{i}(n)}$
\breyta{predicted by}   the
$\delta_{0}$-Poisson-Dirichlet$(\alpha,0)$ coalescent (see
\S~\ref{sec:samplingKPD})   for $n=50$,
$\kappa = 2$, and $c$ and $\alpha$ as shown. The black lines are the
Kingman prediction  $i^{-1}/\sum_{j=1}^{n-1}j^{-1}$  \citep{F95}.
Results from $10^{6}$ experiments.  \breyta{ The scale of the ordinate
 (vertical axis) may vary between graphs.  The graphs are an addition
 to Figure~\ref{fig:kpdbounded}. }    }
\label{fig:dpdexamples}
\end{figure}

\breyta{In this section  we  consider in    Figure~\ref{fig:dpdaddCC} and
Figure~\ref{fig:dpdexamples}  $\overline \varrho_{i}(n)$ (annealed
mean relative branch lengths, recall ~\eqref{eq:functionals})  for     the
$\delta_{0}$-Poisson-Dirichlet$(\alpha,0)$-coalescent.
In addition, in  graphs~\ref{fig:dpdERiNcompareA} and \ref{fig:dpdERiNcompareB}  we
compare  $\overline \varrho_{i}^{N}(n)$ and  $\overline
\varrho_{i}(n)$.  The graphs   are an addition to
Figure~\ref{fig:kpdbounded}.  
 }


\section{Further graphs for time-varying coalescents}\label{sec:further-graphs}

\begin{figure}[htp]
\centering
\captionsetup[subfloat]{labelfont={scriptsize,sf,md,up},textfont={scriptsize,sf}}
\subfloat[$\alpha = 0.01$, $\gamma = 1$]{\includegraphics[scale=0.5]{deltanull_incbeta_expgrowth_graphD-crop}}
\subfloat[$\alpha = 1$, $\gamma = 1$]{\includegraphics[scale=0.5]{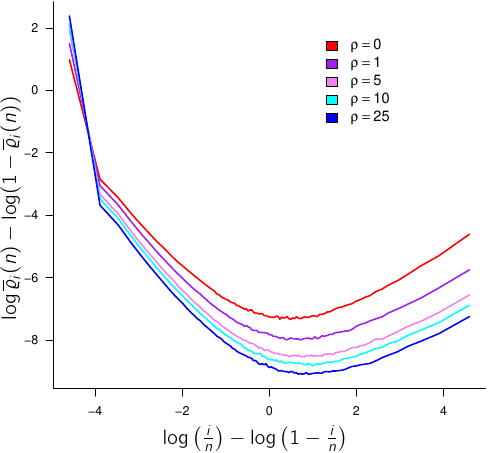}}
\subfloat[$\alpha = 1.5$, $\gamma = 1$]{\includegraphics[scale=0.5]{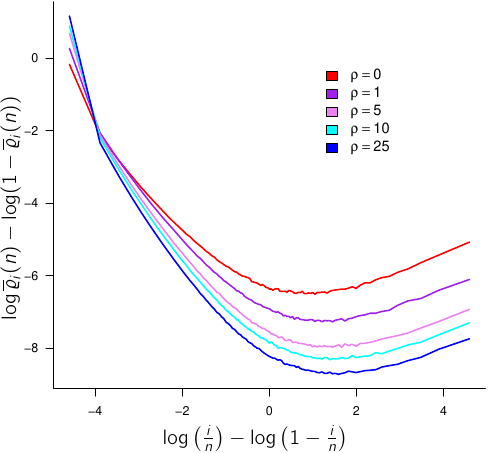}}\\
\subfloat[$\alpha = 0.01$, $\gamma = 0.5$]{\includegraphics[scale=0.5]{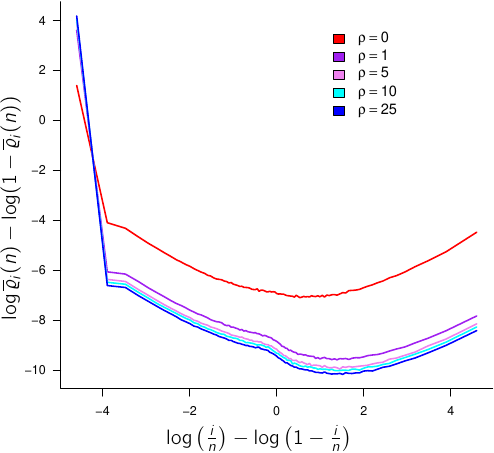}}
\subfloat[$\alpha = 1$, $\gamma = 0.5$]{\includegraphics[scale=0.5]{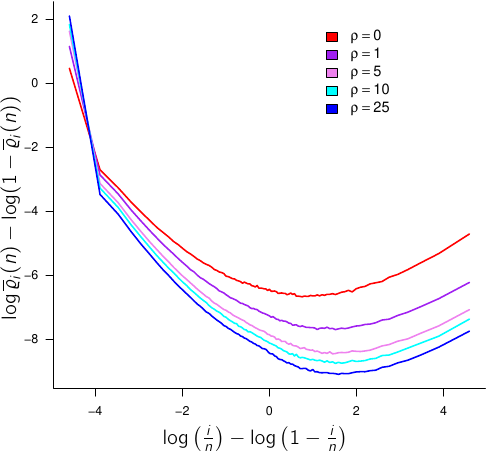}}
\subfloat[$\alpha = 1.5$, $\gamma = 0.5$]{\includegraphics[scale=0.5]{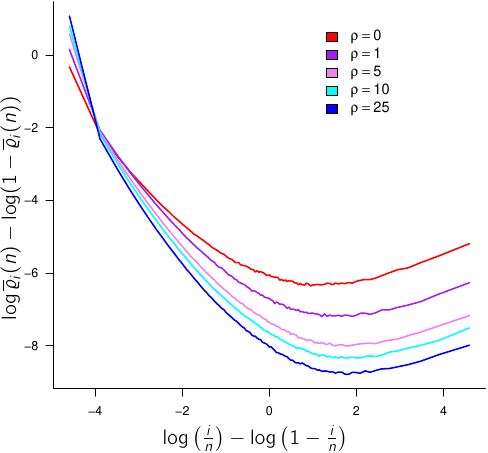}}\\
\subfloat[$\alpha = 0.01$, $\gamma = 0.1$]{\includegraphics[scale=0.5]{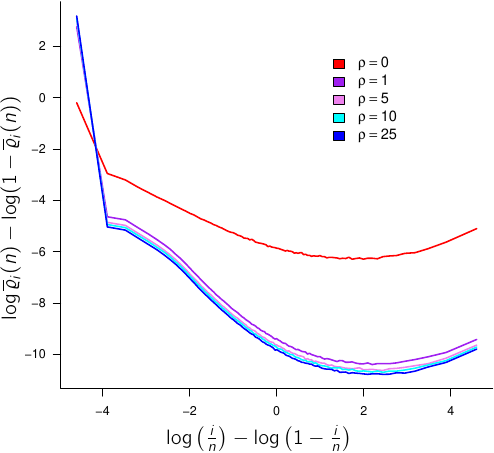}}
\subfloat[$\alpha = 1$, $\gamma = 0.1$]{\includegraphics[scale=0.5]{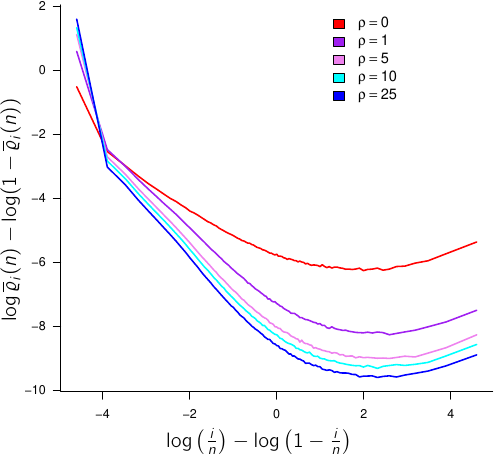}}
\subfloat[$\alpha = 1.5$, $\gamma = 0.1$]{\includegraphics[scale=0.5]{deltanull_incbeta_expgrowth_graphL-crop}}
\caption{Approximations  $\overline \varrho_{i}(n)$
(\breyta{estimates of mean relative branch lengths;}   recall
\eqref{eq:functionals}) \breyta{predicted by }  the time-changed   $\delta_{0}$-Beta$(\gamma,2-\alpha,\alpha)$-coalescent with
time-change $G(t) = \int_{0}^{t}e^{\rho s}{\rm d}s =  \one{\rho > 0}\svigi{e^{\rho
t} - 1}/\rho + \one{\rho = 0}t$  with  $\rho$ as shown,   $\kappa
= 2$, $c=1$, $n=100$.   The scale of the ordinate ($y$-axis) may vary between the
graphs; results from $10^{5}$ experiments.    \breyta{ The scale of the ordinate
 (vertical axis) may vary between graphs. }    The graphs are an addition
to Figure~\ref{fig:deltanullincbetaexpgrowth}.}
\label{fig:addbetaEE}
\end{figure}

\begin{figure}[htp]
\centering
\captionsetup[subfloat]{labelfont={scriptsize,sf,md,up},textfont={scriptsize,sf}}
\subfloat[$\alpha = 0.01$, $c=1$]{\includegraphics[scale=0.5]{deltanull_poissondirichlet_expgrowth_graphA-crop.pdf}}
\subfloat[$\alpha = 0.5$, $c=1$]{\includegraphics[scale=0.5]{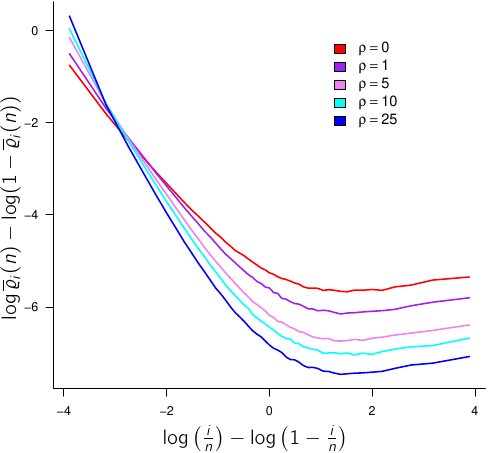}}
\subfloat[$\alpha = 0.99$, $c=1$]{\includegraphics[scale=0.5]{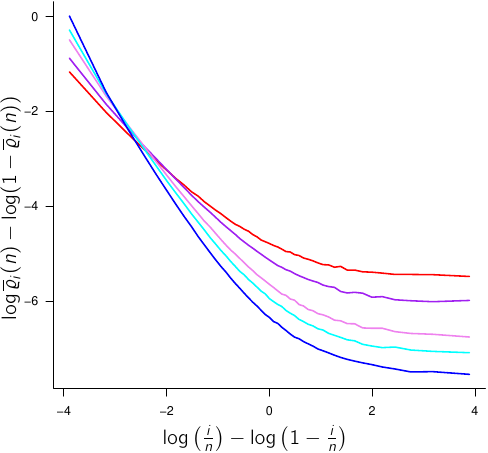}}\\
\subfloat[$\alpha = 0.01$, $c=10$]{\includegraphics[scale=0.5]{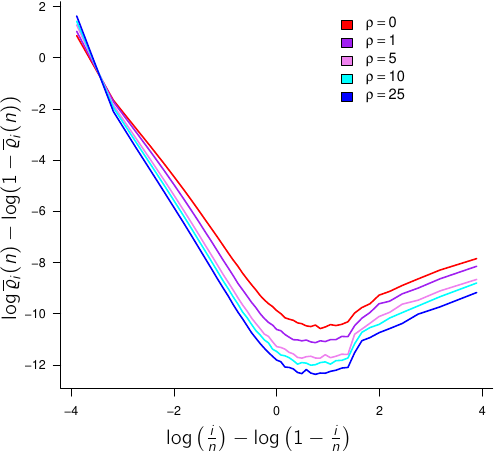}}
\subfloat[$\alpha = 0.5$, $c=10$]{\includegraphics[scale=0.5]{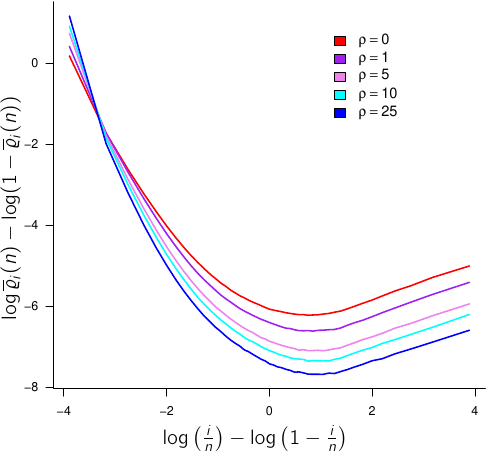}}
\subfloat[$\alpha = 0.99$, $c=10$]{\includegraphics[scale=0.5]{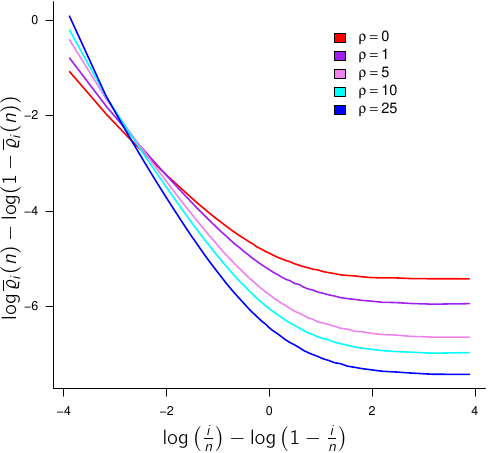}}\\
\subfloat[$\alpha = 0.01$, $c=100$]{\label{fig:adddpdFFg}\includegraphics[scale=0.5]{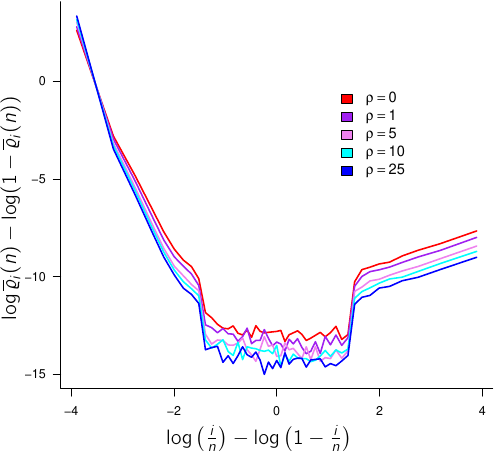}}
\subfloat[$\alpha = 0.5$, $c=100$]{\includegraphics[scale=0.5]{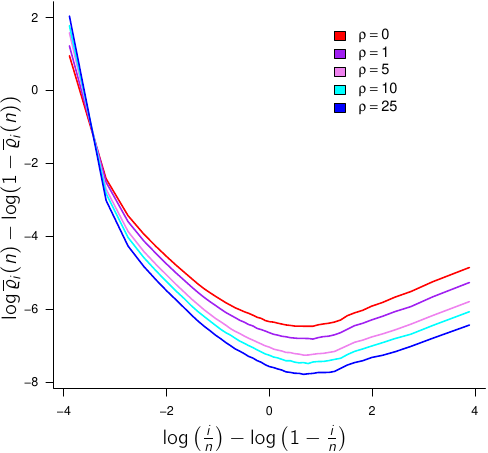}}
\subfloat[$\alpha = 0.99$, $c=100$]{\includegraphics[scale=0.5]{deltanull_poissondirichlet_expgrowth_graphF-crop.pdf}}
\caption{Approximations $\overline \varrho_{i}(n)$
(\breyta{estimates of mean relative branch lengths;} recall
\eqref{eq:functionals})  \breyta{predicted by}  the
time-changed  $\delta_{0}$-Poisson-Dirichlet$(\alpha,0)$-coalescent
with time-change   $G(t) = \int_{0}^{t}e^{\rho s}{\rm d}s =  \one{\rho > 0}\svigi{e^{\rho
t} - 1}/\rho + \one{\rho = 0}t$  with  $\rho$ as shown,  
 $n=50$, $\kappa = 2$; 
results from   $10^{6}$ experiments.  \breyta{ The scale of the ordinate
 (vertical axis) may vary between graphs. }  The graphs are an addition to Figure~\ref{fig:time-changed-deltanull-poissondiri}}
\label{fig:adddpdFF}
\end{figure}

\breyta{This section contains  further  graphs   showing the
joint effects  of  population growth and  sweepstakes reproduction on
the site-frequency spectrum.  Figure~\ref{fig:addbetaEE} contains
graphs  for the  time-changed  
$\delta_{0}$-Beta$(\gamma,2-\alpha,\alpha)$-coalescent,  and
Figure~\ref{fig:adddpdFF} for the time-changed
$\delta_{0}$-Poisson-Dirichlet$(\alpha,0)$-coalescent.
Figure~\ref{fig:addbetaEE} is an addition to
Figure~\ref{fig:deltanullincbetaexpgrowth}, and
Figure~\ref{fig:adddpdFF} to Figure~\ref{fig:time-changed-deltanull-poissondiri}. }

\end{appendices}

\clearpage
\pagebreak
\newpage

\bibliographystyle{plainnat}%
\bibliography{refs}%

@article{eldon24:_genediploid,
      title={Gene genealogies in diploid populations evolving according to sweepstakes reproduction}, 
      author={Bjarki Eldon},
      year={2026},
      eprint={2601.10364},
      journal=arXiv,
      archivePrefix={arXiv},
      primaryClass={q-bio.PE},
      doi={10.48550/arXiv.2601.10364},
      url={https://arxiv.org/abs/2601.10364}
}

@article{https://doi.org/10.1111/mec.16774,
author = {Barfield, Sarah and Davies, Sarah W. and Matz, Mikhail V.},
title = {Evidence of sweepstakes reproductive success in a broadcast-spawning coral and its implications for coral metapopulation persistence},
journal = {Molecular Ecology},
volume = {32},
number = {3},
pages = {696-702},
doi = {https://doi.org/10.1111/mec.16774},
url = {https://onlinelibrary.wiley.com/doi/abs/10.1111/mec.16774},
eprint = {https://onlinelibrary.wiley.com/doi/pdf/10.1111/mec.16774},
year = {2023}
}

@article{Rosen2018,
  title = {Geometry of the Sample Frequency Spectrum and the Perils of Demographic Inference},
  volume = {210},
  ISSN = {1943-2631},
  url = {http://dx.doi.org/10.1534/genetics.118.300733},
  DOI = {10.1534/genetics.118.300733},
  number = {2},
  journal = {Genetics},
  publisher = {Oxford University Press (OUP)},
  author = {Rosen,  Zvi and Bhaskar,  Anand and Roch,  Sebastien and Song,  Yun S},
  year = {2018},
  month = jul,
  pages = {665–682}
}

@article{Bhaskar2014,
  title = {Descartes’ rule of signs and the identifiability of population demographic models from genomic variation data},
  volume = {42},
  ISSN = {0090-5364},
  url = {http://dx.doi.org/10.1214/14-AOS1264},
  DOI = {10.1214/14-aos1264},
  number = {6},
  journal = {The Annals of Statistics},
  publisher = {Institute of Mathematical Statistics},
  author = {Bhaskar,  Anand and Song,  Yun S.},
  year = {2014},
  month = dec 
}

@article{Myers2008,
  title = {Can one learn history from the allelic spectrum?},
  volume = {73},
  ISSN = {0040-5809},
  url = {http://dx.doi.org/10.1016/j.tpb.2008.01.001},
  DOI = {10.1016/j.tpb.2008.01.001},
  number = {3},
  journal = {Theoretical Population Biology},
  publisher = {Elsevier BV},
  author = {Myers,  Simon and Fefferman,  Charles and Patterson,  Nick},
  year = {2008},
  month = may,
  pages = {342–348}
}

@article{Donnelly1995,
  title = {COALESCENTS AND GENEALOGICAL STRUCTURE UNDER NEUTRALITY},
  volume = {29},
  ISSN = {1545-2948},
  url = {http://dx.doi.org/10.1146/annurev.ge.29.120195.002153},
  DOI = {10.1146/annurev.ge.29.120195.002153},
  number = {1},
  journal = {Annual Review of Genetics},
  publisher = {Annual Reviews},
  author = {Donnelly,  Peter and Tavar{\'e},  Simon},
  year = {1995},
  month = dec,
  pages = {401–421}
}

@article{Gnedin2014asymptotics,
  title = {On Asymptotics of the Beta Coalescents},
  volume = {46},
  ISSN = {1475-6064},
  url = {http://dx.doi.org/10.1239/aap/1401369704},
  DOI = {10.1239/aap/1401369704},
  number = {2},
  journal = {Advances in Applied Probability},
  publisher = {Cambridge University Press (CUP)},
  author = {Gnedin,  Alexander and Iksanov,  Alexander and Marynych,  Alexander and M\"{o}hle,  Martin},
  year = {2014},
  month = jun,
  pages = {496–515}
}

@Article{ethier93:_flemin_viot,
  author = 	 {S Ethier and T Kurtz},
  title = 	 {{Fleming-Viot} processes in population genetics},
  journal = 	 {SIAM J Control Optim},
  year = 	 1993,
  volume = 	 31,
  doi={https://doi.org/10.1137/0331019},
  pages = 	 {345--86}}

@Article{fleming79:_some_markov,
  author = 	 {WH Fleming and M  Viot},
  title = 	 {Some measure-valued {Markov} processes in population genetics theory},
  journal = 	 {Indiana University Mathematics Journal},
  year = 	 1979,
  volume = 	 28,
  url = {https://www.jstor.org/stable/24892583},
  pages = 	 {817--43}}

@Misc{tange11:_gnu_paral,
  author = 	 {O Tange},
  title = 	 {{GNU} Parallel -- The Command-Line Power Tool},
  howpublished = {The USENIX Magazine},
  year = 	 2011}

@article{Gnedin2014,
  doi = {10.1239/jap/1417528464},
  year = {2014},
  month = dec,
  publisher = {Cambridge University Press ({CUP})},
  volume = {51},
  number = {A},
  pages = {23--40},
  author = {Alexander Gnedin and Alexander Iksanov and Alexander Marynych},
  title = {{$\Lambda$}-coalescents: a survey},
  journal = {Journal of Applied Probability}
}

@article{freund2021impact,
  title={The impact of genetic diversity statistics on model selection between coalescents},
  author={Freund, Fabian and Siri-J{\'e}gousse, Arno},
  journal={Computational Statistics \& Data Analysis},
  volume={156},
  pages={107055},
  year={2021},
  doi={https://doi.org/10.1016/j.csda.2020.107055},
  publisher={Elsevier}
}

@article{Melfi2018b,
  doi = {10.1016/j.tpb.2018.09.005},
  year = {2018},
  month = dec,
  publisher = {Elsevier {BV}},
  volume = {124},
  pages = {81--92},
  author = {Andrew Melfi and Divakar Viswanath},
  title = {The {Wright}{\textendash}{Fisher} site frequency spectrum as a perturbation of the coalescent's},
  journal = {Theoretical Population Biology}
}

@article{Melfi2018,
  doi = {10.1016/j.tpb.2018.04.001},
  year = {2018},
  month = may,
  publisher = {Elsevier {BV}},
  volume = {121},
  pages = {60--71},
  author = {Andrew Melfi and Divakar Viswanath},
  title = {Single and simultaneous binary mergers in {Wright-Fisher} genealogies},
  journal = {Theoretical Population Biology}
}

@article{CHRISTIE2010,
  title = {Self‐recruitment and sweepstakes reproduction amid extensive gene flow in a coral‐reef fish},
  volume = {19},
  ISSN = {1365-294X},
  url = {http://dx.doi.org/10.1111/j.1365-294X.2010.04524.x},
  DOI = {10.1111/j.1365-294x.2010.04524.x},
  number = {5},
  journal = {Molecular Ecology},
  publisher = {Wiley},
  author = {CHRISTIE,  MARK R. and JOHNSON,  DARREN W. and STALLINGS,  CHRISTOPHER D. and HIXON,  MARK A.},
  year = {2010},
  month = Feb,
  pages = {1042–1057}
}

@article{Baumdicker2021,
  doi = {10.1093/genetics/iyab229},
  url = {https://doi.org/10.1093/genetics/iyab229},
  year = {2021},
  month = dec,
  publisher = {Oxford University Press ({OUP})},
  author = {Franz Baumdicker and Gertjan Bisschop and Daniel Goldstein
                  and Graham Gower and Aaron P Ragsdale and Georgia
                  Tsambos and Sha Zhu and Bjarki Eldon and E Castedo
                  Ellerman and Jared G Galloway and Ariella L
                  Gladstein and Gregor Gorjanc and Bing Guo and Ben
                  Jeffery and Warren W Kretzschmar and Konrad Lohse
                  and Michael Matschiner and Dominic Nelson and
                  Nathaniel S Pope and Consuelo D Quinto-Cort{\'{e}}s
                  and Murillo F Rodrigues and Kumar Saunack and
                  Thibaut Sellinger and Kevin Thornton and Hugo van
                  Kemenade and Anthony W Wohns and Yan Wong and Simon
                  Gravel and Andrew D Kern and Jere Koskela and Peter
                  L Ralph and Jerome Kelleher},
  title = {Efficient ancestry and mutation simulation with msprime $1.0$},
  journal = {Genetics}
}

@article{freund2020cannings,
  title={Cannings models, population size changes and multiple-merger coalescents},
  author={Freund, Fabian},
  journal={Journal of mathematical biology},
  volume={80},
  number={5},
  pages={1497--1521},
  year={2020},
  doi={https://doi.org/10.1007/s00285-020-01470-5},
  publisher={Springer}
}

@Book{athreya06:_measur,
  author = 	 {KB Athreya and SN Lahiri},
  title = 	 {Measure theory and probability theory},
  publisher = 	 {Springer},
  doi = {https://doi.org/10.1007/978-0-387-35434-7},
  year = 	 2006}

@Book{hardy02,
  author = 	 {GH Hardy},
  title = 	 {A course of pure mathematics},
  publisher = 	 {Cambridge University Press},
  year = 	 2002,
  address = 	 {Cambridge, UK},
  doi = {https://doi.org/10.1017/CBO9780511989469},
  edition = 	 {Tenth}}

@article{Diamantidis2024,
  title = {Bursts of coalescence within population pedigrees whenever big families occur},
  volume = {227},
  ISSN = {1943-2631},
  url = {http://dx.doi.org/10.1093/genetics/iyae030},
  DOI = {10.1093/genetics/iyae030},
  number = {1},
  journal = {GENETICS},
  publisher = {Oxford University Press (OUP)},
  author = {Diamantidis,  Dimitrios and Fan,  Wai-Tong (Louis) and Birkner,  Matthias and Wakeley,  John},
  editor = {Jain,  K},
  year = {2024},
  month = feb 
}

@article{PANOV2017379,
title = {Limit theorems for sums of random variables with mixture distribution},
journal = {Statistics \& Probability Letters},
volume = {129},
pages = {379-386},
year = {2017},
issn = {0167-7152},
doi = {https://doi.org/10.1016/j.spl.2017.06.017},
url = {https://www.sciencedirect.com/science/article/pii/S0167715217302213},
author = {Vladimir Panov}
}

@article{10.1214/aoms/1177700291,
author = {Bengt von Bahr and Carl-Gustav Esseen},
title = {{Inequalities for the $r$th Absolute Moment of a Sum of Random Variables, $1 \leqq r \leqq 2$}},
volume = {36},
journal = {The Annals of Mathematical Statistics},
number = {1},
publisher = {Institute of Mathematical Statistics},
pages = {299 -- 303},
year = {1965},
doi = {10.1214/aoms/1177700291},
URL = {https://doi.org/10.1214/aoms/1177700291}
}

@article{freund2023interpreting,
  doi = {10.1371/journal.pgen.1010677},
  url = {https://doi.org/10.1371/journal.pgen.1010677},
  year = {2023},
  month = mar,
  publisher = {Public Library of Science ({PLoS})},
  volume = {19},
  number = {3},
  pages = {e1010677},
  author = {Fabian Freund and Elise Kerdoncuff and Sebastian
                  Matuszewski and Marguerite Lapierre and Marcel
                  Hildebrandt and Jeffrey D. Jensen and Luca Ferretti
                  and Amaury Lambert and Timothy B. Sackton and
                  Guillaume Achaz},
  editor = {Lindi Wahl},
  title = {Interpreting the pervasive observation of {U-shaped} Site Frequency Spectra},
  journal = {{PLOS} Genetics}
}

@book{bertoin2006random,
  title={Random fragmentation and coagulation processes},
  author={Bertoin, Jean},
  volume={102},
  year={2006},
  doi={https://doi.org/10.1017/CBO9780511617768},
  publisher={Cambridge University Press}
}

@article{Kingman1975,
  doi = {10.1111/j.2517-6161.1975.tb01024.x},
  url = {https://doi.org/10.1111/j.2517-6161.1975.tb01024.x},
  year = {1975},
  month = sep,
  publisher = {Wiley},
  volume = {37},
  number = {1},
  pages = {1--15},
  author = {J. F. C. Kingman},
  title = {Random Discrete Distributions},
  journal = {Journal of the Royal Statistical Society: Series B (Methodological)}
}

@article{pitman1995exchangeable,
  title={Exchangeable and partially exchangeable random partitions},
  author={Pitman, Jim},
  journal={Probability theory and related fields},
  volume={102},
  number={2},
  pages={145--158},
  year={1995},
  doi={10.1007/BF01213386},
  publisher={Springer}
}

@article{Etemadi1981,
  doi = {10.1007/bf01013465},
  year = {1981},
  month = feb,
  publisher = {Springer Science and Business Media {LLC}},
  volume = {55},
  number = {1},
  pages = {119--122},
  author = {N. Etemadi},
  title = {An elementary proof of the strong law of large numbers},
  journal = {Zeitschrift fuer Wahrscheinlichkeitstheorie und Verwandte Gebiete}
}

@article{Koskela2018,
  doi = {10.1515/sagmb-2017-0011},
  year  = {2018},
  month = {jun},
  publisher = {Walter de Gruyter {GmbH}},
  volume = {17},
  number = {3},
  author = {Jere Koskela},
  title = {Multi-locus data distinguishes between population growth and multiple merger coalescents},
  journal = {Statistical Applications in Genetics and Molecular Biology}
}

@article{BLS15,
  doi = {10.1214/18-ejp175},
  year  = {2018},
  publisher = {Institute of Mathematical Statistics},
  volume = {23},
  number = {0},
  author = {Matthias Birkner and Huili Liu and Anja Sturm},
  title = {Coalescent results for diploid exchangeable population models},
  journal = {Electronic Journal of Probability}
}

@article{Chvtal1979,
  doi = {10.1016/0012-365x(79)90084-0},
  url = {https://doi.org/10.1016/0012-365x(79)90084-0},
  year = {1979},
  publisher = {Elsevier {BV}},
  volume = {25},
  number = {3},
  pages = {285--287},
  author = {V. Chv{\'{a}}tal},
  title = {The tail of the hypergeometric distribution},
  journal = {Discrete Mathematics}
}

@article{DS05,
  doi = {10.1016/j.spa.2005.04.009},
  year = {2005},
  month = oct,
  publisher = {Elsevier {BV}},
  volume = {115},
  number = {10},
  pages = {1628--1657},
  author = {Rick Durrett and Jason Schweinsberg},
  title = {A coalescent model for the effect of advantageous mutations on the genealogy of a population},
  journal = {Stochastic Processes and their Applications}
}

@article{Eldon2020,
  doi = {10.1146/annurev-genet-021920-095932},
  url = {https://doi.org/10.1146/annurev-genet-021920-095932},
  year = {2020},
  month = nov,
  publisher = {Annual Reviews},
  volume = {54},
  number = {1},
  pages = {213--236},
  author = {Bjarki Eldon},
  title = {Evolutionary Genomics of High Fecundity},
  journal = {Annual Review of Genetics}
}

@article{10.1093/molbev/msaa179,
    author = {Menardo, Fabrizio and Gagneux, Sébastien and Freund, Fabian},
    title = {Multiple Merger Genealogies in Outbreaks of Mycobacterium tuberculosis},
    journal = {Molecular Biology and Evolution},
    volume = {38},
    number = {1},
    pages = {290-306},
    year = {2021},
    month = {01},
    issn = {1537-1719},
    doi = {10.1093/molbev/msaa179},
    url = {https://doi.org/10.1093/molbev/msaa179},
    eprint = {https://academic.oup.com/mbe/article-pdf/38/1/290/35389068/msaa179.pdf},
}

@article{10.1098/rsos.171060,
    author = {Kato, Mamoru and Vasco, Daniel A. and Sugino, Ryuichi and Narushima, Daichi and Krasnitz, Alexander},
    title = {Sweepstake evolution revealed by population-genetic analysis of copy-number alterations in single genomes of breast cancer},
    journal = {Royal Society Open Science},
    volume = {4},
    number = {9},
    pages = {171060},
    year = {2017},
    month = {09},
    issn = {2054-5703},
    doi = {10.1098/rsos.171060},
    url = {https://doi.org/10.1098/rsos.171060},
    eprint = {https://royalsocietypublishing.org/rsos/article-pdf/doi/10.1098/rsos.171060/958679/rsos.171060.pdf},
}

@article {Arnasonsweepstakes2022,
article_type = {journal},
title = {Sweepstakes reproductive success via pervasive and recurrent selective sweeps},
author = {\'Arnason, Einar and Koskela, Jere and Halld\'orsd\'ottir, Katr\'in and Eldon, Bjarki},
editor = {Gagnaire, Pierre-Alexandre and Przeworski, Molly and Freund, Fabian and Gagnaire, Pierre-Alexandre},
volume = 12,
year = 2023,
month = {feb},
pub_date = {2023-02-20},
pages = {e80781},
citation = {eLife 2023;12:e80781},
doi = {10.7554/eLife.80781},
url = {https://doi.org/10.7554/eLife.80781},
journal = {eLife},
issn = {2050-084X},
publisher = {eLife Sciences Publications, Ltd},
}

@Book{ewens79,
  author = 	 {WJ Ewens},
  title = 	 {Mathematical population genetics},
  publisher = 	 {Springer-Verlag},
  year = 	 1979,
  doi={10.1007/978-0-387-21822-9},
  address = 	 {Berlin}}

@article{wakeley2012gene,
	author = {John Wakeley and L\'{e}andra King and Bobbi S Low and Sohini Ramachandran},
	journal = {Genetics},
	number = {4},
	pages = {1433--1445},
	publisher = {Genetics Soc America},
	title = {Gene genealogies within a fixed pedigree, and the robustness of {K}ingman{\rq}s coalescent},
	volume = {190},
	doi={10.1534/genetics.111.135574},
	year = {2012}
}

@article{Minadakis2025,
  title = {Genomic Surveillance and Molecular Evolution of Fungicide Resistance in European Populations of Wheat Powdery Mildew},
  volume = {26},
  ISSN = {1364-3703},
  url = {http://dx.doi.org/10.1111/mpp.70071},
  DOI = {10.1111/mpp.70071},
  number = {3},
  journal = {Molecular Plant Pathology},
  publisher = {Wiley},
  author = {Minadakis,  Nikolaos and Jigisha,  Jigisha and Cornetti,  Luca and Kunz,  Lukas and M\"{u}ller,  Marion C. and Torriani,  Stefano F. F. and Menardo,  Fabrizio},
  year = {2025},
  month = Mar 
}

@article{Jigisha2025,
  title = {Population genomics and molecular epidemiology of wheat powdery mildew in Europe},
  volume = {23},
  ISSN = {1545-7885},
  url = {http://dx.doi.org/10.1371/journal.pbio.3003097},
  DOI = {10.1371/journal.pbio.3003097},
  number = {5},
  journal = {PLOS Biology},
  publisher = {Public Library of Science (PLoS)},
  author = {Jigisha, Jigisha and Ly, Jeanine and Minadakis, Nikolaos
                  and Freund, Fabian and Kunz, Lukas and Piechota,
                  Urszula and Akin, Beyhan and Balmas, Virgilio and
                  Ben-David, Roi and Bencze, Szilvia and Bourras,
                  Salim and Bozzoli, Matteo and Cotuna, Otilia and
                  Couleaud, Gilles and Cséplő, Mónika and Czembor,
                  Paweł and Desiderio, Francesca and D\"{o}rnte, Jost
                  and Dreiseitl, Antonín and Feechan, Angela and
                  Gadaleta, Agata and Gauthier, Kevin and Giancaspro,
                  Angelica and Giove, Stefania L. and
                  Handley-Cornillet, Alain and Hubbard, Amelia and
                  Karaoglanidis, George and Kildea, Steven and Koc,
                  Emrah and Liatukas, Žilvinas and Lopes, Marta S. and
                  Mascher, Fabio and McCabe, Cathal and Miedaner,
                  Thomas and Martínez-Moreno, Fernando and Nellist,
                  Charlotte F. and Okoń, Sylwia and Praz, Coraline and
                  Sánchez-Martín, Javier and Sărăţeanu, Veronica and
                  Schulz, Philipp and Schwartz, Nathalie and Seghetta,
                  Daniele and Martel, Ignacio Solís and Švarta, Agrita
                  and Testempasis, Stefanos and Villegas, Dolors and
                  Widrig, Victoria and Menardo, Fabrizio},
  editor = {Kamoun,  Sophien},
  year = {2025},
  month = May,
  pages = {e3003097}
}

@article{HM11b,
	author = {T Huillet and M M\"{o}hle},
	journal = {Stoch Models},
	pages = {521--554},
	title = {Population genetics models with skewed fertilities: forward and backward analysis},
	volume = 27,
	doi={https://doi.org/10.1080/15326349.2011.593411},
	year = 2011
}

@Manual{Rsystem,
    title = {R: A Language and Environment for Statistical Computing},
    author = {{R Core Team}},
    organization = {R Foundation for Statistical Computing},
    address = {Vienna, Austria},
    year = {2026},
    doi={10.32614/R.manuals},
    url = {https://www.R-project.org/},
  }

@book{feng2010poisson,
  title={The {Poisson-Dirichlet} distribution and related topics: models and asymptotic behaviors},
  author={Feng, Shui},
  year={2010},
  doi={10.1007/978-3-642-11194-5},
  publisher={Springer Science \& Business Media}
}

@article{M98,
	author = {Martin M{\"{o}}hle},
	journal = {Adv Appl Prob},
	pages = {493--512},
	title = {A convergence theorem for {Markov} chains arising in population genetics and the coalescent with selfing},
	volume = 30,
	doi={https://doi.org/10.1239/aap/1035228080},
	year = 1998
}

@article{JFB:JFB13143,
author = {Waples, R. S.},
title = {Tiny estimates of the ${N_e}/{N}$ ratio in marine fishes: Are they real?},
journal = {Journal of Fish Biology},
volume = {89},
number = {6},
publisher = {Blackwell Publishing Ltd},
issn = {1095-8649},
doi = {10.1111/jfb.13143},
pages = {2479--2504},
year = {2016}
}

@article{schweinsberg03,
	author = {J Schweinsberg},
	journal = {Stoch Proc Appl},
	pages = {107--139},
	title = {Coalescent processes obtained from supercritical {G}alton-{W}atson processes},
	volume = 106,
	doi={10.1016/S0304-4149(03)00028-0},
	year = 2003
}

@article{S00,
	author = {J Schweinsberg},
	journal = {Electron J Probab},
	pages = {1--50},
	title = {Coalescents with simultaneous multiple collisions},
	volume = 5,
	doi={10.1214/EJP.v5-68},
	year = 2000
}

@article{MS03,
	author = {M M{\"{o}}hle and S Sagitov},
	journal = {J Math Biol},
	pages = {337--352},
	title = {Coalescent patterns in diploid exchangeable population models},
	volume = 47,
	doi={10.1007/s00285-003-0218-6},
	year = 2003
}

@article{S03,
	author = {S Sagitov},
	journal = {J Appl Probab},
	pages = {839--854},
	title = {Convergence to the coalescent with simultaneous multiple mergers},
	volume = 40,
	doi={10.1239/jap/1067436085},
	year = 2003
}

@article{Li1998,
	author = {Gang Li and Dennis Hedgecock},
	doi = {10.1139/f97-312},
	journal = {Can. J. Fish. Aquat. Sci.},
	month = {apr},
	number = {4},
	pages = {1025--1033},
	publisher = {Canadian Science Publishing},
	title = {Genetic heterogeneity, detected by {PCR}-{SSCP}, among samples of larval {P}acific oysters ( \emph{Crassostrea gigas} ) supports the hypothesis of large variance in reproductive success },
	volume = {55},
	year = {1998}
}

@article{EBBF2015,
	author = {B Eldon and M Birkner and J Blath and F Freund},
	journal = {Genetics},
	pages = {841--856},
	title = {Can the site-frequency spectrum distinguish exponential population growth from multiple-merger coalescents?},
	volume = 199,
	doi={https://doi.org/10.1534/genetics.114.173807},
	year = 2015
}

@article{K82,
title = {The coalescent},
journal = {Stochastic Processes and their Applications},
volume = {13},
number = {3},
pages = {235-248},
year = {1982},
issn = {0304-4149},
doi = {https://doi.org/10.1016/0304-4149(82)90011-4},
url = {https://www.sciencedirect.com/science/article/pii/0304414982900114},
author = {J.F.C. Kingman}
}

@article{K82b,
	author = {J F C Kingman},
	journal = {J App Probab},
	pages = {27--43},
	title = {On the genealogy of large populations},
	volume = {{19A}},
	doi={10.2307/3213548},
	year = 1982
}

@article{H83b,
	author = {R R Hudson},
	journal = {Theor Popul Biol},
	pages = {183--201},
	title = {Properties of a neutral allele model with intragenic recombination},
	volume = 23,
	doi={https://doi.org/10.1016/0040-5809(83)90013-8},
	year = 1983
}

@article{T83,
	author = {F Tajima},
	journal = {Genetics},
	pages = {437--460},
	title = {Evolutionary relationships of {DNA} sequences in finite populations},
	volume = 105,
	doi={https://doi.org/10.1093/genetics/105.2.437},
	year = 1983
}

@article{Birkner2024,
  title = {The joint fluctuations of the lengths of the {Beta}({2 -  $\alpha$}, {$\alpha$})-coalescents},
  volume = {34},
  ISSN = {1050-5164},
  url = {http://dx.doi.org/10.1214/23-AAP1964},
  DOI = {10.1214/23-aap1964},
  number = {1A},
  journal = {The Annals of Applied Probability},
  publisher = {Institute of Mathematical Statistics},
  author = {Birkner,  Matthias and Dahmer,  Iulia and Diehl,  Christina S. and Kersting,  G\"{o}tz},
  year = {2024},
  month = feb 
}

@article{DKW2014,
	author = {I Dahmer and G Kersting and A Wakolbinger},
	journal = {Comb Prob Comp},
	pages = {1010--1027},
	title = {The total external length of {Beta}-coalescents},
	volume = 23,
	doi={https://doi.org/10.1017/S0963548314000297},
	year = 2014
}

@article{W1975,
	author = {G A Watterson},
	journal = {Theor Pop Biol},
	pages = {256--276},
	title = {On the number of segregating sites in genetical models without recombination},
	volume = 7,
	doi={10.1016/0040-5809(75)90020-9},
	year = 1975
}

@inproceedings{H94,
	address = {London},
	author = {D Hedgecock},
	booktitle = {Genetics and evolution of Aquatic Organisms},
	editor = {A Beaumont},
	pages = {1222--1344},
	publisher = {Chapman and Hall},
	title = {Does variance in reproductive success limit effective population sizes of marine organisms?},
	isbn = {978-0-412-49370-6},
	year = 1994
}

@article{Blath2016,
	author = {Jochen Blath and Mathias Christensen Cronj\"{a}ger and Bjarki Eldon and Matthias Hammer},
	doi = {10.1016/j.tpb.2016.04.002},
	journal = {Theoretical Population Biology},
	month = {aug},
	pages = {36--50},
	publisher = {Elsevier {BV}},
	title = {The site-frequency spectrum associated with {$\Xi$}-coalescents},
	url = {http://dx.doi.org/10.1016/j.tpb.2016.04.002},
	volume = {110},
	year = {2016}
}

@article{HP11,
	author = {D Hedgecock and A I Pudovkin},
	journal = {Bull Marine Science},
	pages = {971--1002},
	title = {Sweepstakes reproductive success in highly fecund marine fish and shellfish: a review and commentary},
	volume = 87,
	doi={10.5343/bms.2010.1051},
	year = 2011
}

@incollection{H82,
	address = {New York},
	author = {D Hedgecock and M Tracey and K Nelson},
	booktitle = {The {B}iology of {C}rustacea: {E}mbryology, {M}orphology, and {G}enetics},
	editor = {L G Abele},
	pages = {297--403},
	publisher = {Academic Press},
	title = {Genetics},
	volume = 2,
	isbn={0121064050},
	year = 1982
}

@incollection{B94,
	address = {New York},
	author = {A T Beckenbach},
	booktitle = {Non-neutral {E}volution},
	editor = {B Golding},
	pages = {188--198},
	publisher = {Chapman \& Hall},
	title = {Mitochondrial haplotype frequencies in oysters: neutral alternatives to selection models},
	doi={10.1007/978-1-4615-2383-3},
	year = 1994
}

@article{Vendrami2021,
  title = {Sweepstake reproductive success and collective dispersal produce chaotic genetic patchiness in a broadcast spawner},
  volume = {7},
  ISSN = {2375-2548},
  url = {http://dx.doi.org/10.1126/sciadv.abj4713},
  DOI = {10.1126/sciadv.abj4713},
  number = {37},
  journal = {Science Advances},
  publisher = {American Association for the Advancement of Science (AAAS)},
  author = {Vendrami,  David L. J. and Peck,  Lloyd S. and Clark,  Melody S. and Eldon,  Bjarki and Meredith,  Michael and Hoffman,  Joseph I.},
  year = {2021},
  month = sep 
}

@article{AH2015,
	author = {Einar \'{A}rnason and Katr\'{\i}n Halld\'{o}rsd\'{o}ttir},
	doi = {10.7717/peerj.786},
	journal = {{PeerJ}},
	pages = {e786},
	publisher = {{PeerJ}},
	title = {Nucleotide variation and balancing selection at the \emph{Ckma} gene in {A}tlantic cod: analysis with multiple merger coalescent models },
	volume = {3},
	year = {2015}
}

@article{A04,
	author = {E \'{A}rnason},
	journal = {Genetics},
	pages = {1871--1885},
	title = {Mitochondrial cytochrome \emph{b} variation in the high-fecundity {A}tlantic cod: trans-{A}tlantic clines and shallow gene genealogy},
	doi={https://doi.org/10.1534/genetics.166.4.1871},
	volume = 166,
	year = 2004
}

@article{EW06,
	author = {B Eldon and J Wakeley},
	journal = {Genetics},
	pages = {2621--2633},
	title = {Coalescent processes when the distribution of offspring number among individuals is highly skewed},
	volume = 172,
	doi={https://doi.org/10.1534/genetics.105.052175},
	year = 2006
}

@book{W07,
  author = {John Wakeley},
  title = {Coalescent Theory: An Introduction},
  year = {2009},
  publisher = {Roberts \& Company Publishers},
  address = {Greenwood Village},
  isbn = {0-9747077-5-9},
  language = {English},
}

@article{SW08,
	author = {O Sargsyan and J Wakeley},
	journal = {Theor Pop Biol},
	pages = {104--114},
	title = {A coalescent process with simultaneous multiple mergers for approximating the gene genealogies of many marine organisms},
	volume = 74,
	doi={10.1016/j.tpb.2008.04.009},
	year = 2008
}

@article{HM12,
	author = {T Huillet and M M\"{o}hle},
	journal = {Theor Popul Biol},
	pages = {5--14},
	title = {On the extended {M}oran model and its relation to coalescents with multiple collisions},
	volume = 87,
	doi={https://doi.org/10.1016/j.tpb.2011.09.004},
	year = 2013
}

@incollection{BB09,
	author = {M Birkner and J Blath},
	booktitle = {Trends in stochastic analysis},
	editor = {J Blath and P M\"{o}rters and M Scheutzow},
	pages = {329--363},
	publisher = {Cambridge University Press},
	title = {Measure-valued diffusions, general coalescents and population genetic inference},
	doi={https://doi.org/10.1017/CBO9781139107020},
	year = 2009
}

@article{Kelleher2016,
	author = {Jerome Kelleher and Alison M Etheridge and Gilean McVean},
	doi = {10.1371/journal.pcbi.1004842},
	editor = {Yun S. Song},
	journal = {{PLOS} Computational Biology},
	month = {may},
	number = {5},
	pages = {e1004842},
	publisher = {Public Library of Science ({PLoS})},
	title = {Efficient Coalescent Simulation and Genealogical Analysis for Large Sample Sizes},
	url = {http://dx.doi.org/10.1371/journal.pcbi.1004842},
	volume = {12},
	year = {2016}
}

@article{BBE13,
	author = {M Birkner and J Blath and B Eldon},
	journal = {Genetics},
	pages = {255--290},
	title = {An ancestral recombination graph for diploid populations with skewed offspring distribution},
	volume = 193,
	doi={https://doi.org/10.1534/genetics.112.144329},
	year = 2013
}

@article{Eldon2026biorxiv,
  title = {Beta-coalescents when sample size is large},
  url = {http://dx.doi.org/10.64898/2025.12.30.697022},
  DOI = {10.64898/2025.12.30.697022},
  journal={bioRxiv},
  publisher = {openRxiv},
  author = {Chetwynd-Diggle,  Jonathan A and Eldon, Bjarki},
  year = {2026},
  month = jan 
}

@article{BBE2013a,
	author = {M Birkner and J Blath and B Eldon},
	journal = {Genetics},
	pages = {1037--1053},
	title = {Statistical properties of the site-frequency spectrum associated with {$\Lambda$}-coalescents},
	volume = 195,
	doi={https://doi.org/10.1534/genetics.113.156612},
	year = 2013
}

@article{WT2003,
	author = {J Wakeley and T Takahashi},
	journal = {Mol Biol Evol},
	pages = {208--2013},
	title = {Gene genealogies when the sample size exceeds the effective size of the population},
	volume = 20,
	doi={10.1093/molbev/msg024},
	year = 2003
}

@article{BBMST09,
	author = {M Birkner and J Blath and M M\"{o}hle and M Steinr\"{u}cken and J Tams},
	journal = {ALEA Lat. Am. J. Probab. Math. Stat.},
	pages = {25--61},
	title = {A modified lookdown construction for the {X}i-{F}leming-{V}iot process with mutation and populations with recurrent bottlenecks},
	volume = 6,
	doi={https://doi.org/10.34657/1856},
	year = 2009
}

@article{BCS2014,
	author = {A Bhaskar and AG Clark and YS Song},
	journal = {PNAS},
	pages = {2385--2390},
	title = {Distortion of genealogical properties when the sample size is very large},
	volume = 111,
	doi={10.1073/pnas.1322709111},
	year = 2014
}

@article{SD2005,
	author = {J Schweinsberg and R Durrett},
	journal = {Ann Appl Probab},
	title = {Random partitions approximating the coalescence of lineages during a selective sweep},
	volume = {1591--1651},
	doi={10.1214/105051605000000430},
	year = 2005
}

@article{DK99,
	author = {P Donnelly and T G Kurtz},
	journal = {Ann Probab},
	pages = {166--205},
	title = {Particle Representations for Measure-Valued Population Models},
	volume = 27,
	doi={https://doi.org/10.1214/aop/1022677258},
	year = 1999
}

@article{P99,
	author = {J Pitman},
	journal = {Ann Probab},
	pages = {1870--1902},
	title = {Coalescents with multiple collisions},
	volume = 27,
	doi={10.1214/aop/1022874819},
	year = 1999
}

@article{S99,
	author = {S Sagitov},
	journal = {J Appl Probab},
	pages = {1116--1125},
	title = {The general coalescent with asynchronous mergers of ancestral lines},
	volume = 36,
	doi={10.1239/jap/1032374759},
	year = 1999
}

@article{MS01,
	author = {M M\"{o}hle and S Sagitov},
	journal = {Ann Probab},
	pages = {1547--1562},
	title = {A classification of coalescent processes for haploid exchangeable population models},
	volume = 29,
	doi={https://doi.org/10.1214/aop/1015345761},
	year = 2001
}

@article{F95,
	author = {{Y-X} Fu},
	journal = {Theor Popul Biol},
	pages = {172--197},
	title = {Statistical Properties of Segregating Sites},
	volume = 48,
	doi={https://doi.org/10.1006/tpbi.1995.1025},
	year = 1995
}

@article{BBC05,
	author = {M Birkner and J Blath and M Capaldo and A M Etheridge and M M\"{o}hle and J Schweinsberg and A Wakolbinger},
	journal = {Electron. J. Probab},
	pages = {303--325},
	title = {Alpha-stable branching and {Beta}-coalescents},
	volume = 10,
	doi={https://doi.org/10.1214/EJP.v10-241},
	year = 2005
}

@article{B09,
	author = {N Berestycki},
	journal = {Ensaios Math\'{e}maticos},
	pages = {1--193},
	title = {Recent progress in Coalescent Theory},
	volume = 16,
	doi={10.21711/217504322009/em161},
	year = 2009
}

@article{wright31:_evolut_mendel,
    author = {Wright, Sewall},
    title = {EVOLUTION IN MENDELIAN POPULATIONS},
    journal = {Genetics},
    volume = {16},
    number = {2},
    pages = {97-159},
    year = {1931},
    month = {03},
    issn = {1943-2631},
    doi = {10.1093/genetics/16.2.97},
    url = {https://doi.org/10.1093/genetics/16.2.97},
    eprint = {https://academic.oup.com/genetics/article-pdf/16/2/97/35081059/genetics0097.pdf},
}

@Article{fisher22,
  author = 	 {RA Fisher},
  title = 	 {On the dominance ratio},
  journal = 	 {Proc Royal Society Edinburgh},
  year = 	 1923,
  volume = 	 42,
  doi={10.1017/S0370164600023993},
  pages = 	 {321--341}}

@PhdThesis{mccloskey65,
  author = 	 {JW McCloskey},
  title = 	 {A model for the distribution of individuals by species in an environment},
  school = 	 {Michigan State University},
  url =          {https://doi.org/doi:10.25335/3syn-6j67},
  year = 	 1965}

@Book{engen78:_abund,
  author = 	 {S Engen},
  title = 	 {Stochastic abundance models: with emphasis on biological communities and species diversity},
  publisher = 	 {Chapman-Hall},
  year = 	 1978,
  doi={https://doi.org/10.1007/978-94-009-5784-8},
  isbn={94-009-5784-X},
  address = 	 {London}}

@article{K1969,
	author = {M Kimura},
	journal = {Genetics},
	pages = {893--903},
	title = {The number of heterozygous nucleotide sites maintained in a finite population due to steady flux of mutations},
	volume = 61,
	doi={https://doi.org/10.1093/genetics/61.4.893},
	year = 1969
}

@article{BBS08,
	author = {J Berestycki and N Berestycki and J Schweinsberg},
	journal = {Ann Inst H Poincar\'{e} Probab Statist},
	pages = {214--238},
	title = {Small-time behavior of {beta} coalescents},
	volume = 44,
	doi={https://doi.org/10.1214/07-AIHP103},
	year = 2008
}

@article{Niwa2016,
	author = {Hiro-Sato Niwa and Kazuya Nashida and Takashi Yanagimoto},
	doi = {10.1093/icesjms/fsw070},
	journal = {{ICES} Journal of Marine Science: Journal du Conseil},
	month = {may},
	number = {9},
	pages = {2181--2189},
	publisher = {Oxford University Press ({OUP})},
	title = {Reproductive skew in Japanese sardine inferred from {DNA} sequences},
	url = {http://dx.doi.org/10.1093/icesjms/fsw070},
	volume = {73},
	year = {2016}
}

@article{BBS07,
	author = {J Berestycki and N Berestycki and J Schweinsberg},
	journal = {Ann Probab},
	pages = {1835--1887},
	title = {Beta-coalescents and continuous stable random trees},
	doi={https://doi.org/10.1214/009117906000001114},
	volume = 35,
	year = 2007
}

@article{hedgecock_etal07_marinebiology,
	author = {D Hedgecock and S Launey and A I Pudovkin and Y Naciri and S Lapegue and F Bonhomme},
	journal = {Marine Biology},
	pages = {1173--1182},
	title = {Small effective number of parents ({$N_b$}) inferred for a naturally spawned cohort of juvenile {European} flat oysters \textit{{Ostrea} edulis}},
	type = {Journal Article},
	volume = {150},
	doi={10.1007/s00227-006-0441-y},
	url={http://dx.doi.org/10.1007/s00227-006-0441-y},
	year = {2007}
}

@article{mohle2000total,
  title={Total variation distances and rates of convergence for ancestral coalescent processes in exchangeable population models},
  author={M{\"o}hle, Martin},
  journal={Advances in Applied Probability},
  pages={983--993},
  year={2000},
  doi={10.1017/S0001867800010417},
  publisher={JSTOR}
}

@article{mohle1998robustness,
  title={Robustness results for the coalescent},
  author={M{\"o}hle, Martin},
  journal={Journal of Applied Probability},
  volume={35},
  number={2},
  pages={438--447},
  year={1998},
  doi={10.1017/S0021900200015060},
  publisher={Cambridge University Press}
}

\end{document}